\documentclass[oneside,10pt]{article}

\usepackage[b5paper]{geometry}
\usepackage{amsfonts,amsmath,latexsym,amssymb}
\usepackage{theorem}
\usepackage{mathrsfs,upref} 
\usepackage{mathptmx}
\usepackage{mathrsfs}
\usepackage{wasysym}
\usepackage{graphicx}

\allowdisplaybreaks

\newtheorem{theorem}{Theorem}

\theoremstyle{definition}

\newtheorem{example}{Example}

\newtheorem{remark}{Remark}

\title{Families of self-inverse functions and dilogarithm identities}
\author{{\bf Lauri Alha} \\
email: lare1962@gmail.com.}

\begin{document}


\maketitle
\begin{abstract}
 We introduce a self-inverse function via an integral equivalent to a two-term combination of dilogarithms. We refer to this function as a 
 \emph{fundamental form}, since there is a family of extensions of this function that satisfy similar self-inverse and symmetric properties. 
 We also construct a family of functions generalizing the fundamental form via two auxiliary parameters, which we refer to as \emph{shape} 
 and \emph{scale factors}. Through new integration techniques, we introduce and prove a variety of dilogarithm identities and evaluations 
 for dilogarithm ladders and for two-term dilogarithm combinations. The functions $\gemini_{a}^{b}(x)$ we introduce are referred to as 
 \emph{gemini functions} and may be seen as providing a broad framework in the derivation of and application of dilogarithm identities.

\vspace{0.1in}

\end{abstract}

\section{Introduction}
 The polylogarithm function $\operatorname{Li}_{m}(z) = \sum_{n=1}^{\infty} \frac{z^n}{n^m}$ provides one of the most important 
 generalizations of the Riemann zeta function $\zeta(m) = \sum_{n=1}^{\infty} \frac{1}{n^m}$ from both number-theoretic perpsectives 
 and in terms of the uses of polylogarithms within special functions theory and related areas of computational mathematics and 
 mathematical analysis. The special function $\operatorname{Li}_{2}(z)$ is referred to as the \emph{dilogarithm} and may be seen as the 
 most basic instance of a higher logarithm function. Recent research related to the behavior of and structural properties associated with 
 the dilogarithm is representative of the variety of different disciplines in which the $\operatorname{Li}_{2}$ function arises. In this 
 direction, we highlight the work of Jaipong et al.\ \cite{JaipongLangTanTee2023} on the evaluation of series involving 
 $\operatorname{Li}_{2}$ using geometric properties associated with hyperbolic cylinders with holonomy in $\text{PSL}(2, \mathbb{R})$, the 
 connection discovered by Bridgeman~\cite{Bridgeman2021} between $\operatorname{Li}_{2}$ and the solutions to Pell's equation 
 $x^2 - ny^2 = \pm 1$, the derivations of dilogarithm identities due to Nakanishi \cite{Nakanishi2024} using paths in scattering diagrams 
 associated with cluster algebras, and the connections established by Freed and Neitzke \cite{FreedNeitzke2023} between the dilogarithm 
 and abelian Chern--Simons theory. The introduction of and application of dilogarithm relations, together with techniques that may be 
 used to generate new dilogarithm relations, are motivated by the many areas inside and outside of mathematics in which the 
 $\operatorname{Li}_{2}$ function plays a central role. 

\begin{figure}[htbp!]
\begin{center}
\item
\includegraphics[width=9cm, height=6cm]{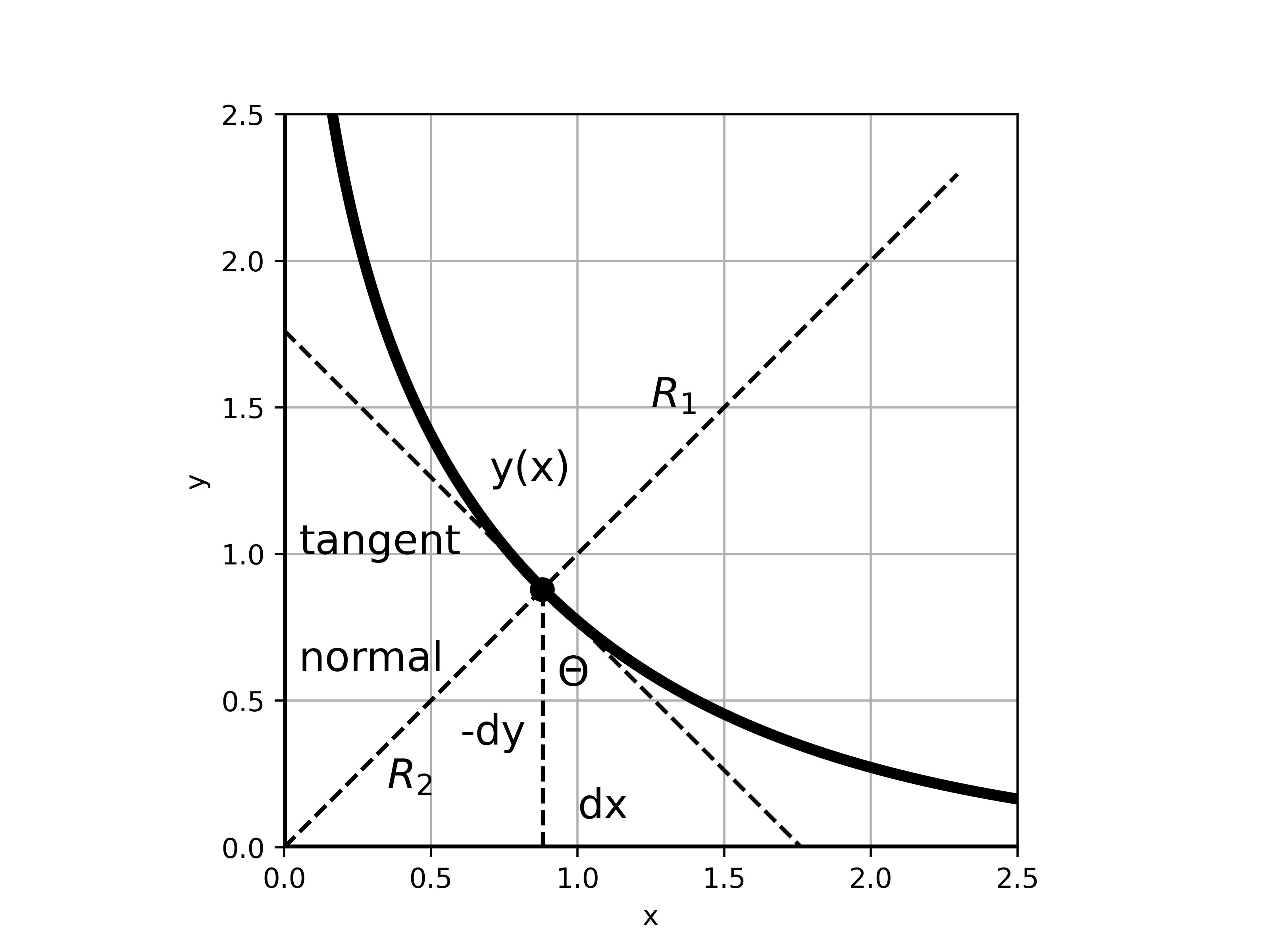}
\end{center}
\caption{\label{Figure1} A self-inverse function derived using \eqref{splitimply}.}
\end{figure}

 Our techniques for generating dilogarithm identities rely on a family of functions that we refer to as \emph{gemini functions}. To begin 
 with, we construct a function that may be thought of as being based on the parallel postulate in the Euclidean plane, as clarified below. 
 The $y$-coordinates for the graph of this function correspond to values of the form $d$ such that
\begin{equation}\label{parallelpost}
 \Pi(d) = 2\arctan\big(e^{-d}\big), 
\end{equation}
 with $d$ giving the vertical distance between the function and the $x$-axis, and where the tangential angle $\theta$ of this function 
 always corresponds to the parallel angle $\Pi(d)$. This formula may be seen as illustrating the concept of the \emph{angle of parallelism} 
that plays an important role in hyperbolic geometry. An elementary description of this topic can be found in Anderson's monograph on 
hyperbolic geometry \cite{And_05}. An illustration related to a derivation of \eqref{parallelpost} is shown in Figure \ref{Figure1}. The 
function referenced above can be evaluated by solving a simple separable first-order differential equation, in the following manner. 
According to Figure \ref{Figure1}, we can write $\theta=2\arctan{(e^{-y})}=\arctan{\big(-\frac{dx}{dy}\big)}$, and this gives us that 

\begin{equation}\label{splitimply}
 \tan\big[2\arctan(e^{-y})\big] = -\frac{dx}{dy}. 
\end{equation}

We thus obtain that $ x + C = \ln{\big(\frac{1+e^{y}}{1-e^{y}}\big)}$ through the application of an appropriate integral operator to 
\eqref{splitimply}, i.e., so that $ y=\ln{\big(\frac{1+e^{x+C}}{1-e^{x+C}}\big)}$, and we proceed to set $C = 0$. Informally, the assigning 
of a value to $C$ can be thought of as producing a shifting along the $x$-axis. We thus obtain a function that is symmetric about $y = x$ 
in the first quadrant, as illustrated in Figure \ref{Figure1}. This function is a self-inverse function, which enables us to derive the 
 equivalent relations 
\begin{equation}\label{xybiconditional}
 x=\ln{\left(\frac{1+e^{-y}}{1-e^{-y}}\right)} \iff y=\ln\left(\frac{1+e^{-x}}{1-e^{-x}}\right)=\ln\left(\frac{e^{x}+1}{e^{x}-1}\right)
\end{equation}
 via an interchange of the $x$-and $y$-coordinates. 
 
 We refer  either formulation
 of the two equivalent relations in \eqref{xybiconditional}
 as a \emph{fundamental form}, 
     and this can be seen as serving as a foundation for the 
  gemini functions  introduced in this paper. 
   The function of $x$ 
 on the right-hand side of the biconditional equivalence in \eqref{xybiconditional} has 
 also been considered in the works by Romakina \cite{Rom_18} and Basmajian \cite{Bas_93}, from hyperbolic geometry-based 
 perspectives. Observe the trivial factor given in exponent term in the denominator. The purpose of this factor will be introduced in the 
 upcoming section. This value 1 is also the subscript value for the gemini sign, denoting the fundamental form of a gemini function, with 
\begin{multline*}
\gemini_{1}{(x)} =\ln{\left(\frac{1+1 \cdot e^{-x}}{1-e^{-x}}\right)}=\ln\left[\coth\left(\frac{x}{2}\right)\right]=
2\operatorname{arctanh}{(e^{-x})} = \\ 
\operatorname{arcsinh}{\left[\frac{1}{\sinh(x)}\right]}=\ln\left[\frac{1+\cosh(x)}{\sinh(x)}\right]. 
\end{multline*}
 
 The symmetric property of gemini functions given by such functions being self-inverse may be considered in relation how integrals of 
 gemini functions consist of two dilogarithm terms, excluding the cases whereby the shape factor is $-1$ or $0$.
This feature plays the key role in our study in this paper. Equivalent definitions for the dilogarithm function $\text{Li}_{2}$
are below given in \eqref{Li2definition}, with
\begin{equation}\label{Li2definition}
\operatorname{Li_{2}}(x)=\sum_{k=1}^{\infty}\frac{x^{k}}{k^{2}}=-\int_{0}^{x}\frac{\ln(1-t)dt}{t}
\end{equation}
for arguments $x$ such that $|x| < 1$, whereas the \emph{polylogarithm} $\text{Li}_{s+1}$ is such that 
$$ \text{Li}_{s+1}(x) = \sum_{k=1}^{\infty} \frac{x^k}{k^{s+1}} = \frac{x (-1)^s}{s!} \int_{0}^{1} \frac{\ln^s(t)}{1 - t x} \, dt $$ 
for positive integers $s$ and $x$ as above. For background on the dilogarithm and polylogarithm functions and the importance of
these special functions within many different areas of mathematics, we refer to monographs by Zagier \cite{Zag_07} and by Lewin 
\cite{Lew_58,Lewin1981,Lewin1991Structural}. 

The integral of the fundamental form is shown in \eqref{integraltwoterm}, with
\begin{equation}\label{integraltwoterm}
\int \gemini_{1}{(x)} \, dx=\int \ln\left(\frac{1+e^{-x}}{1-e^{-x}}\right) \, 
 dx =\operatorname{Li_{2}}(-e^{-x}) - 
\operatorname{Li_{2}}(e^{-x})+C. 
\end{equation}
The total area $A_{tot}$ bounded by the fundamental form and the positive coordinate axes is finite and it is given by
$$ A_{tot} = \int_{0}^{\infty} \gemini_{1}{(x)} \, dx=\frac{\pi^2}{4}. $$

\subsection{Motivating results}\label{subsectionmotivating}
 There is a long and rich history concerning the applications of closed-form evaluations for two-term combinations of dilogarithmic values 
 with algebraic argument. There are notable such applications within conformal field theory \cite{Bytsko1999,Byt01_99}, and remarkable 
 and topology-based derivations of two-term dilogarithm evaluations are due to Khoi \cite{Kho14}, who used a knot-theoretic argument 
 involving Seifert surfaces to prove that
\begin{equation}\label{displayKhoimotivate}
L\left( \frac{1}{\phi(\phi + \sqrt{\phi})} \right) - L\left( \frac{\phi}{\phi + \sqrt{\phi}} \right) = \frac{\pi^2}{20}, 
\end{equation}
 for the golden ratio $\phi = \frac{1 + \sqrt{5}}{2}$ and for the Rogers dilogarithm function $L(z) := \text{Li}_{2}(z) + \frac{1}{2} \ln z \ln(1 - 
 z)$ \cite{Rogers1907}. Through the use of the family of gemini functions introduced in this paper, we introduced closed-form evaluations 
 reminiscent of the Khoi formula in \eqref{displayKhoimotivate}, such as the evaluation 
\begin{multline*}
\text{Li}_2\left(\frac{\sqrt{\phi }+1}{\phi ^2}\right)+\text{Li}_2\left(\phi ^3-\phi ^{5/2}\right) = \frac{17 \pi ^2}{60} 
-\ln \left(\phi ^{5/2}-2 \phi \right) \ln \left(\phi ^{7/2}-\phi ^3\right) - \\ \frac{1}{2} \ln \left(\phi ^{7/2}-\phi ^3\right)
\ln \left(\phi^{11/2}+2 \phi ^{7/2}+\phi ^6+\phi ^4\right)-\frac{11}{8} \ln ^2(\phi ) 
\end{multline*}
 that we prove in Section \ref{subsectionphimotivate}, along with the evaluation 
\begin{multline*}
 \operatorname{Li_{2}}\left(\frac{\sqrt{\phi^{7} + 
 3\phi^{5}}-\phi^{2}-3\phi}{2}\right)-\operatorname{Li_{2}}\left(\frac{1+\sqrt{4\phi-3}}{2\phi}\right) = \\ 
 -\frac{\pi^{2}}{10}+\ln^{2}(\phi)-\ln\left(\frac{\phi^{2}+1+\sqrt{4\phi^{3}-
3\phi^{2}}}{2}\right)\ln\left(\frac{2\phi-1+\sqrt{4\phi-3}}{2\phi}\right)
\end{multline*}
 that we prove in Section \ref{secondphimotivate}, and together with the evaluation 
\begin{multline*}
 \operatorname{Li_{2}}\left(\frac{1}{2}\phi-\frac{1}{2}\sqrt{\frac{\phi^{2}+2}{\phi^{3}}}\right)-
\operatorname{Li_{2}}\left(\frac{1}{2\phi^{2}}-\frac{1}{2}\sqrt{\frac{\phi^{2}+2}{\phi^{3}}}\right) = \\ 
\frac{\pi^{2}}{10}+\ln^{2}(\phi)+2\ln(\phi)\ln\left(-\frac{1}{2\phi}+\frac{1}{2}\sqrt{\frac{\phi^{2}+2}{\phi}}\right) 
\end{multline*}
 that we prove in Section \ref{negativephisec}, and together with the evaluation 
\begin{multline*}
 \operatorname{Li_{2}}\left(\frac{1}{2}\sqrt{\phi^{2}+3\phi}-\frac{\phi^{2}+1}{2\phi}\right)-\operatorname{Li_{2}}\left(\frac{1}{2\phi^{2}}
 +\frac{1}{2}\sqrt{\frac{\phi^{2}+2}{\phi^{3}}}\right)=\\
 -\frac{\pi^{2}}{15}+\ln^{2}(\phi)-\ln\left(\frac{\phi^{2}}{2}+\frac{1}{2}\sqrt{\frac{\phi^{2}+2}{\phi}}\right)\ln\left(\frac{1}{2}\phi
 +\frac{1}{2}\sqrt{\frac{\phi^{2}+2}{\phi^{3}}}\right) 
\end{multline*}
that we prove in Section \ref{motivatefourth}, 
 and together with the evaluation 
\begin{multline*}
 \operatorname{Li}_2\left(1-\frac{1+\sqrt{1+\frac{4}{\phi }}}{2 \phi }\right)-\operatorname{Li}_2\left(\frac{1-\sqrt{1+\frac{4}{\phi }}}{2 \phi }\right) = \\ 
 -\frac{1}{2} \log ^2(\phi )+\log \left(\frac{2}{\sqrt{\frac{4}{\phi }+1}+1}\right) \log (\phi )+\frac{\pi ^2}{15}, 
\end{multline*}
 that we prove in Section \ref{framework}. 
 We show how these closed-form evaluations are closely related to the Khoi formula in 
\eqref{displayKhoimotivate}, and this gives us a way of obtaining families of generalizations of Khoi's formula. 

 The concept of a dilogarithm \emph{ladder} was pioneered by Lewin, with reference to Lewin's monograph on structural properties of
polylogarithms \cite{Lewin1991Structural}. The many research works by Lewin et al.\ 
 \cite{AbouzahraLewin1985,AbouzahraLewin1986,AbouzahraLewinXiao1987,CohenLewinZagier1992,Lewin1982,Lewin1993,Lewin1984,Lewin1986,Lewin1991collection} concerning polylogarithm ladders paved the way toward deep results due to Zagier linking polylogarithm ladders and higher Bloch groups
\cite{Zag_07}. Formally, a \emph{ladder}, in the context of higher logarithm functions, of weight $n$ and index $N$ and base $u$ may be defined so that
\begin{equation}\label{ladderdefinition}
L_{n}(N, u) := \frac{\text{Li}_{n}(u^{N})}{N^{n-1}} - \left( 
\sum_{r=0}^{N-1} \frac{A_{r} \text{Li}_{n}(u^r)}{r^{n-1}} + \frac{A_{0} \ln^{n}(u)}{n!} \right) 
\end{equation}
 for $A_{r} \in \mathbb{Q}$ and for an algebraic number $u$ and where the right-hand side of \eqref{ladderdefinition} may be 
 expressed as $$ \sum_{m=2}^{n} \frac{D_{m} \zeta(m) \ln^{n-m}(u)}{(n-m)!} $$ for $D_{m} \in \mathbb{Q}$ \cite[p.\ 
 6]{Lewin1991Structural}. There is a rich history surrounding the number-theoretic uses of polylogarithm ladders, and an especially 
 notable number-theoretic result on polylogarithm ladders is given by the sixteenth-order ladder introduced by Cohen, Lewin, and 
 Zagier \cite{CohenLewinZagier1992}. This motivates the development of techniques in the construction of valid ladder relations, and our 
 gemini function-based techniques introduced in this paper are versatile in terms of generating and proving valid ladder relations. In this 
 regard, letting $\theta$ denote the second smallest Pisot number, i.e., the positive root of the quartic trinomial $x^{4}-x^{3}-1=0$, we 
 introduce the ladder relation
\begin{multline*}
\operatorname{Li}_2\left(\frac{1}{\theta ^{14}}\right)-2 \operatorname{Li}_2\left(\frac{1}{\theta ^7}\right)-2 \operatorname{Li}_2\left(\frac{1}{\theta^5}\right)+\operatorname{Li}_2\left(\frac{1}{\theta ^4}\right) + \\
2 \operatorname{Li}_2\left(\frac{1}{\theta ^2}\right) + 
2 \operatorname{Li}_2\left(\frac{1}{\theta}\right) = 
\frac{\pi ^2}{3}-5 \log ^2\left(\frac{1}{\theta }\right), 
\end{multline*}
as in Theorem \ref{Pisotladdertheorem} below, together with many further ladder relations. 

\section{Generalized gemini functions}
\begin{figure}[htbp!]
\begin{center}
 \includegraphics[width=9cm, height=6cm]{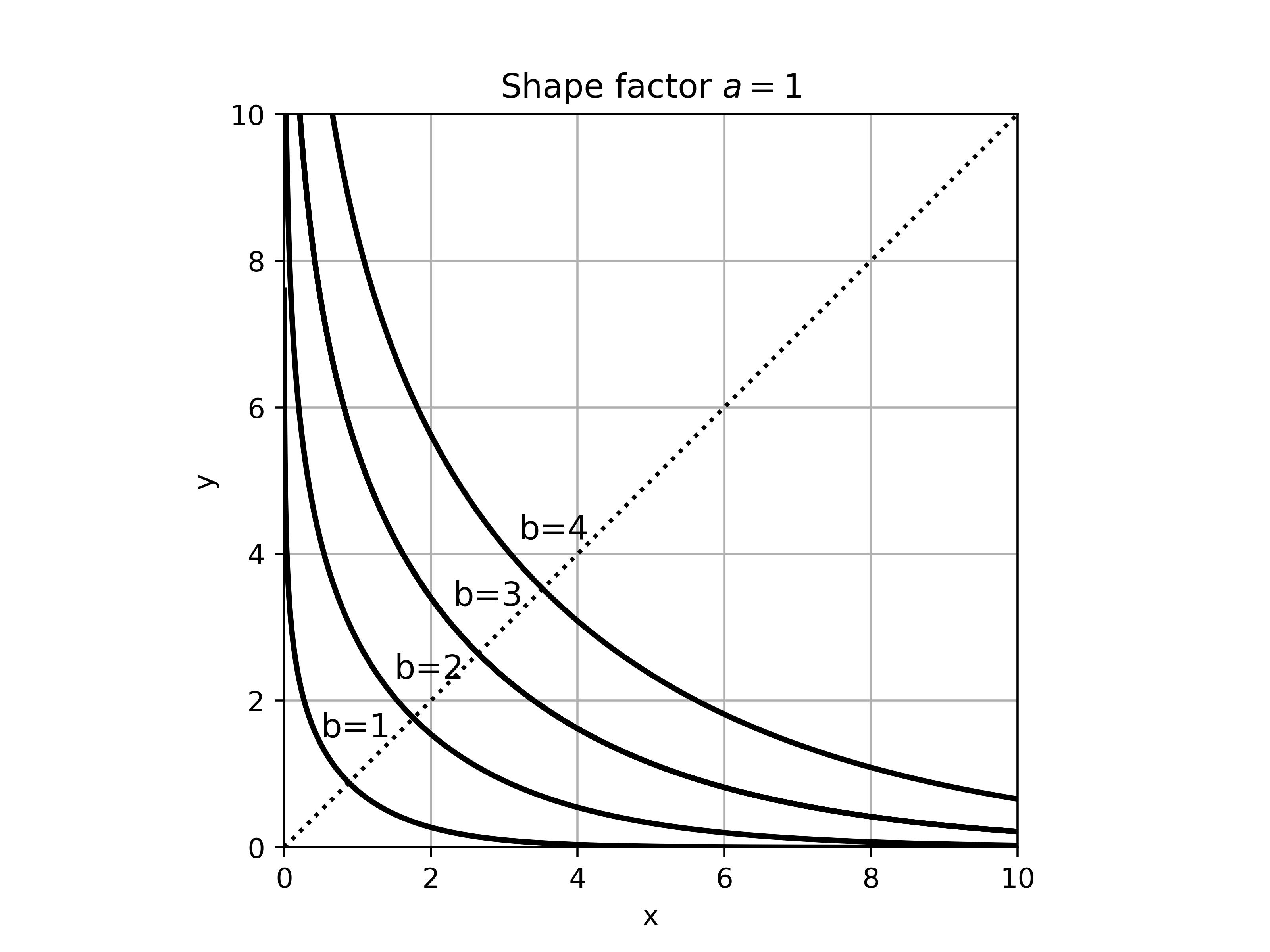}
\end{center}
\caption{\label{Figure2} These given graphs illustrate gemini functions for various values of $b$ and with $a = 1$.}
\end{figure}

 The multiplicative parameter involved in the exponential term in the denominator of the gemini function's fundamental form 
 may be seen as having the effect of changing the steepness of the function, but the function itself retains the self-inverse feature. 
Hence, this parameter is nominated as a shape factor by denoting it with the symbol $a$ further on in this paper. We made another
respective trial by adding a second parameter in the exponents and its reciprocal value in front of the whole function formula. This 
second parameter $b$ scales the size of the function without deforming its shape, and the obtained function is still self-inverse and 
retains its symmetry. We refer to this new function, equipped with the parameters $a$ and $b$, as a \emph{generalized form} of a gemini
function, and we refer to the parameter $b$ as a \emph{scale factor}, and this leads us to the 
 generalized form 
\begin{equation}\label{fullgemini}
\gemini_{a}^b{(x)}=b\ln\left(\frac{1+ae^{-\frac{x}{b}}}{1-e^{-\frac{x}{b}}}\right) 
\end{equation}
 of the gemini function, for parameters $a$ and $b$, noting that the superscript $b$ may be omitted for the $b = 1$ case. Four gemini 
 function graphs with different scale factors are shown in Figure \ref{Figure2}. The shape factor $a$ is equal to $1$ for all these 
 four functions. 

\subsection{Derivation of a five-term gemini identity}
 The integral of the general form of a gemini function is shown in \eqref{integrategemini}. Observe that the total area
increases proportionally to $b^{2}$ for $b>1$, with
\begin{equation}\label{integrategemini}
\int\gemini_{a}^b{(x)} \, dx=\int b\ln\left(\frac{1+ae^{-\frac{x}{b}}}{1-e^{-\frac{x}{b}}}\right) \, dx=
b^{2}\left[\operatorname{Li_{2}}(-ae^{-\frac{x}{b}})-\operatorname{Li_{2}}(e^{-\frac{x}{b}})\right] + C. 
\end{equation}
This may be seen as the first step toward our derivation of the following 5-term dilogarithm relation that is to be heavily
exploited in our work. We refer to this relation as the \emph{five-term gemini identity}. 

\begin{theorem}\label{finalfivegemini}
 The relation 
\begin{multline*}
 \operatorname{Li_{2}}\left(-\frac{a}{x}\right) - 
\operatorname{Li_{2}}\left(\frac{1}{x}\right)+\frac{\pi^{2}}{6}-\operatorname{Li_{2}}\left(-a\right)-
\ln(x)\ln\left(\frac{x+a}{x-1}\right) = \\ -\operatorname{Li_{2}}\left(-a \cdot \frac{x-1}{x+a}\right)+
\operatorname{Li_{2}}\left(\frac{x-1}{x+a}\right)
\end{multline*}
 holds for complex $a$ and $x \not\in \{ 0, 1 \}$ such that $x + a \neq 0$ and such that the arguments of the 
 dilogarithmic expressions given above 
 are such that the associated power series converge absolutely. 
\end{theorem}

\textit{Proof.} We make use of the regions of the plane highlighted in in Figure \ref{Figure4}, with $b = 1$. The apex areas of the form $A_{a}$ 
 are equal when $A_{r}$ is subtracted from the first integral, so that
 \begin{equation}\label{subtractintegral}
\int_{0}^{x_{1}}\gemini_{a}{(x)} \, dx - A_{r}=\int_{x_{2}}^{\infty}\gemini_{a}{(x)} \, dx. 
\end{equation} 
 Let the integration limits $x_1$ and $x_2$ in \eqref{subtractintegral} be such that 
$$x_{1}=\ln(x) \ \ \ \text{and} \ \ \ 
x_{2}=\ln\left(\frac{1+ae^{-\ln(x)}}{1-e^{-\ln(x)}}\right)=\ln\left(\frac{x+a}{x - 1}\right).$$
We thus obtain the area $A_{r}=x_{1}x_{2}=\ln(x)\ln\left(\frac{x+a}{x - 1}\right)$, so that
 \begin{equation}\label{expresslnx}
\int_{0}^{\ln(x)}\left(\frac{1+ae^{-x}}{1-e^{-x}}\right) \, dx-\ln(x)\ln\left(\frac{x+a}{x -1}\right)=
\int_{\ln(\frac{x+a}{x-1})}^{\infty}\left(\frac{1+ae^{-x}}{1-e^{-x}}\right) \, dx, 
\end{equation}
 The first integral in \eqref{expresslnx} may be expressed so that 
$$ \bigg|_{0}^{\ln(x)}\left[\operatorname{Li_2}(-ae^{-x})-\operatorname{Li_2}(e^{-x})\right]-\ln(x)\ln\left(\frac{x+a}{x-1}\right), $$
and this gives us an equivalent version of the desired result, subject to the given conditions. 

\begin{figure}[htbp!]
\begin{center}
 \includegraphics[width=6cm, height=6cm]{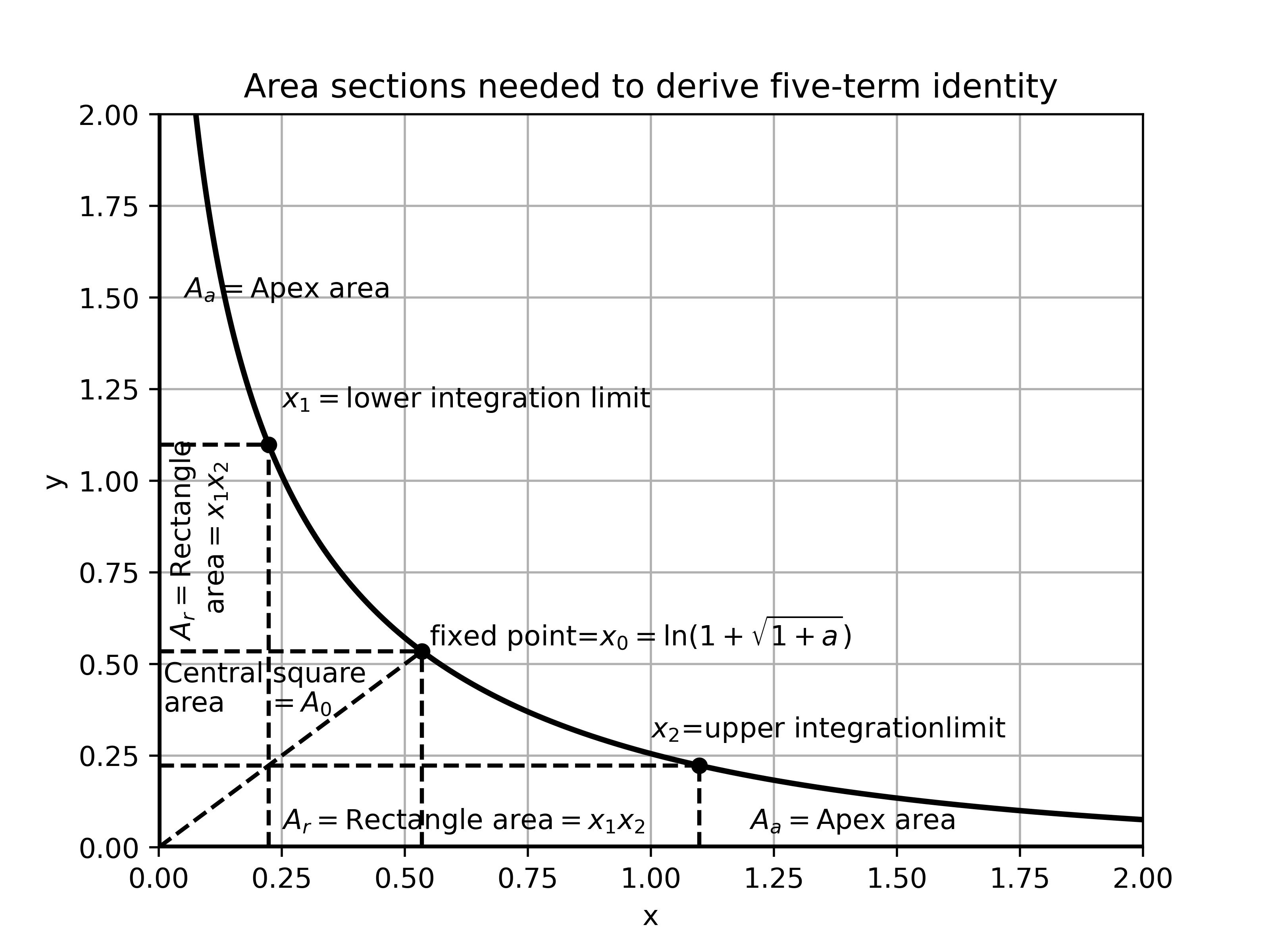}
\end{center}
\caption{\label{Figure3} This plot illustrates the curve of a $\gemini_{-\frac{1}{2}}(x)$-function and required area components needed to
derive dilogarithm identities introduced in this paper.}
\end{figure}

An equivalent version of this identity, for the $a = 1$ case, was recently given in the work of Hakimoglu-Brown \cite{Hak_25}. This derived 
identity reduces down to four-term identity, when the shape factor is equal to +1, because in this case, the third dilogarithm term becomes
a constant value, i.e., $-\operatorname{Li_{2}}(-1 \cdot e^{0})=-\operatorname{Li_{2}(-1)}=\frac{\pi^{2}}{12}$. The valid domain for the shape
factor is such that $a\in[-1,\infty)$. If the shape factor $\ a \ne +1$, then this identity becomes totally different and it enables us to
generate couple of new dilogarithm identities. The five-term gemini-identities obtained from $\gemini_{1}(x)$ and $\gemini_{a}(x)$ at
$x_{1}=\ln(a)$ yield always to one and the same identity. We will deal this issue later on in this paper. 

\subsection{On the derivation of a three-term gemini identity}
 For the fixed point $x_{0}$ of the common integration limit, the five-term identity is reduced to a three-term identity, and this 
 leads us to the following. 

\begin{theorem}\label{firstfixedpoint}
The vanishing
\begin{multline*}
\operatorname{Li_2}\left(-\frac{a}{1+\sqrt{1+a}}\right)-\operatorname{Li_2}\left(\frac{1}{1+\sqrt{1+a}}\right)-
\frac{1}{2}\operatorname{Li_{2}}(-a) + \\ 
\frac{\pi^{2}}{12}-\frac{1}{2}\ln^{2}(1+\sqrt{1+a})=0.
\end{multline*}
 holds for complex $ a \neq 0$ such that the arguments of the above dilogarithmic expressions are such that the associated power series 
 converge absolutely. 
\end{theorem}

\textit{Proof.} For the fixed point $x_0$, the integration limits on the left-hand (resp.\ right-hand) side are from zero to $x_{0}$ (resp.\ from $x_0$ to 
 infinity). The fixed point $x_{0}$ and the shape factor $a$ satisfy $$\ln\left(\frac{1+ae^{-x_{0}}}{1-e^{-x_{0}}}\right)=x_{0} \Rightarrow 
 x_{0}=\ln(x)=\ln(1+\sqrt{1+a}).$$ The area of the associated middle square is such that $$A_{0}=A_{r}=x_{1}x_{2}=x^{2}_{0}=\ln^{2}(1 + 
 \sqrt{1+a}).$$ This leads us to $$\int_{0}^{x_{0}}\gemini_{a}{(x)} \, dx-A_{0}=\int_{x_{0}}^{\infty}\gemini_{a}(x) \, dx. $$ Consequently, 
 we have that $$\int_{0}^{x_{0}}\ln\left(\frac{1+ae^{-x}}{1-e^{-x}}\right) \, dx-x_{0}^{2}=\int_{x_{0}}^{\infty}\ln\left(\frac{1+ae^{-x}}{1 - 
 e^{-x}}\right) \, dx,$$ so that $$\bigg|_{0}^{x_{0}}\left[\operatorname{Li_2}(-ae^{-x})-\operatorname{Li_2}(e^{-x})\right]-x_{0}^{2} = 
\bigg|_{x_{0}}^{\infty}\left[\operatorname{Li_2}(-ae^{-x})-\operatorname{Li_2}(e^{-x})\right], $$ giving us that $$\operatorname{Li_2}\left( 
 -\frac{a}{x_{0}}\right)-\operatorname{Li_2}\left(\frac{1}{x_{0}}\right)+\frac{\pi^{2}}{12}-\frac{1}{2}x_{0}^{2} -\frac{1}{2} 
 \operatorname{Li_{2}}\left(-a\right)=0, $$ giving us an equivalent formulation of the desired result. 

The relation in Theorem \ref{firstfixedpoint} provides a new three-term single variable dilogarithm identity. We refer to this as the 
 \emph{first fixed-point gemini identity}, noting that an identity of a similar nature is to later be derived in our work, 
and where the arguments of the dilogarithm terms are expressed with the aid of the fixed point values, i.e., with $x_{0}$. 
A quite similar identity is introduced in the Lewin's monograph on dilogarithms and associated functions \cite{Lew_58}. The
derivation of this second three-term gemini-identity is analogous with respect to the previous derivation. Now, we need to
express the shape factor as a function of the argument of the fixed point $x$, where $x_{0}=\ln(x)$. Hence, $x=1+\sqrt{1+a}
\Rightarrow a=(x-1)^{2}-1=x^{2}-2x$. We thus obtain that
\begin{equation}\label{threetermfromfixed}
\operatorname{Li_{2}}\left(2-x\right)-\operatorname{Li_{2}}\left(\frac{1}{x}\right)-\frac{1}{2}\operatorname{Li_{2}}\left(2x-x^{2}\right)+
\frac{\pi^{2}}{12}-\frac{1}{2}\ln^{2}(x)=0, \quad x>1 
\end{equation}

\subsection{A degenerate form of a gemini function}
As already explained, the acceptable domain for the shape factor is defined in such a way that $a \geq-1$. Next, we deal with the
gemini function equipped with $a=0$. The exponential term vanishes in the nominator. Hence, this kind of gemini function is called
a degenerate form of a gemini function. The formula for a degenerate gemini function and its integral are given by
$$ \int\gemini_{0}{\left(x\right)} \, dx=\int\ln\left(\frac{1}{1-e^{-x}}\right) \, dx=-\operatorname{Li_{2}}(e^{-x}) + C. $$

\begin{figure}[htbp!]
\begin{center}
 \includegraphics[width=6cm, height=6cm]{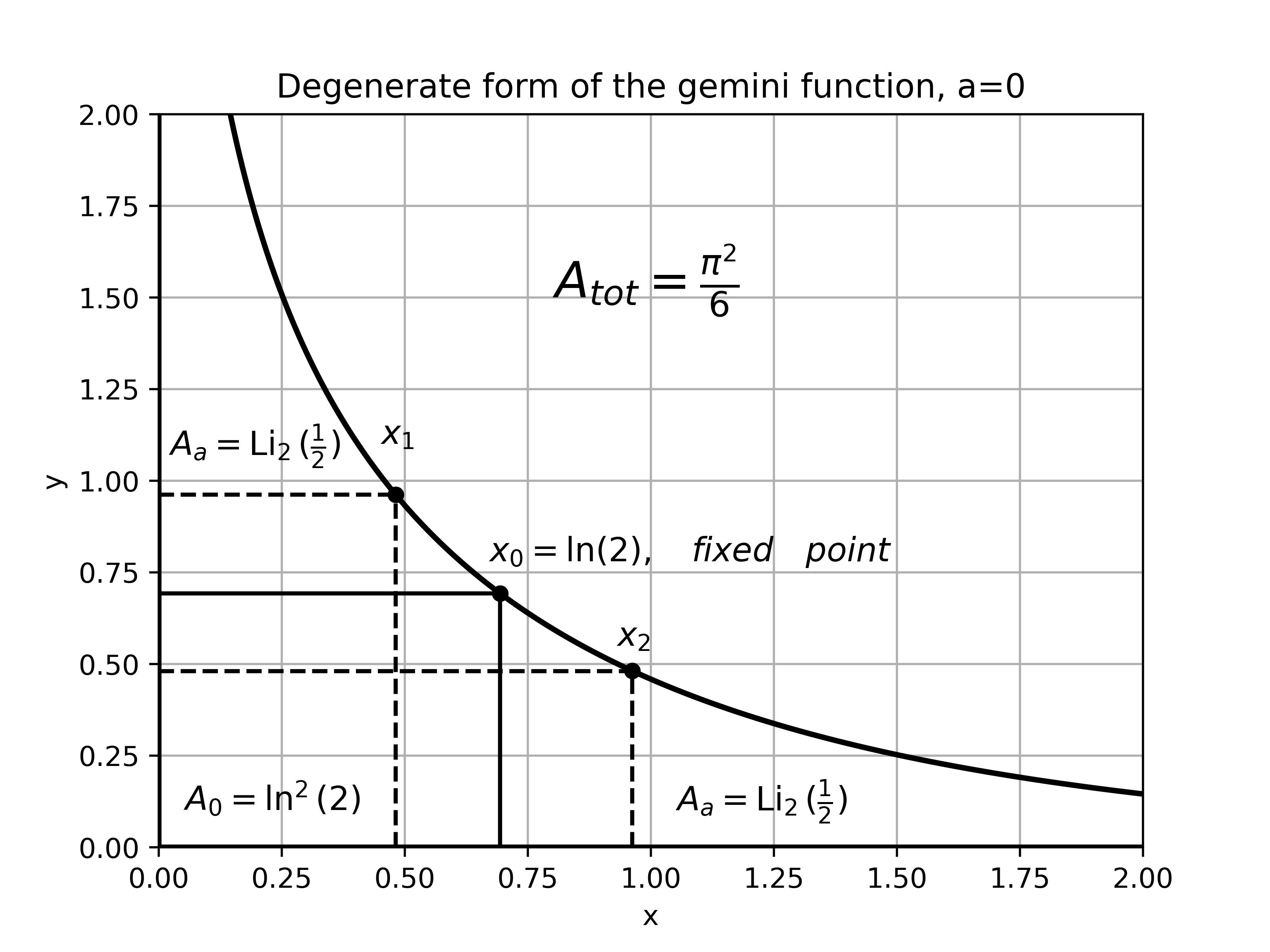}
 \end{center}
\caption{\label{Figure4} The graph of a degenerate gemini function and the schematics of area sections needed to derive the reflection
formula are illustrated in this figure. The fixed point of a degenerate form is given by $x_{0}=\ln(1+\sqrt{1+0})=\ln(2)$. The area of
a middle square is such that $A_{0}=\ln^{2}(2)$. The corresponding two apex areas are equal, which are given by
$A_{a}=\frac{1}{2}[A_{tot}-A_{0}]=\frac{1}{2}[\frac{\pi^{2}}{6}-\ln^{2}(2)]=\frac{\pi^{2}}{12}-\frac{1}{2}\ln^{2}(2)=
\operatorname{Li}_{2}\left(\frac{1}{2}\right)$.}
\end{figure}

The reflection identity, which is also called Euler's identity, is easy to derive by applying the degenerate gemini function. The
graphics in Figure \ref{Figure4} illustrates the area sections needed to build the equation for this identity. The relation between
the integration limits is such that $x_{1}=\ln(x)$ and $x_{2}=\ln(\frac{x}{x -1})$ for $x>1$. Hence, we write
$$\int_{0}^{x_{1}}\gemini_{0}{(x)} \, dx-x_{1}x_{2}=\int_{x_{2}}^{\infty}\gemini_{0}{(x)} \, dx, $$ which implies that 
$$\int_{0}^{\ln(x)}\ln\left(\frac{1}{1-e^{-x}}\right) \, dx-\ln(x)\ln\left(\frac{x}{x-1}\right)=
\int_{\ln(\frac{x}{x - 1})}^{\infty}\ln\left(\frac{1}{1-e^{-x}}\right) \, dx, $$
which, in turn, gives us that 
$$\bigg|_{0}^{\ln(x)}-\operatorname{Li_2}(e^{-x})-\ln(x)\ln\left(\frac{x}{x - 1}\right)=
\bigg|_{\ln(\frac{x}{x-1})}^{\infty}-\operatorname{Li_2}(e^{-x}), $$
so that 
$$\frac{\pi^{2}}{6}-\operatorname{Li_{2}}\left(\frac{1}{x}\right)-\ln(x)\ln\left(\frac{x}{x-1}\right)=
\operatorname{Li_{2}}\left(\frac{x-1}{x}\right),$$
and we thus arrive at the reflection formula

\begin{eqnarray}\label{reflectionformula}
\operatorname{Li_{2}}\left(\frac{1}{x}\right)+\operatorname{Li_{2}}\left(1-\frac{1}{x}\right)=\frac{\pi^{2}}{6}-
\ln(x)\ln\left(\frac{x}{x -1}\right), \quad x>1.
\end{eqnarray}

This identity can be simply obtained also by setting the shape factor $a$ equal to zero in Theorem \ref{finalfivegemini}.
The exact value of the $\operatorname{Li_{2}}(\frac{1}{2})$ can be calculated simply by using this degenerate gemini function.
The area sections of a degenerate gemini function have interesting values, e.g. {$A_{tot}=\int_{0}^{\infty}\gemini_{0}(x) \, dx=
\frac{\pi^{2}}{6}$, $A_{0}=\ln^{2}(2)$ and $A_{a}=\int_{\ln(2)}^{\infty}\gemini_{0}(x) \, dx=\operatorname{Li_{2}}(\frac{1}{2})$ as drawn
in Figure \ref{Figure5}. Hence, $\operatorname{Li_{2}}(\frac{1}{2})=\frac{\pi^{2}}{12}-\frac{1}{2}\ln^{2}(2)$, which is one of the eighth known
exact real values of a dilogarithm. The seven others are: $\operatorname{Li_{2}}(0)=0, \quad \operatorname{Li_{2}}(1)=\frac{\pi^{2}}{6},
\quad \operatorname{Li_{2}}(-1)=-\frac{\pi^{2}}{12}, \quad \operatorname{Li_{2}}(-\frac{1}{\phi})=
-\frac{\pi^{2}}{15}+\frac{1}{2}\ln^{2}(\phi),\quad \operatorname{Li_{2}}(\frac{1}{\phi})=\frac{\pi^{2}}{10}-\ln^{2}(\phi),\quad
\operatorname{Li_{2}}(\frac{1}{\phi^{2}})=\frac{\pi^{2}}{15}-\ln^{2}(\phi)$ and $\operatorname{Li_{2}}(-\phi)=-\frac{\pi^{2}}{10}-\ln^{2}(\phi)$.

\subsection{The inversion identity, Landen's identity, and the duplication formula}

Some of the most basic identities involving the dilogarithm function, such as the inversion formula 

\begin{eqnarray}\label{inversionformula} 
\operatorname{Li_{2}}(-x)+\operatorname{Li_{2}}\left(-\frac{1}{x}\right)+\frac{\pi^{2}}{6}+\frac{1}{2}\ln^{2}(x)=0, \quad x>1, 
\end{eqnarray}

can be obtained using gemini functions. Omitting details, by rotating the degenerate gemini function counterclockwise by an
angle of $\frac{\pi}{4}$, we find that $x_{2}=x_{1}\cos(\theta)-y_{1}\cos(\theta)$ and $y_{2}=x_{1}\sin(\theta)+y_{1}\cos(\theta)$, 
and, since $ \sin(\theta)=\cos(\theta)=\frac{1}{\sqrt{2}}$, we obtain 
$$x_{2}=x=\ln(t)\frac{1}{\sqrt{2}}-\ln\left(\frac{t}{t-1}\right)\frac{1}{\sqrt{2}} \Rightarrow x\sqrt{2}=\ln(t-1) \Rightarrow t=e^{x\sqrt{2}}+1, $$ 
and a substitution of the form $t=e^{x\sqrt{2}}+1$ can then be used to obtain that 
$$\int\gemini_{0}^{rot}(x) \, dx=\int\frac{1}{\sqrt{2}}\ln\left(2\cosh(x\sqrt{2})+2\right) \, 
 dx=\operatorname{Li_{2}}\left(-e^{-x\sqrt{2}}\right)+
\frac{1}{2}x^{2}+C. $$
A similar approach can be used to derive \emph{Landen's formula}

\begin{eqnarray}\label{Landensformula}
\operatorname{Li_{2}}\left(\frac{1}{1+x}\right)-\operatorname{Li_{2}}(-x)=\frac{\pi^{2}}{6}-\frac{1}{2}\ln(1+x)\ln\left(\frac{1+x}{x^{2}}\right),
\quad x>0. 
\end{eqnarray}

Omitting details, by setting $x_{11}=\ln(x)$, and by then letting the upper integration limit of the degenerate gemini function be such that
$x_{12}=\ln\left(\frac{x}{x-1}\right)$, so that the corresponding rotated $x$-coordinates are such $x_{21}=-\frac{1}{\sqrt{2}}\ln(x-1)$ and
$x_{22}=\frac{1}{\sqrt{2}}\ln\left(\frac{1}{x-1}\right)$, and by then using the relation 
\begin{multline*}
\int_{\ln(x)}^{\ln(\frac{x}{x-1})}\gemini_{0}(x) \, dx = \\ \operatorname{Li_{2}}\left(\frac{1}{x}\right)-
\operatorname{Li_{2}}\left(\frac{x}{x-1}\right)=2\operatorname{Li_{2}}\left(\frac{1}{x}\right)-
\frac{\pi^{2}}{6}+\ln(x)\ln\left(\frac{x}{x-1}\right). 
\end{multline*}
Gemini functions can also be used to derive the famous duplication formula such that

\begin{eqnarray}\label{duplicationformula}
\operatorname{Li_{2}}\left(\frac{1}{x}\right)+\operatorname{Li_2}\left(-\frac{1}{x}\right)=
\frac{1}{2}\operatorname{Li_{2}}\left(\frac{1}{x^{2}}\right).
\end{eqnarray}

\subsection{Further three-term relations}
 In this section, two three-term dilogarithm identities are derived starting from the five-term gemini-identity. The cancellation 
 of two terms out of the five is based on the selection of the suitable initial values. 

\begin{theorem}\label{mainthreeterm}
 The vanishing 
\begin{multline*}
\operatorname{Li_{2}}\left(\frac{a-1}{a^{2}}\right)+\operatorname{Li_{2}}\left(\frac{1}{a^{2}-a+1}\right)-
\operatorname{Li_{2}}\left(\frac{a}{a^{2}-a+1}\right) + \\ 
\ln\left(\frac{a}{a-1}\right)\ln\left(\frac{a^{2}}{a^{2}-a+1}\right)= 0 
\end{multline*} 
 holds for complex $a \not\in \{ 0, 1, \sqrt[3]{-1}, -(-1)^{2/3} \}$ such that the arguments of the above dilogarithmic expressions are such 
 that the associated power series converge absolutely. 
\end{theorem}
 
\textit{Proof.} Let the shape factor be such that $a=-\frac{1}{a}$ for $a>1$. Hence, the integration limits are such that $x_{1}=\ln(\frac{a}{a-1})$
and $x_{2}=\ln(\frac{a^{2}-a+1}{a})$. Now, the integration limits are expressed as a function of the scale factor $a$. By substituting
these initial values in \eqref{finalfivegemini}, we can write
\begin{multline*}
\operatorname{Li_{2}}\left(\frac{a-1}{a^{2}}\right)-\operatorname{Li_{2}}\left(\frac{a-1}{a}\right)+\frac{\pi^{2}}{6}-
\operatorname{Li_{2}}\left(\frac{1}{a}\right)-\ln\left(\frac{a}{a-1}\right)\ln\left(\frac{a^{2}-a+1}{a}\right)= \\ 
-\operatorname{Li_{2}}\left(\frac{1}{a^{2}-a+1}\right)+\operatorname{Li_{2}}\left(\frac{a}{a^{2}-a+1}\right). 
\end{multline*}
We proceed to apply the reflection identity to the second dilogarithm term, so that the resultant term is equal to the negative of
the third dilogarithmic term, and hence the conversion given by 
$$-\operatorname{Li_{2}}\left(\frac{a-1}{a}\right)=\operatorname{Li_{2}}\left(\frac{1}{a}\right)+\ln\left(\frac{a}{a-1}\right)\ln(a)-
\frac{\pi^{2}}{6}, $$
and this gives us the desired result. 

A similar cancellation method, relative to our proof of Theorem \ref{mainthreeterm}, leads us toward the three-term dilogarithm
identity highlighted below. 

\begin{theorem}\label{threetermconstant}
The relation 
\begin{multline*}
\operatorname{Li_{2}}\left(\frac{a}{(a+1)^{2}}\right)+\operatorname{Li_{2}}\left(\frac{1}{a^{2}+a+1}\right)-
\operatorname{Li_{2}}\left(\frac{a+1}{a^{2}+a+1}\right) + \\ 
\ln\left(\frac{a+1}{a}\right)\ln\left(\frac{(a+1)^{2}}{a^{2}+a+1}\right)=0, 
\end{multline*}
for complex $a$ such that $a \not\in \{ 0, -1, -\sqrt[3]{-1}, (-1)^{2/3} \}$ and such that the series corresponding to the dilogarithmic
terms given above are convergent. 
\end{theorem}

\begin{figure}[htbp!]
\begin{center}
 \resizebox{9.0cm}{!}
 {\includegraphics{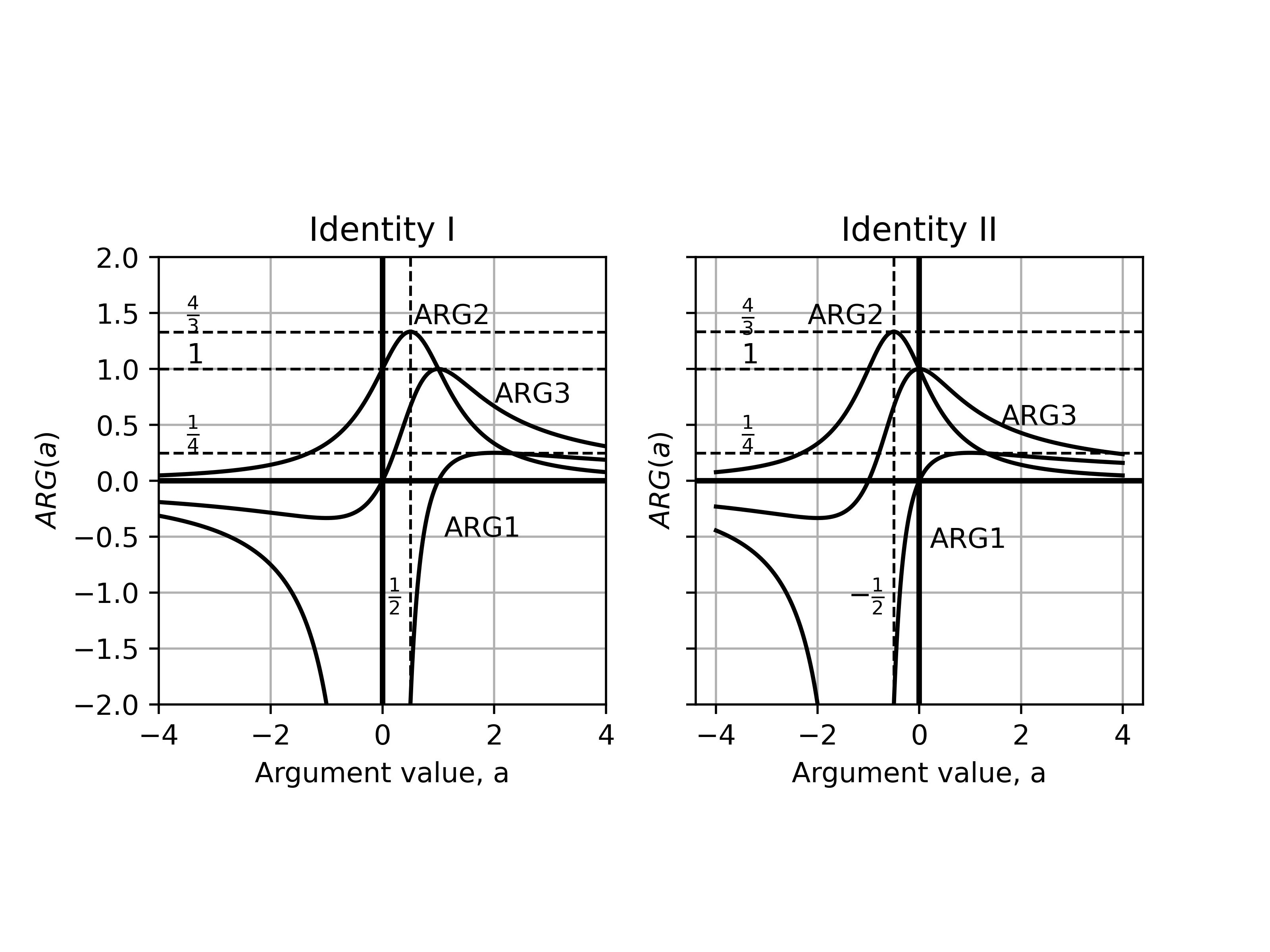}}
\end{center}
\caption{\label{Figure5} The graphs on the left side illustrates the arguments of all three terms of the first identity as a function of the shape
factor $a$. The graphs on the right side illustrates the arguments of all three terms of the second identity as a function of the shape
factor $a$.}
\end{figure}

\textit{Proof.} Letting the shape factor be such that $a>0$, we proceed to let the integration limits be such that $x_{1}=\ln\left(\frac{a+1}{a}\right)$ and
$x_{2}=\ln\left(a^{2}+a+1\right)$. Now, the integration limits are also expressed as functions of the scale factor $a$ like earlier. By
substituting these initial values in Theorem \ref{finalfivegemini}, we obtain 
\begin{multline*}
\operatorname{Li_{2}}\left(-\frac{a^{2}}{a+1}\right)-\operatorname{Li_{2}}\left(\frac{a}{a^{2}+1}\right)+\frac{\pi^{2}}{6}-
\operatorname{Li_{2}}(-a)-\ln\left(\frac{a+1}{a}\right)\ln\left(a^{2}+a+1\right) = \\ 
-\operatorname{Li_{2}}\left(-\frac{a}{a^{2}+a+1}\right)+
\operatorname{Li_{2}}\left(\frac{1}{a^{2}+a+1}\right), 
\end{multline*}
which is equivalent to the desired result. 

According to Figure \ref{Figure5}, it is easy to realize that these two three-term identities are otherwise similar, but they just differ from 
each other by a unit translation along the $x$-axis. The arguments of the both identities are plotted here for the becoming purposes,
because these identities are functional partially also in the complex domain. The maximum value of the first terms in the both identity
is only {$\frac{1}{4}$}. The limiting value is 1, if we purely deal with in the real domain. We will briefly investigate this behavior
of gemini-identities in the complex domain later on in this paper. 

\subsection{A framework for generating two-term dilogarithm identities}\label{framework}
 Theorems \ref{finalfivegemini}, \ref{firstfixedpoint}, 
 and \ref{mainthreeterm} may be seen as providing a framework 
 for the generation of new two-term dilogarithm identities, as we demonstrate in the current subsection. 

\begin{theorem}\label{firsttwotermtheorem}
 The closed-form evaluation 
\begin{multline*}
\operatorname{Li_{2}}\left(\frac{1}{2\phi}+\frac{1}{2}\sqrt{\frac{\phi^{2}+3}{\phi^{3}}}\right)-
\operatorname{Li_{2}}\left(\frac{1}{2\phi}-\frac{1}{2}\sqrt{\frac{\phi^{2}+3}{\phi^{3}}}\right)=\frac{\pi^{2}}{10}+\frac{1}{2}\ln^{2}(\phi) - \\ 
 \ln\left(\frac{1}{2\phi}+\frac{1}{2}\sqrt{\frac{\phi^{2}+3}{\phi^{3}}}\right)
\ln\left(\frac{\phi^{2}+1}{2\phi^{2}}-\frac{1}{2}\sqrt{\frac{\phi^{2}+3}{\phi^{3}}}\right)-
\ln(\phi)\ln\left(\frac{1}{2}\sqrt{\phi^{2}+4\phi}-\frac{1}{2}\phi\right)
\end{multline*}
 holds. 
\end{theorem}

\textit{Proof.} Our strategy is to apply Theorem \ref{mainthreeterm} by fixing the first term, by setting the argument of the initial term as $-\frac{1}{\phi}$. 
 Explicitly, by setting $\frac{a-1}{a^{2}}=-\frac{1}{\phi}$, we find that $ a = -\frac{-\phi-\sqrt{\phi^{2}+4\phi}}{2}$, so that Theorem 
 \ref{mainthreeterm} gives us the vanishing relation 
\begin{multline*}
 \operatorname{Li_{2}}\left(-\frac{1}{\phi}\right)+
\operatorname{Li_{2}}\left(\frac{\phi^{2}+1}{2\phi^{2}}-\frac{1}{2}\sqrt{\frac{\phi^{2}+3}{\phi^{3}}}\right) - \\ 
\operatorname{Li_{2}}\left(\frac{1}{2\phi}-\frac{1}{2}\sqrt{\frac{\phi^{2}+3}{\phi^{3}}}\right)-
\ln(\phi)\ln\left(\frac{\sqrt{\phi^{2}+4\phi}-\phi}{2}\right)=0, 
\end{multline*}
 so that the known closed form 
 for $\operatorname{Li_{2}\left(-\frac{1}{\phi}\right)} = -\frac{\pi^2}{15} + \frac{\log^2 \phi}{2}$ 
 gives us the desired result. 

 Relative to our derivation of Theorem \ref{firsttwotermtheorem}, 
 a simialr approach relying on the known closed-form evaluation
 for $\operatorname{Li}_{2}\big( \frac{1}{\phi^2} \big) = \frac{\pi^2}{15} - \log^2 \phi$ 
 leads us to the following result. 

\begin{theorem}
 The closed-form evaluation 
\begin{multline*}
 \operatorname{Li_{2}}\left(\frac{\sqrt{17\phi^{2}+11\phi}-3\phi-\phi^{2} }{2} \right)+ 
\operatorname{Li_{2}}\left(\frac{\phi^{2}+1}{2\phi^{2}}-\frac{1}{2}\sqrt{\frac{\phi^{2}+3}{\phi^{3}}}\right) = \\ 
\frac{\pi^{2}}{15}-\ln^{2}(\phi)
 -\ln\left(\frac{1}{2}+\frac{1}{2}\sqrt{\frac{\phi^{2}+3}{\phi}}\right)\ln\left(\frac{\phi^{2}+1}{2\phi^{2}}+\frac{1}{2}\sqrt{\frac{\phi^{2}+3}{\phi^{3}}}\right)
\end{multline*} 
 holds. 
\end{theorem}

\textit{Proof.} Applying Theorem \ref{mainthreeterm}, we set the argument of the third term as $\frac{1}{\phi^2}$.
Explicitly, by setting $\frac{a}{a^2-a + 1}=\frac{1}{\phi^{2}}$, we obtain that $ a = \frac{\phi^{2}+1}{2}+\frac{1}{2}\sqrt{\phi^{2}+4\phi}$.
 Inserting this specified value of $a$ into Theorem \ref{mainthreeterm}, 
 and exploiting the known closed form for $\operatorname{Li}_{2}\big( \frac{1}{\phi^2} \big)$, 
 we obtain an equivalent version of the desired result. 

 Similarly, using an analytic continuation for $\phi$ 
 together with the known closed form for $\operatorname{Li}_{2}(\phi)$
 (equivalent to $\operatorname{Li}_{2}\big( \frac{1}{\phi} \big) = \frac{\pi^2}{10} - \log^2 \phi$ by the inversion formula), 
 we obtain the following. 

\begin{theorem}\label{thirdtwotermtheorem}
 The closed-form evaluation 
$$ \operatorname{Li}_{2}\left( \frac{1}{2} + \frac{i}{2} \sqrt{\frac{1 + \phi^2}{\phi^4}} \right) + 
 \operatorname{Li}_{2}\left( \frac{\phi^2 - i \sqrt{1 + \phi^2}}{2} \right) = \\ 
 \frac{23 \pi ^2}{150}-\frac{\log ^2(\phi )}{2}-\frac{4}{5} i \pi \log (\phi ) $$ holds.
\end{theorem}

\textit{Proof.} Setting the third term of Theorem \ref{mainthreeterm} to be equal
 to $\phi$, so that $\frac{a}{a^2-a+1}=\phi \Rightarrow a=(-1)^{\frac{1}{5}}=e^{\frac{i\pi}{5}}$, 
 by setting the specified value for $a$ into Theorem
 \ref{mainthreeterm}, we obtain the vanishing relation 
\begin{multline*}
 \operatorname{Li_{2}}\left(\frac{1}{2}+\frac{1}{2}i\sqrt{\frac{\phi^{2}+1}{\phi^{4}}}\right)+
\operatorname{Li_{2}}\left(\frac{1}{2}\phi^{2}-\frac{1}{2}i\sqrt{\phi^{2}+1}\right)-
\operatorname{Li_{2}}\left(\phi\right)+ \\ 
\ln\left(\frac{1}{2}-\frac{1}{2}i\sqrt{4\phi+3}\right)
\ln\left(\frac{1}{2}\phi^{2}+\frac{1}{2}i\sqrt{\phi^{2}+1}\right)=0, 
\end{multline*}
 so that the closed form for $\operatorname{Li_{2}}\left(\phi\right)$ gives us the desired result. 

 Many further results of a similar nature can be obtained, through the application of our results as in Theorems 
 \ref{finalfivegemini}--\ref{mainthreeterm} together with the major functional equations for $\operatorname{Li}_{2}$, such as
\begin{equation}\label{majorfunctional}
 \operatorname{Li}_{2}(-x) - \operatorname{Li}_{2}(1-x) + \frac{1}{2} \operatorname{Li}_{2}(1-x^2) 
 = -\frac{\pi^2}{12} - \log x \log(x+1). 
\end{equation}
 For example, by applying \eqref{majorfunctional} to the first term of Theorem \ref{mainthreeterm}, it can be shown that 
\begin{multline*}
 \operatorname{Li}_2\left(\frac{a-1}{a^2}-1\right)+\operatorname{Li}_2\left(\frac{1}{a^2-a+1}\right) - \\
 \operatorname{Li}_2\left(\frac{a}{a^2-a+1}\right)+\frac{1}{2}
 \operatorname{Li}_2\left(\frac{(a-1) \left(2 a^2-a+1\right)}{a^4}\right) 
\end{multline*}
 and 
$$ -\log \left(1-\frac{a-1}{a^2}\right) \log \left(2-\frac{a-1}{a^2}\right)-\log \left(\frac{a}{a-1}\right) \log
 \left(\frac{a^2}{a^2-a+1}\right)-\frac{\pi ^2}{12}, $$ are equal. 
 By then setting $a = \frac{1}{4} \left(-1-\sqrt{5}-\sqrt{14+10 \sqrt{5}}\right)$, 
 this gives us an equivalent version of final two-term dilogarithm evaluation highlighted as a motivating example 
 in Section \ref{subsectionmotivating}. 
 
\section{Application examples of gemini identities in the real domain}
 We have now rederived the four main well known dilogarithm identities by applying the properties of the gemini functions, which
are the reflection, the inversion, Landen's and the duplication formula. Next, we introduce the suitability of the gemini-identities
for evaluating exact values for certain dilogarithms, two-term value identities and ladders. 

\subsection{On the evaluation of $\operatorname{Li_{2}}(\frac{1}{\phi^{2}})$}
 Let $\phi$ denote the golden ratio, with $\phi=\frac{1+\sqrt{5}}{2}$. By setting the shape factor so that $a=+\phi^{2}$,
and by setting the lower integration limit so that $x_{1}=\ln(\phi)$, the upper integration limit is $x_{2} = \ln(\phi^{4})$, so that the 
five-term gemini-identity gives us that 
\begin{multline*}
\operatorname{Li_{2}}\left(-\frac{\phi^{2}}{\phi}\right) - \operatorname{Li_{2}}\left(\frac{1}{\phi}\right)+\frac{\pi^{2}}{6}-
\operatorname{Li_{2}}(-\phi^{2})-\ln(\phi)\ln(\phi^{4}) = \\ 
-\operatorname{Li_{2}}\left(-\frac{\phi^{2}}{\phi^{4}}\right)+
\operatorname{Li_{2}}\left(\frac{1}{\phi^{4}}\right), 
\end{multline*}
which, in turn, implies that 
\begin{multline*}
 \operatorname{Li_{2}}(-\phi) - 
 \operatorname{Li_{2}}\left(
 \frac{1}{\phi} \right) + 
 \frac{\pi^{2}}{6}+\operatorname{Li_{2}}\left(-\frac{1}{\phi^{2}}\right)+
\frac{\pi^{2}}{6}+\frac{1}{2}\ln^{2}(\phi^{2})-4\ln^{2}(\phi)= \\ 
 -\operatorname{Li_{2}}\left(-\frac{1}{\phi^{2}}\right)+\operatorname{Li_{2}}\left(\frac{1}{\phi^{4}}\right), 
\end{multline*}
so that routine rearrangements allow us to derive the closed-form evaluation such that 
$$ \operatorname{Li_{2}}\left(\frac{1}{\phi^{2}}\right)=\frac{\pi^{2}}{15}-\ln^{2}(\phi).$$

\subsection{On Legendre's chi-function evaluated at $\frac{1}{\phi^{3}}$}

There have been a number of recent publications concerning the derivation of closed forms for
two-term dilogarithm identities through 
a variety of different methods \cite{Ade_24,Campbell202122,Lima2024,Ste22}. If the arguments in the two-term value identity are equal, 
but having opposite signs, then this identity can be simply represented with the 
aid of the Legendre's chi-function \cite[\S1.8]{Lewin1981}, which may be defined so that 
$$\operatorname{\chi_{2}}\left(x\right)=\frac{1}{2}\left[\operatorname{Li_{2}}(x)-\operatorname{Li_{2}}(-x)\right].$$ 
The closed-form evaluation for $\operatorname{\chi_{2}}\left(\frac{1}{\phi^{3}}\right)$ is classical known \cite[\S1.8]{Lewin1981}
and follows from one of the standard functional equations for $\chi_2$, but the construction of alternative proofs for the specified
closed form is motivated by a recent Fourier--Legendre-based proof of this closed form given by Campbell \cite{Campbell202122} and by
how such alternative proofs could be useful in the construction of further closed forms for two-term dilogarithm relations. 

\begin{theorem}
The closed-form evaluation 
$ \operatorname{\chi_{2}}\big(\frac{1}{\phi^{3}} \big)=\frac{\pi^{2}}{24}-\frac{3}{4}\ln^{2}(\phi)$
holds \cite[p.\ 19]{Lewin1981}.
\end{theorem}

\textit{Proof.} We apply five-term gemini-identity related to the $\gemini_{1}(x)$-function, setting $a = 1$. Letting the associated limits of
integration be such that $x_{1}=\ln(\phi)$ and $x_{2}=\ln\left(\frac{1 + 
\phi}{1-\phi}\right)=\ln(\phi^{3})$, 
and by inputting these values in the five-term gemini-identity, we find that 
\begin{multline*}
\operatorname{Li_{2}}\left(-\frac{1}{\phi}\right) - 
\operatorname{Li_{2}}\left(\frac{1}{\phi}\right)+\frac{\pi^{2}}{6}-\operatorname{Li_{2}}(-1)
-\ln(\phi)\ln\left(\frac{\phi+1}{\phi-1}\right) = \\ 
-\operatorname{Li_{2}}\left(-\frac{\phi-1}{\phi+1}\right)+\operatorname{Li_{2}}\left(\frac{\phi-1}{\phi+1}\right).
\end{multline*}
Consequently, we have that 
$$\operatorname{Li_{2}}\left(-\frac{1}{\phi}\right)-\operatorname{Li_{2}}\left(\frac{1}{\phi}\right)+\frac{\pi^{2}}{4}-3\ln^{2}(\phi)= 
-\operatorname{Li_{2}}\left(-\frac{1}{\phi^{3}}\right)+\operatorname{Li_{2}}\left(\frac{1}{\phi^{3}}\right), $$
and this gives us the desired result from the known closed forms for $\text{Li}_{2}\big( -\frac{1}{\phi} \big)$ and for
$\text{Li}_{2}\big( \frac{1}{\phi} \big)$. 

\subsection{On the Legendre chi-function evaluated at $\sqrt{2} - 1$}

The classically known closed-form evaluation for $\chi_2$ at $\sqrt{2}-1$ \cite[p.\ 19]{Lewin1981} has been considered recently by a 
number of authors \cite{Campbell202122,Lim12,Ste22}, inspired by the work of Lima \cite{Lim12} concerning a double integral 
employed by Beukers et al.\ to solve the Basel problem. We derive the closed form for $\chi_2(\sqrt{2}-1)$ with the aid of the first
three-term fixed point gemini-identity. As noted in a recent survey concerning two-term dilogarithm relations \cite{Campbell2022},
the classically known closed form for $\text{Li}_{2}\big( \sqrt{2} - 1 \big) - \text{Li}_{2}\big( 1 - \sqrt{2} \big)$ can be obtained,
in equivalent ways, from results given by Bytsko \cite{Bytsko1999,Byt01_99} and by Richmond and Szekeres \cite{RichmondSzekeres1981}. 

\begin{theorem}
The closed-form evaluation 
$\operatorname{\chi_{2}}\big(\sqrt{2}-1\big)= 
\frac{\pi^{2}}{16}-\frac{1}{4}\ln^{2}(\sqrt{2} - 1)$
holds \cite[p.\ 19]{Lewin1981}. 
\end{theorem}

\textit{Proof.}For the gemini function $\gemini_{1}(x)$, with $a = 1$, we set $x_{0}=\sqrt{2}-1$. Inputting the given values into the first fixed-point
identity in Theorem \ref{firstfixedpoint}, we find that 
 $$ \operatorname{Li_2}\left(-\frac{1}{1+\sqrt{2}}\right)-\operatorname{Li_2}\left(\frac{1}{1+\sqrt{2}}\right)-
\frac{1}{2}\operatorname{Li_{2}}(-1)+\frac{\pi^{2}}{12}-\frac{1}{2}\ln^{2}(1+\sqrt{2}) = 0, $$
and the desired result then follows by evaluating $\text{Li}_{2}(-1)$ in closed form. 

\begin{remark}
From the classically known evaluation for $ \operatorname{\chi_{2}}\big( \sqrt{2} - 1 \big)$ together with Landen's identity, we can
recover a relationship between $\text{Li}_{2}\big( \frac{1}{\sqrt{2}} \big)$ and $\text{Li}\big( \sqrt{2} - 1 \big)$ considered by 
Bytsko \cite{Byt01_99} in the context of of applications in conformal field theory and subsequently considered by Lima \cite{Lim12}. 
\end{remark}

\subsection{On a relation between $\operatorname{Li_{2}}\big(\sqrt{2}-1\big)$ and $\operatorname{Li_{2}}\big(\frac{2-\sqrt{2}}{4}\big)$}
 Using the five-term gemini identity, we have obtained a closed-form evaluation for a combination 
 of $\operatorname{Li_{2}}\big(\sqrt{2}-1\big)$
 and $\operatorname{Li_{2}}\big(\frac{2-\sqrt{2}}{4}\big)$, as below. This two-term dilogarithm relation appears to be original. 

\begin{theorem}
The closed-form evaluation 
\begin{multline*}
6\operatorname{Li_{2}}\left(\sqrt{2}-1\right)+\operatorname{Li_{2}}\left(\frac{2-\sqrt{2}}{4}\right) = 
\frac{11\pi^2}{24}-\frac{3}{8}\ln^2(2)- \\ 
\frac{3}{2}\ln^{2}(\sqrt{2}+1)-\frac{3}{2}\ln(2)\ln(2+\sqrt{2})+\ln(\sqrt{2}+1)\left[\frac{1}{2}\ln(2)+\ln(2+\sqrt{2})\right]
\end{multline*}
holds. 
\end{theorem}

\textit{Proof.} We apply the five-term gemini identity, by setting $a=-2\sqrt{2}+2$ and $x_{1}=\ln(4-2\sqrt{2})$ and $x_{2}=\ln(2)$, so that 
\begin{multline*}
\operatorname{Li_{2}}\left(\frac{2\sqrt{2}-2}{4-2\sqrt{2}}\right)-\operatorname{Li_{2}}\left(\frac{1}{4-2\sqrt{2}}\right)
-\operatorname{Li_{2}}\left(2\sqrt{2}-2\right)+\frac{\pi^{2}}{6} - \\ 
\ln\left(4-2\sqrt{2}\right)\ln(2) = -\operatorname{Li_{2}}\left(\sqrt{2}-1\right)+\operatorname{Li_{2}}\left(\frac{1}{2}\right). 
\end{multline*}
Consequently, we have that 
\begin{multline*}
\operatorname{Li_{2}}\left(\frac{1}{\sqrt{2}}\right)-\operatorname{Li_{2}}\left(\frac{2+\sqrt{2}}{4}\right) - 
\operatorname{Li_{2}}(2\sqrt{2}-2) + \frac{\pi^{2}}{12} - \\ \ln(4-2\sqrt{2})\ln(2)= 
-\operatorname{Li_{2}}(\sqrt{2}-1)-\frac{1}{2}\ln^{2}(2). 
\end{multline*}
By applying the reflection identity in \eqref{reflectionformula} to the third term, we find that 
$$ -\operatorname{Li_{2}}\left(2\sqrt{2}-2\right)=\operatorname{Li_{2}}\left(3-2\sqrt{2}\right)-\frac{\pi^{2}}{6}+\ln\left(2\sqrt{2}-
2\right)\ln\left(3-2\sqrt{2}\right), $$
and each side may be expressed as 
\begin{multline*}
\operatorname{Li_{2}}\left((\sqrt{2}-1)^{2}\right)-\frac{\pi^{2}}{6}+
\ln\left(2\sqrt{2}-2\right)\ln\left(3-2\sqrt{2}\right) = \\ 2\operatorname{Li_{2}}\left(\sqrt{2}-1\right)+
2\operatorname{Li_{2}}\left(1-\sqrt{2}\right)- 
\frac{\pi^{2}}{6}+\ln\left(2\sqrt{2}-2\right)\ln\left(3-\sqrt{2}\right).
\end{multline*}
On the other hand, from the relation 
$$\operatorname{Li_{2}}\left(1-\sqrt{2}\right)=\operatorname{Li_{2}}\left(\sqrt{2}-1\right)-\frac{\pi^{2}}{8}+\frac{1}{2}\ln^{2}\left(\sqrt{2}+
1\right), $$ we find that 
\begin{multline*}
-\operatorname{Li_{2}}\left(2\sqrt{2}-2\right)=4\operatorname{Li_{2}}\left(\sqrt{2}-1\right)-\frac{5\pi^{2}}{12} + \\ 
\ln^{2}\left(\sqrt{2}+1\right)
+\ln\left(2\sqrt{2}-2\right)\ln\left(3-\sqrt{2}\right).
\end{multline*}
From the relation between $\operatorname{Li_{2}}(\frac{1}{\sqrt{2}})$ and $\operatorname{Li_{2}}(\sqrt{2}-1)$, we obtain the desired result. 
 
\subsection{A new two-term dilogarithm identity involving $\phi$}\label{subsectionphimotivate}
 Using the gemini function $\gemini_{\phi}(x)$ together with the five-term identity for the gemini function we introduce a proof of
the first of the new and motivating results involving $\phi$ highlighted in Section 
\ref{subsectionmotivating}

\begin{theorem}
The closed-form evaluation 
\begin{multline*}
\operatorname{Li}_2\left(\frac{\sqrt{\phi }+1}{\phi ^2}\right)+\operatorname{Li}_2\left(\phi ^3-\phi ^{5/2}\right) = 
\frac{17 \pi ^2}{60} 
-\ln \left(\phi ^{5/2}-2 \phi \right) \ln \left(\phi ^{7/2} - 
\phi ^3\right) - \\ \frac{1}{2} \ln \left(\phi ^{7/2}-\phi ^3\right) \ln \left(\phi
^{11/2}+2 \phi ^{7/2}+\phi ^6+\phi ^4\right)-\frac{11}{8} \ln ^2(\phi ) 
\end{multline*}
holds. 
\end{theorem}

\textit{Proof.} Setting $a=\phi$, and setting the integration
limits so that $x_{1}=\ln(\sqrt{\phi})$ and $x_{2}=\ln\left(\frac{\phi+\sqrt{\phi}}{\sqrt{\phi}-1}\right)$, the
five-term identity with the substituted initial values gives us that 
\begin{multline*}
\operatorname{Li_{2}}\left(-\frac{\phi}{\sqrt{\phi}}\right)-\operatorname{Li_{2}}\left(\frac{1}{\sqrt{\phi}}\right)+\frac{\pi^{2}}{6}
-\operatorname{Li_{2}}\left(-\phi\right)-\ln\left(\sqrt{\phi}\right)\ln\left(\frac{\sqrt{\phi}+\phi}{\sqrt{\phi}-1}\right) = \\ 
-\operatorname{Li_{2}}\left(-\phi\frac{\sqrt{\phi}-1}{\sqrt{\phi}+\phi}\right)+
\operatorname{Li_{2}}\left(\frac{\sqrt{\phi}-1}{\sqrt{\phi}+\phi}\right). 
\end{multline*}
Consequently, we have that 
\begin{multline*}
\operatorname{Li_{2}}\left(-\sqrt{\phi}\right)-\operatorname{Li_{2}}\left(\frac{1}{\sqrt{\phi}}\right)+\frac{\pi^{2}}{3}
+\operatorname{Li_{2}}\left(-\frac{1}{\phi}\right)+\frac{1}{2}\ln^{2}(\phi) - \\
\ln\left(\sqrt{\phi}\right)\ln\left(\frac{\sqrt{\phi}+\phi}{\sqrt{\phi}-1}\right) = 
-\operatorname{Li_{2}}\left(-\phi\frac{\sqrt{\phi}-1}{\sqrt{\phi}+\phi}\right)+
\operatorname{Li_{2}}\left(\frac{\sqrt{\phi}-1}{\sqrt{\phi}+\phi}\right), 
\end{multline*}
 which, in turn, gives us that 
\begin{multline*}
-\operatorname{Li_{2}}\left(-\frac{1}{\sqrt{\phi}}\right)-\operatorname{Li_{2}}\left(\frac{1}{\sqrt{\phi}}\right)+\frac{\pi^{2}}{6}
+\operatorname{Li_{2}}\left(-\frac{1}{\phi}\right)+\frac{3}{8}\ln^{2}(\phi) - \\
\ln\left(\sqrt{\phi}\right)\ln\left(\frac{\sqrt{\phi}+\phi}{\sqrt{\phi}-1}\right) = 
-\operatorname{Li_{2}}\left(-\phi\frac{\sqrt{\phi}-1}{\sqrt{\phi}+\phi}\right)+
\operatorname{Li_{2}}\left(\frac{\sqrt{\phi}-1}{\sqrt{\phi}+\phi}\right). 
\end{multline*}
From the consequent relation 
\begin{multline*}
-\frac{1}{2}\operatorname{Li_{2}}\left(\frac{1}{\phi}\right)+\frac{\pi^{2}}{10}+\frac{7}{8}\ln^{2}(\phi)-
\frac{1}{2}\ln(\phi)\ln\left(\frac{\sqrt{\phi}+\phi}{\sqrt{\phi}-1}\right) = \\ 
-\operatorname{Li_{2}}\left(-\phi\frac{\sqrt{\phi}-1}{\sqrt{\phi}+\phi}\right)+
\operatorname{Li_{2}}\left(\frac{\sqrt{\phi}-1}{\sqrt{\phi}+\phi}\right),
\end{multline*}
the desired result then follows from the closed form for $\text{Li}_{2}\big( \frac{1}{\phi} \big)$. 

\subsection{Applying Theorem \ref{mainthreeterm} with the second term fixed by $\frac{1}{\phi}$}\label{secondphimotivate}
 We proceed to apply Theorem \ref{mainthreeterm} to prove the second out of the new and motivating results involving
 $\phi$ highlighted in Section \ref{subsectionmotivating}. 

\begin{theorem}
The closed-form evaluation 
\begin{multline*}
\operatorname{Li_{2}}\left(\frac{\sqrt{\phi^{7} + 3\phi^{5}}-\phi^{2}-3\phi}{2}\right)-\operatorname{Li_{2}}\left(\frac{1 + 
 \sqrt{4\phi-3}}{2\phi}\right) = \\ 
 -\frac{\pi^{2}}{10}+\ln^{2}(\phi)-\ln\left(\frac{\phi^{2}+1+\sqrt{4\phi^{3}-
3\phi^{2}}}{2}\right)\ln\left(\frac{2\phi-1+\sqrt{4\phi-3}}{2\phi}\right)
\end{multline*}
holds. 
\end{theorem}

\textit{Proof.} From the three-term relation in Theorem \ref{mainthreeterm}, we find that 
$$\operatorname{Li_{2}}\left(\frac{a-1}{a^{2}}\right)+\operatorname{Li_{2}}\left(\frac{1}{a^{2}-a+1}\right)-
\operatorname{Li_{2}}\left(\frac{a}{a^{2}-a+1}\right)+\ln\left(\frac{a}{a-1}\right)\ln\left(\frac{a^{2}}{a^{2}-a+1}\right)=0. $$
We set the argument of the second term as $\frac{1}{\phi}$, i.e., so that we obtain
$\operatorname{Li_{2}}\left(\frac{1}{\phi}\right)$. Solving for the parameter $a$, we obtain $a = \frac{1\pm\sqrt{4\phi-3}}{2}$,
giving us that 
\begin{multline*}
\operatorname{Li_{2}}\left(-\frac{\phi^{2}+3\phi}{2}+\frac{\sqrt{\phi^{7}+3\phi^{5}}}{2}\right)+\operatorname{Li_{2}}\left(\frac{1}{\phi}\right)-
\operatorname{Li_{2}}\left(\frac{1+\sqrt{4\phi-3}}{2\phi}\right) + \\ 
\ln\left(\frac{\phi^{2}+1+\sqrt{4\phi^{3}-3\phi^{2}}}{2}\right)\ln\left(\frac{2\phi-1+\sqrt{4\phi-3}}{2\phi}\right)=0, 
\end{multline*}
and this gives us an equivalent formulation of the desired result. 

\subsection{Applying Theorem \ref{mainthreeterm} with the first term fixed by $-\phi$}\label{negativephisec}
 We proceed to prove the third out of the new and motivating results from Section \ref{subsectionmotivating} involving $\phi$. 

\begin{theorem}
The closed-form evaluation 
\begin{multline*}
\operatorname{Li_{2}}\left(\frac{1}{2}\phi-\frac{1}{2}\sqrt{\frac{\phi^{2}+2}{\phi^{3}}}\right)-
\operatorname{Li_{2}}\left(\frac{1}{2\phi^{2}}-\frac{1}{2}\sqrt{\frac{\phi^{2}+2}{\phi^{3}}}\right) = \\ 
\frac{\pi^{2}}{10}+\ln^{2}(\phi)+2\ln(\phi)\ln\left(-\frac{1}{2\phi}+\frac{1}{2}\sqrt{\frac{\phi^{2}+2}{\phi}}\right) 
\end{multline*}
holds. 
\end{theorem}

\textit{Proof.} Using the three-term relation in Theorem \ref{mainthreeterm}, let the first term be such that
$$\operatorname{Li_{2}}\left(\frac{a-1}{a^{2}}\right)=\operatorname{Li_{2}}(-\phi) \Rightarrow \frac{a-1}{a^{2}}=-\phi \Rightarrow
a=\pm\frac{\sqrt{4\phi+1}-1}{2\phi}.$$ 
Now, set $a=-\frac{\sqrt{4\phi+1}+1}{2\phi}$. We thus obtain that 
\begin{multline*}
\operatorname{Li_{2}}(-\phi)+\operatorname{Li_{2}}\left(\frac{1}{2}\phi-\frac{1}{2}\sqrt{\frac{\phi^{2}+2}{\phi^{3}}}\right)
-\operatorname{Li_{2}}\left(\frac{1}{2\phi^{2}}-\frac{1}{2}\sqrt{\frac{\phi^{2}+2}{\phi^{3}}}\right) + \\ 
\ln\left(\frac{1}{\phi^{2}}\right)\ln\left(-\frac{1}{2\phi}+\sqrt{\frac{\phi^{2}+2}{\phi}}\right)=0, 
\end{multline*}
and this is equivalent to the desired result. 

\subsection{Applying Theorem \ref{mainthreeterm} with the second term fixed by $\frac{1}{\phi^{2}}$}\label{motivatefourth}

We proceed to derive another two-term identity with the aid of known value of $\operatorname{Li_{2}}\left(\frac{1}{\phi^{2}}\right)$. 
This gives us a proof of another motivating result highlighted in the Introduction. 

\begin{theorem}
The closed-form evaluation 
\begin{multline*}
\operatorname{Li_{2}}\left(\frac{1}{2}\sqrt{\phi^{2}+3\phi}-\frac{\phi^{2}+1}{2\phi}\right)-\operatorname{Li_{2}}\left(\frac{1}{2\phi^{2}}
+\frac{1}{2}\sqrt{\frac{\phi^{2}+2}{\phi^{3}}}\right)= \\
-\frac{\pi^{2}}{15}+\ln^{2}(\phi)-\ln\left(\frac{\phi^{2}}{2}+\frac{1}{2}\sqrt{\frac{\phi^{2}+2}{\phi}}\right)\ln\left(\frac{1}{2}\phi
+\frac{1}{2}\sqrt{\frac{\phi^{2}+2}{\phi^{3}}}\right) 
\end{multline*}
holds. 
\end{theorem}

\textit{Proof.} Setting $\frac{1}{\phi^{2}}$ in the second argument of Theorem \ref{mainthreeterm}, we obtain $\frac{1}{a^{2}-a+1}=\frac{1}{\phi^{2}}
\Rightarrow a=\frac{1\pm\sqrt{\phi^{2}+3\phi}}{2}$, we obtain the vanishing relation 
\begin{multline*}
\operatorname{Li_{2}}\left(-\frac{\phi^{2}+1}{2\phi}+\frac{1}{2}\sqrt{\phi^{2} + 
3\phi}\right)+\operatorname{Li_{2}}\left(\frac{1}{\phi^{2}}\right)
-\operatorname{Li_{2}}\left(\frac{1}{2\phi^{2}}+\frac{1}{2}\sqrt{\frac{\phi^{2}+2}{\phi^{3}}}\right)+ \\ 
\quad \quad \quad \quad\ln\left(\frac{\phi^{2}}{2}+\frac{1}{2}\sqrt{\frac{\phi^{2}+2}{\phi}}\right)\ln\left(\frac{1}{2}\phi+
\frac{1}{2}\sqrt{\frac{\phi^{2}+2}{\phi^{3}}}\right)=0, 
\end{multline*}
and this gives us an equivalent version of the desired result. 

\subsection{A proof of Khoi's two-term dilogarithm relation}\label{subsectionKhoi}
 As below, we formulate an alternative and elementary derivation of Khoi's two-term dilogarithm relation displayed in
\eqref{displayKhoimotivate}. This provides another solution to an open problem given by Khoi \cite{Kho14} that was previously solved,
in a different way, by Lima \cite{Lima2017}. 

\begin{theorem}
Khoi's two-term dilogarithm relation in \eqref{displayKhoimotivate} holds. 
\end{theorem}

\textit{Proof.} We apply the fixed-point identity in Theorem \ref{firstfixedpoint}, by applying the $\gemini_{-\frac{1}{\phi^{2}}}(x)$-function,
so that the shape factor is such that $a=-\frac{1}{\phi^{2}}$ and the corresponding fixed point is given by
$x=\ln(1+\sqrt{1+a})=\ln\left(1+\sqrt{1-\frac{1}{\phi^{2}}}\right)=\ln\left(1+\frac{1}{\sqrt{\phi}}\right)$. By Theorem \ref{firstfixedpoint},
we find that 
\begin{multline*}
\operatorname{Li_{2}}\left(\frac{1}{\phi^{2}}\cdot\frac{1}{1+\frac{1}{\sqrt{\phi}}}\right)-
\operatorname{Li_{2}}\left(\frac{1}{1+\frac{1}{\sqrt{\phi}}}\right) - \\ 
\frac{1}{2}\operatorname{Li_{2}}\left(\frac{1}{\phi^{2}}\right)+
\frac{\pi^{2}}{12}-\frac{1}{2}\ln^{2}\left(1+\frac{1}{\sqrt{\phi}}\right)=0. 
\end{multline*}
Consequently, we have that 
$$ \operatorname{Li_{2}}\left(1-\frac{1}{\sqrt{\phi}}\right)-\operatorname{Li_{2}}\left(\frac{1}{1+\frac{1}{\sqrt{\phi}}}\right)
-\frac{\pi^{2}}{30}+\frac{1}{2}\ln^{2}(\phi)+\frac{\pi^{2}}{12}-\frac{1}{2}\ln^{2}\left(1+\frac{1}{\sqrt{\phi}}\right)$$
vanishes, so that 
$$ \operatorname{Li_{2}}\left(1-\frac{1}{\sqrt{\phi}}\right)-\operatorname{Li_{2}}\left(\frac{1}{1+\frac{1}{\sqrt{\phi}}}\right)=
-\frac{\pi^{2}}{20}-\frac{1}{2}\ln^{2}(\phi)+\frac{1}{2}\ln^{2}\left(1+\frac{1}{\sqrt{\phi}}\right), $$
and this is equivalent to the desired result. 

\subsection{An alternative proof of Adegoke and Frontczak's formula}

We have also applied our methods to recover a two-term dilogarithm relation given recently by identity in the quite recent paper
by Adegoke and Frontczak \cite{Ade_24} and reproduced below 

\begin{theorem} 
(Adegoke and Frontczak, 2024) The two-term dilogarithm relation such that 
$$\operatorname{Li_{2}}\left(\frac{\phi}{\sqrt{5}}\right)+\operatorname{Li_{2}}\left(\frac{1}{\sqrt{5}\phi}\right) = 
\frac{\pi^{2}}{6}-\ln\left(\frac{\phi}{\sqrt{5}}\right)\ln\left(\frac{1}{\sqrt{5}\phi}\right) $$
holds \cite{Ade_24}. 
\end{theorem}

\textit{Proof.} By writing $F_{n} $ in place of the $n^{\text{th}}$ entry in the Fibonacci sequence, we have that $\phi^{n+1} = \phi F_{n+1}+F_{n}$, so that
$\frac{F_{n+1}}{\phi^{n}}=1-\frac{F_{n}}{\phi^{n+1}}$, and this can be exploited using the reflection identity for $\text{Li}_{2}$, with 
$$ \operatorname{Li_{2}}\left(\frac{F_{n+1}}{\phi^{n}}\right)=-\operatorname{Li_{2}}\left(1-\frac{F_{n+1}}{\phi^{n}}\right)+\frac{\pi^{2}}{6}-
\ln\left(\frac{F_{n+1}}{\phi^{n}}\right)\ln\left(1-\frac{F_{n+1}}{\phi^{n}}\right), $$
which implies that 
$$\operatorname{Li_{2}}\left(\frac{F_{n+1}}{\phi^{n}}\right)=-\operatorname{Li_{2}}\left(\frac{F_{n}}{\phi^{n+1}}\right)+\frac{\pi^{2}}{6}-
\ln\left(\frac{F_{n+1}}{\phi^{n}}\right)\ln\left(\frac{F_{n}}{\phi^{n+1}}\right).$$
By the Binet formula $F_{n}=\frac{\phi^{n}-(-\phi)^{n}}{2\phi-1}$, we obtain the identity such that 
\begin{multline*}
\operatorname{Li_{2}}\left(\frac{\phi^{n+1}-(-\phi)^{-n-1}}{2\phi^{n+1}-\phi^{n}}\right)+
\operatorname{Li_{2}}\left(\frac{\phi^{n}-(-\phi)^{-n}}{2\phi^{n+2}-\phi^{n+1}}\right) = \\ 
\frac{\pi^{2}}{6}-\ln\left(\frac{\phi^{n+1}-(-\phi)^{-n-1}}{2\phi^{n + 
1}-\phi^{n}}\right)\ln\left(\frac{\phi^{n}-(-\phi)^{-n}}{2\phi^{n+2}-
\phi^{n+1}}\right). 
\end{multline*}
We proceed to define the dilogarithm argument values given above as $n$ approaches infinity. These are given by 
$$\lim\limits_{n \to \infty}\frac{\phi^{n+1}-(-\phi)^{-n-1}}{2\phi^{n+1}-\phi^{n}}=\frac{\phi}{\sqrt{5}} $$ and 
$$ \lim\limits_{n \to \infty}\frac{\phi^{n}-(-\phi)^{-n}}{2\phi^{n+2}-\phi^{n+1}}=\frac{1}{\sqrt{5}\phi}.$$ 
By inserting these limiting values into the above identity derived via the Binet formula, we obtain the desired result. 

\subsection{A dilogarithm relation involving the argument $-\frac{1}{\phi^2}$}

To the below of our knowledge, the two-term dilogarithm relation given below is original. 

\begin{theorem}
The closed-form evaluation 
$$ \operatorname{Li_{2}}\left(-\frac{1}{\phi^{2}}\right)+\operatorname{Li_{2}}\left(\frac{\phi^{2}+1}{5\phi^{2}}\right) 
=-\frac{1}{2}\ln^{2}\left(\frac{\phi^{2}+1}{\phi^{2}}\right) 
$$ holds. 
\end{theorem}

\textit{Proof.} We apply our five-term gemini identity, setting $a=-\frac{1}{\phi^{2}}$, $x_{1}=\ln\left(\frac{\sqrt{5}}{\phi}\right)$ and $x_{2}
=\ln(\phi^{2}).$ This gives us that 
\begin{multline*}
\operatorname{Li_{2}}\left(\frac{1}{\phi^{2}}\cdot\frac{\phi}{\sqrt{5}}\right)-\operatorname{Li_{2}}\left(\frac{\phi}{\sqrt{5}}\right)-
\operatorname{Li_{2}}\left(\frac{1}{\phi^{2}}\right)+\frac{\pi^{2}}{6} - 
\\ \ln(\phi^{2})\ln\left(\frac{\sqrt{5}}{\phi}\right)= 
-\operatorname{Li_{2}}\left(\frac{1}{\phi^{4}}\right)+\operatorname{Li_{2}}\left(\frac{1}{\phi^{2}}\right). 
\end{multline*}
Consequently, we obtain that 
\begin{multline*}
\operatorname{Li_{2}}\left(\frac{1}{\sqrt{5}\phi}\right)-\operatorname{Li_{2}}\left(\frac{\phi}{\sqrt{5}}\right)+\frac{\pi^{2}}{6}-
2\ln(\phi)\ln\left(\frac{\sqrt{5}}{\phi}\right) = \\
-\operatorname{Li_{2}}\left(\frac{1}{\phi^{4}}\right)+
2\operatorname{Li_{2}}\left(\frac{1}{\phi^{2}}\right)=-2\operatorname{Li_{2}}\left(-\frac{1}{\phi^{2}}\right).
\end{multline*}
We proceed to rewrite the second term using the formulation of Landen's identity in \eqref{Landensformula}. Writing $\frac{\sqrt{5}}{\phi}=
\frac{\phi^{2}+1}{\phi^{2}}$, we obtain that 
$$\operatorname{Li_{2}}\left(\frac{1}{\sqrt{5}\phi}\right)=\operatorname{Li_{2}}\left(\frac{\phi^{2}+1}{5\phi^{2}}\right)$$
and that 
$$-\operatorname{Li_{2}}\left(\frac{\phi}{\sqrt{5}}\right)=-\operatorname{Li_{2}}\left(-\frac{1}{\phi^{2}}\right)-\frac{\pi^{2}}{6}+
\frac{1}{2}\ln\left(1+\frac{1}{\phi^{2}}\right)\ln\left(\phi^{4}+\phi^{2}\right).$$
Consequently, we obtain 
\begin{multline*}
\operatorname{Li_{2}}\left(\frac{\phi^{2}+1}{5\phi^{2}}\right)-\operatorname{Li_{2}}\left(-\frac{1}{\phi^{2}}\right)-\frac{\pi^{2}}{6} + 
\frac{1}{2}\ln\left(1+\frac{1}{\phi^{2}}\right)\ln(\phi^{4}+\phi^{2}) + \\
\frac{\pi^{2}}{6}-2\ln(\phi)\ln\left(\frac{\phi^{2}+1}{\phi^{2}}\right)=-2\operatorname{Li_{2}}\left(-\frac{1}{\phi^{2}}\right).
\end{multline*}
As a consequence, we have that 
\begin{multline*}
\operatorname{Li_{2}}\left(-\frac{1}{\phi^{2}}\right)+\operatorname{Li_{2}}\left(\frac{\phi^{2}+1}{5\phi^{2}}\right)+ \\ 
\frac{1}{2}\ln\left(\frac{\phi^{2}+1}{\phi^{2}}\right)\ln(\phi^{4}+\phi^{2})-2\ln(\phi)\ln\left(\frac{\phi^{2}+1}{\phi^{2}}\right)=0, 
\end{multline*}
and this is equivalent to the desired result. 

\subsection{On a modified version of \eqref{finalfivegemini}}

With regard to our proof of Theorem \ref{finalfivegemini}, by transforming the third dilogarithmic term using Landen's identity
\eqref{Landensformula}, we obtain that 
\begin{multline*}
\operatorname{Li_{2}}\left(-\frac{a}{x}\right)-\operatorname{Li_{2}}\left(\frac{1}{x}\right)+\frac{\pi^{2}}{6}-
\operatorname{Li_{2}}\left(-a\right)-\ln(x)\ln\left(\frac{x+a}{x-1}\right) + \\ 
\operatorname{Li_{2}}\left(-a \cdot \frac{x-1}{x+a}\right)-\operatorname{Li_{2}}\left(\frac{x-1}{x+a}\right)=0. 
\end{multline*}
Consequently, 
\begin{align}
\begin{split} 
\operatorname{Li_{2}}\left(-\frac{a}{x}\right)-\operatorname{Li_{2}}\left(\frac{1}{x}\right)+\frac{\pi^{2}}{3}-
\operatorname{Li_{2}}\left(\frac{1}{a+1}\right)-\frac{1}{2}\ln(a+1)\ln\left(\frac{a+1}{a^{2}}\right) - \\ 
\ln(x)\ln\left(\frac{x+a}{x-1}\right)+\operatorname{Li_{2}}\left(-a \cdot \frac{x-1}{x+a}\right)-
\operatorname{Li_{2}}\left(\frac{x-1}{x+a}\right)=0. 
\end{split}\label{splitforinversion}
\end{align}
This leads us toward the four-term dilogarithm identity highlighted below. 

\begin{theorem}\label{theorem4from5}
The relation 
\begin{multline*}
-\operatorname{Li_{2}}\left(-\frac{2a+1}{a^{2}}\right)-\operatorname{Li_{2}}\left(\frac{a}{2a+1}\right)- 
2\operatorname{Li_{2}}\left(\frac{1}{a+1}\right) + 
\operatorname{Li_{2}}\left(-\frac{a}{a+1}\right)+\frac{\pi^{2}}{6}- \\ 
\frac{1}{2}\ln^{2}\left(\frac{2a+1}{a^{2}}\right) - 
\frac{1}{2}\ln(a+1)\ln\left(\frac{a+1}{a^{2}}\right)-\ln(a+1)\ln\left(\frac{2a+1}{a}\right)=0 
\end{multline*}
holds for $a \in \mathbb{C} \setminus \{ 0, -\frac{1}{2}, -1 \}$ such that the series associated with the dilogarithmic expressions
above converge. 
\end{theorem}

\textit{Proof.} We transform the first dilogarithmic term in \eqref{splitforinversion} according to the inversion formula \eqref{inversionformula}
such that $\operatorname{Li_{2}}\left(-\frac{a}{x}\right)=-\operatorname{Li_{2}}\left(-\frac{x}{a}\right)-\frac{\pi^{2}}{6}-
\frac{1}{2}\ln^{2}(\frac{a}{x})$. We then let the argument of the third term in \eqref{splitforinversion} be equal to the argument
value of the fifth term in \eqref{splitforinversion}. Hence, the corresponding value for $x$ can be calculated according to the relation
$\frac{1}{a+1}=\frac{x-1}{x+a} \Rightarrow x=\frac{2a+1}{a}$. From our obtained value of $x$, we obtain the desired four-term identity. 

Theorem \ref{theorem4from5} allows us to obtain two-term dilogarithm identities that recalls the remarkable two-term dilogarithm due to
Khoi and reproduced in Section \ref{subsectionKhoi}. 

\begin{example}
By taking the second dilogarithmic in Theorem \ref{theorem4from5} and transforming this expression using Euler's formula in
\eqref{reflectionformula}, and by then rewriting the fourth term according to Landen's formula in \eqref{Landensformula} in such a way so that 
the resultant terms can be combined, we find that the resultant second term is given by 
$$-\operatorname{Li_{2}}\left(\frac{a}{2a+1}\right)=\operatorname{Li_{2}}\left(\frac{a+1}{2a+1}\right)-\frac{\pi^{2}}{6}+
\ln\left(\frac{a}{2a+1}\right)\ln\left(\frac{a+1}{2a+1}\right)$$ 
and that the resultant fourth term is given by 
$$\operatorname{Li_{2}}\left(-\frac{a}{a+1}\right)=\operatorname{Li_{2}}\left(\frac{a+1}{2a+1}\right)-\frac{\pi^{2}}{6}+
\frac{1}{2}\ln\left(\frac{2a+1}{a+1}\right)\ln\left(\frac{2a^{2}+3a+1}{a^{2}}\right).$$
By then setting $a = \phi$, we obtain the closed-form evaluation such that 
\begin{multline*}
2\operatorname{Li_{2}}\left(\frac{1}{\phi\sqrt{5}}\right) - 
 \operatorname{Li_{2}}\left(\frac{\sqrt{5}}{\phi^{2}}\right)+\frac{\pi^{2}}{30} + 
 \frac{1}{8}\ln^{2}(5) + \\ 
\frac{1}{2}\ln(\phi)\ln\left(\frac{125}{\phi^{7}}\right)+\ln(\phi^{2}+1)\ln\left(\frac{\sqrt{\phi^{2}+1}}{\phi^{2}}\right) = 0. 
\end{multline*}
\end{example}

\section{Further identities obtained via gemini identities}
 This section deals with mathematical constants that we apply to obtain dilogarithm evaluations and ladders and related propertities 
 of gemini functions. 

\subsection{On Ramanujan's two-term dilogarithm identities}
 Ramanujan discovered a number of remarkable evaluations for two-term combinations of dilogarithms with rational arguments
 \cite[pp.\ 323--326]{Berndt1994}. These discoveries due to Ramanujan may be proved, through our gemini function-based methods, 
 as below. 

\begin{theorem}
Ramanujan's evaluation 
$$\operatorname{Li_{2}}\left(-\frac{1}{2} \right) + 
\frac{1}{6}\operatorname{Li_{2}}\left(\frac{1}{9}\right)=-\frac{\pi^{2}}{18}+\ln(2)\ln(3)-\frac{1}{2}\ln^{2}(2) - 
 \frac{1}{3}\ln^{2}(3)$$ holds. 
\end{theorem}

We select the initial values so that $a=-\frac{1}{3}$ and $x_{1}=\ln(\frac{4}{3}) \Rightarrow x_{2}=\ln(3)$, giving us that 
\begin{multline*}
 \operatorname{Li_{2}}\left( \frac{1}{3}\cdot\frac{3}{4} \right) - \operatorname{Li_{2}}\left( \frac{3}{4} \right) + \frac{\pi^{2}}{6} - 
 \operatorname{Li_{2}}\left( \frac{1}{3} \right) - \ln\left( 
 \frac{4}{3} \right)\ln(3) = \\ 
-\operatorname{Li_{2}}\left(
 \frac{1}{3}\cdot\frac{1}{3} \right) + 
\operatorname{Li_{2}}\left(
 \frac{1}{3} \right), 
\end{multline*}
 so that 
$$\operatorname{Li_{2}}\left( \frac{1}{4} \right) - 
 \operatorname{Li_{2}}\left(
 \frac{3}{4} \right) + 
 \frac{\pi^{2}}{6}-\ln\left(
 \frac{4}{3} \right) 
 \ln(3)=-\operatorname{Li_{2}}\left(
 \frac{1}{9} \right) 
 +2\operatorname{Li_{2}}\left(
 \frac{1}{3} \right), $$ 
which, in turn, gives us that 
$$2\operatorname{Li_{2}}\left(\frac{1}{4}\right)+\ln^{2}\left(\frac{4}{3}\right)=-\operatorname{Li_{2}}\left(\frac{1}{9}\right)
+2\operatorname{Li_{2}}\left(\frac{1}{3}\right). $$
 Continuing in a similar fashion, the relation
 $$2\operatorname{Li_{2}}\left(\frac{1}{4}\right)=-\operatorname{Li_{2}}\left(\frac{1}{9}\right)+2\operatorname{Li_{2}}(-2)+\frac{\pi^{2}}{3}-
\ln\left(\frac{4}{3}\right)\ln\left(\frac{9}{4}\right)$$ implies that
$$ 4\operatorname{Li_{2}}\left(\frac{1}{2}\right)+4\operatorname{Li_{2}}\left(-\frac{1}{2}\right)=-\operatorname{Li_{2}}\left(\frac{1}{9}\right)
-2\operatorname{Li_{2}}\left(-\frac{1}{2}\right)-
\ln^{2}(2)+\ln\left(\frac{4}{3}\right)\ln\left(\frac{9}{4}\right), $$
 and thus  $$6\operatorname{Li_{2}}\left(-\frac{1}{2}\right)+\operatorname{Li_{2}}\left(\frac{1}{9}\right)=-\frac{\pi^{2}}{3}+\ln^{2}(2) + 
 \ln\left(\frac{4}{3}\right)\ln\left(\frac{9}{4}\right),  $$ and this implies that 
 $$\operatorname{Li_{2}}\left(-\frac{1}{2}\right)+\frac{1}{6}\operatorname{Li_{2}}\left(\frac{1}{9}\right)=-\frac{\pi^{2}}{18}+\ln(2)\ln(3)-
\frac{1}{2}\ln^{2}(2)-\frac{1}{3}\ln^{2}(3), $$
 as desired. 

 The first and the last equations are exactly the same as expected. Below are two other Ramanujan's identities, which can also be proved
with the aid of the five-term gemini-identity. The initial values for the first identity must be as follows, $a=-\frac{2}{3}$, $x_{1}=\ln(2)$ and
$x_{2}=\ln(\frac{4}{3})$. The initials for the lower one must be in such a way that $a=+\frac{9}{8}$, $x_{1}=\ln(\frac{81}{64})$ and $x_{2}=\ln(9)$, and we obtain that 
 $$\operatorname{Li_{2}}\left(\frac{1}{4} \right) + 
 \frac{1}{3}\operatorname{Li_{2}}\left(\frac{1}{9} \right) = 
 \frac{\pi^{2}}{18}-2\ln^{2}(2)+\ln(2)\ln(3)-\frac{2}{3}\ln^{2}(3) $$
 and that 
 $$\operatorname{Li_{2}}\left(-\frac{1}{8} \right) + 
 \operatorname{Li_{2}}\left(\frac{1}{9}\right) = 
 -\frac{1}{2}\ln^{2}\left(\frac{9}{8} \right). $$

 We can derive a respective two-term identity by applying the fixed-point identity shown in Theorem \ref{firstfixedpoint}. This is a special
case related to gemini functions. By setting $a=+3$ and $x_{0}=x_{1}=x_{2}=\ln(3)$, we get a very simple two-term identity. The shape
factor and the fixed point have a curious connection, i.e., $x_{0}=\ln(1+\sqrt{1+a})=\ln(1+\sqrt{1+3})=\ln(3)$. The number 3 is also
the so-called 0-addinacci constant, which will be discussed later in this publication. We thus obtain the  vanishing relation 
\begin{multline*}
  \operatorname{Li_2}\left(
  -\frac{3}{1+\sqrt{1+3}}\right) - 
 \operatorname{Li_2}\left(\frac{1}{1+\sqrt{1+3}}\right) -  \\ 
\frac{1}{2}\operatorname{Li_{2}}(-3)+\frac{\pi^{2}}{12}-\frac{1}{2}\ln^{2}(1+\sqrt{1+3})=0, 
\end{multline*}
 and this implies that  
 $$\operatorname{Li_2(-1)}-\operatorname{Li_2}\left(\frac{1}{3}\right)-\frac{1}{2}\operatorname{Li_{2}}(-3)+\frac{\pi^{2}}{12}-
\frac{1}{2}\ln^{2}(3)=0,$$
 which, in  turn, implies that 
 $$\operatorname{Li_2}\left(\frac{1}{3}\right)+\frac{1}{2}\operatorname{Li_{2}}(-3)+\frac{1}{2}\ln^{2}(3)=0.$$ 

\subsection{Identities related to $\operatorname{Li_{2}}\left(\frac{1}{\phi^{3}}\right)$}
 Previously, we have already drawn the connection between $\operatorname{Li_{2}}(\frac{1}{\phi^{3}})$ and
$\operatorname{Li_{2}}(-\frac{1}{\phi^{3}})$. The value of $\operatorname{Li_{2}}(\frac{1}{\phi^{3}})$ trivially connects
several others dilogarithm values to arguments including the golden ratio $\phi$. How these interconnected terms behave together is discussed
next. However, the arithmetic properties of the golden ratio enables an easy formulation of identities derived next. In other words, these
identities can be derived trivially, but then a few successive transformations have to be made. Hence, we call these identities semi-trivial,
since they cannot be directly derived for this purpose just by substituting appropriate values in basic identities. 
 We begin with 
 $$\operatorname{Li_{2}}\left( 
 \frac{1}{\phi^{3}} \right) = 
 -\operatorname{Li_{2}}\left( 
 \frac{2}{\phi^{2}} \right) + 
 \frac{\pi^{2}}{6}-
\ln\left( 
 \frac{1}{\phi^{3}} \right) 
 \ln\left(
 \frac{2}{\phi^{2}}\right) $$
 and 
 $$\operatorname{Li_{2}}\left( 
 \frac{1}{\phi^{3}} \right) = 
 \operatorname{Li_{2}}(-2\phi) + 
 \frac{\pi^{2}}{6}-
\frac{3}{2}\ln(\phi)\ln\left(\frac{\phi}{4}\right), $$
 and these results follow in a direct way
 from Euler's reflection formula in 
 \eqref{reflectionformula} and the Landen identity in 
\eqref{Landensformula}. 
 Combining the above identities, we find that 
 $$\operatorname{Li_{2}}\left(\frac{2}{\phi^{2}}\right)+\operatorname{Li_{2}}\left(-2\phi\right)=-\frac{9}{2}\ln^{2}(\phi). $$
 Next, we obtain a relation involving 
 $\operatorname{Li_{2}}(\frac{1}{\phi^{3}})$ and
$\operatorname{Li_{2}}(\frac{\phi}{2})$. In this direction, 
 we apply $\gemini_{-\frac{1}{\phi}}(x)$-function. Let us set the
integration limits in such a way that $x_{1}=\ln(\frac{2}{\phi})$ and $x_{2}=\ln(\phi^{2})$. Hence, we can write
\begin{multline*}
 \operatorname{Li_{2}}\left(
 \frac{1}{\phi}\cdot\frac{\phi}{2} \right) - 
 \operatorname{Li_{2}}\left(
 \frac{\phi}{2} \right) + 
 \frac{\pi^{2}}{6}
 - \operatorname{Li_{2}}\left(
 \frac{1}{\phi} \right) 
 - \ln(\phi^{2})\ln\left(
 \frac{2}{\phi} \right) = \\ 
 -\operatorname{Li_{2}}\left(
 \frac{1}{\phi}\cdot\frac{1}{\phi^{2}} \right) + 
 \operatorname{Li_{2}}\left( 
 \frac{1}{\phi^{2}} \right), 
\end{multline*} 
 and this implies that 
\begin{multline*}
 \operatorname{Li_{2}}\left(
  \frac{1}{2} \right) - 
  \operatorname{Li_{2}}\left(
  \frac{\phi}{2} \right) + 
  \frac{\pi^{2}}{6}-\frac{\pi^{2}}{10}+\ln^{2}(\phi)
-\ln(\phi^{2})\ln\left(
  \frac{2}{\phi} \right)   = \\ 
  -\operatorname{Li_{2}}\left(
  \frac{1}{\phi^{3}} \right) + 
  \frac{\pi^{2}}{15}-\ln^{2}(\phi), 
\end{multline*}
 so that  $$\operatorname{Li_{2}}\left(\frac{1}{\phi^{3}}\right)-\operatorname{Li_{2}}\left(\frac{\phi}{2}\right)=
-\frac{\pi^{2}}{12}+\frac{1}{2}\ln^{2}(2)+2\ln(2)\ln(\phi)-4\ln^{2}(\phi). $$
Next, we derive a trivial two-term identity between $\operatorname{Li_{2}}(\frac{\phi}{2})$ and
$\operatorname{Li_{2}}(\frac{1}{2\phi^{2}})$ with the aid of a reflection identity \eqref{reflectionformula}, 
 with $$\operatorname{Li_{2}}\left(
 \frac{\phi}{2} \right) = 
  -\operatorname{Li_{2}}\left(
  \frac{1}{2\phi^{2}} \right) + 
  \frac{\pi^{2}}{6}-
\ln(\frac{\phi}{2})\ln\left(
  \frac{1}{2\phi^{2}} \right).$$ 

Now, we know that both $\operatorname{Li_{2}}(\frac{1}{2\phi^{2}})$ and $\operatorname{Li_{2}}(\frac{2}{\phi^{2}})$ are related
to $\operatorname{Li_{2}}(\frac{1}{\phi^{3}})$. Next, we derive their mutual connection. Until now, we have always fixed one of the
integration limits and the shape factor first, which we have used to determine the other integration limit. In this case, we define the shape
factor with the help of integration limits by making a good guess. Let $x_{1}=\ln(\frac{\phi^{2}}{2})$ and $x_{2}=\ln(\phi^{2})$. Hence,
we can calculate the shape factor $a$ as follows, $\frac{x_{2}+a}{x_{2}-1}=x_{1} \Rightarrow \frac{\phi^{2}+a}{\phi^{2}-1}=
\frac{1}{2}\phi^{2} \Rightarrow a=\frac{1}{2}\phi^{3}-\phi^{2}=-\frac{1}{2}$. Now, we can apply
 the five-term identity in 
 \eqref{finalfivegemini}
 using the $\gemini_{-\frac{1}{2}}(x)$-function, so that 
\begin{multline*}
 \operatorname{Li_{2}}\left(
 \frac{1}{2}\cdot\frac{2}{\phi^{2}} \right) - 
 \operatorname{Li_{2}}\left(
 \frac{2}{\phi^{2}} \right) + 
 \frac{\pi^{2}}{6}
-\operatorname{Li_{2}}\left(
 \frac{1}{2} \right) - 
 2\ln(\phi)\ln\left(
 \frac{\phi^{2}}{2} \right) = \\ 
-\operatorname{Li_{2}}\left(
 \frac{1}{2}\cdot\frac{1}{\phi^{2}} \right) + 
 \operatorname{Li_{2}}\left(
 \frac{1}{\phi^{2}} \right), 
\end{multline*}
 which implies that 
\begin{multline*}
 \operatorname{Li_{2}}\left(
 \frac{1}{\phi^{2}} \right) - 
 \operatorname{Li_{2}}\left(
 \frac{2}{\phi^{2}} \right) + 
 \frac{\pi^{2}}{6}
-\frac{\pi^{2}}{12}+\frac{1}{2}\ln^{2}(2)-2\ln(\phi)\ln\left(
 \frac{\phi^{2}}{2} \right) = \\ 
-\operatorname{Li_{2}}\left(
 \frac{1}{2\phi^{2}} \right) + 
 \operatorname{Li_{2}}\left(
 \frac{1}{\phi^{2}} \right), 
\end{multline*}
 which, in turn, implies that 
$$ \operatorname{Li_{2}}\left(\frac{1}{2\phi^{2}}\right)-\operatorname{Li_{2}}\left(\frac{2}{\phi^{2}}\right) = 
 -\frac{\pi^{2}}{12}-\frac{1}{2}\ln^{2}(2)+
 2\ln(\phi)\ln\left(\frac{\phi^{2}}{2}\right). $$
 Next, we apply the $\gemini_{\phi}(x)$-function with the integration limits $x_{1}=\ln(\phi)$ and $x_{2} = 
 \ln(2\phi^{2})$. We get the following
five-term gemini-identity \eqref{finalfivegemini}, as shown below.

\vskip 0.1in

$\operatorname{Li_{2}}(-\frac{\phi}{\phi})-\operatorname{Li_{2}}(\frac{1}{\phi})+\frac{\pi^{2}}{6}
-\operatorname{Li_{2}}(-\phi)-\ln(\phi)\ln(2\phi^{2})=-\operatorname{Li_{2}}(-\frac{\phi}{2\phi^{2}})
+\operatorname{Li_{2}}(\frac{1}{2\phi^{2}}) \Rightarrow$
\vskip 0.01in
$-\frac{\pi^{2}}{12}-\frac{\pi^{2}}{10}+\ln^{2}(\phi)+\frac{\pi^{2}}{10}+\ln^{2}(\phi)+\frac{\pi^{2}}{6}-\ln(\phi)\ln(2\phi^{2})=
-\operatorname{Li_{2}}(-\frac{1}{2\phi})+\operatorname{Li_{2}}(\frac{1}{2\phi^{2}}) \Rightarrow$

\begin{eqnarray} 
 \operatorname{Li_{2}}\left(\frac{1}{2\phi^{2}}\right)-\operatorname{Li_{2}}\left(-\frac{1}{2\phi}\right)=\frac{\pi^{2}}{12}-\ln(\phi)\ln(2)
\end{eqnarray}

The list below includes all the known simplest real valued connections of the $\operatorname{Li_{2}}(\frac{1}{\phi^{3}})$-term,
which can be derived with the aid of gemini-identities. The last one of the identities is derived with the aid of
Theorem \ref{firstfixedpoint} in such a way that the fixed point $x_{0}=\ln(2\phi^{2})$ and the respective shape factor
$a=(2\phi^{2}-1)^{2}-1=(\phi^{3})^{2}-1=(\phi^{6}-1)=(\phi^{3}+1)(\phi^{3}-1)=2\phi^{2} \cdot 2\phi=4\phi^{3}$.

\vskip 0.1in

$\operatorname{Li_{2}}(\frac{1}{\phi^{3}})=\operatorname{Li_{2}}(-\frac{1}{\phi^{3}})+\frac{\pi^{2}}{12}
-\frac{3}{2}\ln^{2}(\phi)$
\vskip 0.01in
$\operatorname{Li_{2}}(\frac{1}{\phi^{3}})=-\operatorname{Li_{2}}(-\phi^{3})-\frac{\pi^{2}}{12}-6\ln^{2}(\phi)$
\vskip 0.01in
$\operatorname{Li_{2}}(\frac{1}{\phi^{3}})=\frac{1}{4}\operatorname{Li_{2}}(\frac{1}{\phi^{6}})+\frac{\pi^{2}}{24}-\frac{3}{4}\ln^{2}(\phi)$
\vskip 0.01in
$\operatorname{Li_{2}}(\frac{1}{\phi^{3}})=\operatorname{Li_{2}}(\frac{\phi}{2})-\frac{\pi^{2}}{12}
+\frac{1}{2}\ln^{2}(2)+2\ln(2)\ln(\phi)-4\ln^{2}(\phi) \quad$ (ID4)
\vskip 0.01in
$\operatorname{Li_{2}}(\frac{1}{\phi^{3}})=\operatorname{Li_{2}}(-2\phi)+\frac{\pi^{2}}{6}-\frac{3}{2}\ln(\phi)\ln(\frac{\phi}{4})$
\vskip 0.01in
$\operatorname{Li_{2}}(\frac{1}{\phi^{3}})=-\operatorname{Li_{2}}(-\frac{1}{2\phi})-\frac{1}{2}\ln^{2}(2\phi)
-\frac{3}{2}\ln(\phi)\ln(\frac{\phi}{4}) \quad$ (ID6)
\vskip 0.01in
$\operatorname{Li_{2}}(\frac{1}{\phi^{3}})=-\operatorname{Li_{2}}(\frac{2}{\phi^{2}})+\frac{\pi^{2}}{6}
-\ln(\frac{1}{\phi^{3}})\ln(\frac{2}{\phi^{2}})$
\vskip 0.01in
$\operatorname{Li_{2}}(\frac{1}{\phi^{3}})=-\operatorname{Li_{2}}(\frac{1}{2\phi^{2}})+\frac{\pi^{2}}{12}-2\ln^{2}(2\phi)
+5\ln(2)\ln(2\phi)-\frac{7}{2}\ln^{2}(2)$
\vskip 0.01in
$\operatorname{Li_{2}}(\frac{1}{\phi^{3}})=-\frac{1}{4}\operatorname{Li_{2}}(\frac{4}{\phi^{3}})+\frac{\pi^{2}}{12}
-\frac{3}{4}\ln^{2}(\phi)-\frac{3}{2}\ln(\phi)\ln(\frac{1}{4}\phi^{3})$
\vskip 0.01in
$\operatorname{Li_{2}}(\frac{1}{\phi^{3}})=\frac{1}{4}\operatorname{Li_{2}}(-4\phi^{3})+\frac{\pi^{2}}{12}
-\frac{3}{4}\ln^{2}(\phi)+3\ln(2)\ln(\phi)$

\vskip 0.1in

By combining two of the above identities marked as (ID4) and (ID6), we can formulate the following two-term identity, which is
also introduced in the paper of Adegoke and Frontczak \cite{Ade_24}.

\begin{eqnarray}
 \operatorname{Li_{2}}\left(\frac{1}{2}\phi\right)+\operatorname{Li_{2}}\left(-\frac{1}{2\phi}\right)=
 \frac{\pi^{2}}{12}+2\ln^{2}(\phi)-\ln^{2}(2)
\end{eqnarray}

\subsection{Dilogarithmic relations involving the plastic constant}

In this Section, two-term identities and one ladder are derived related to the plastic constant $P$. First, we apply the
three-term identity given in Theorem \ref{threetermconstant}, to derive relations involving $\operatorname{Li_{2}}\left(\frac{1}{P}\right)$.
The three-term cancellation identity we obtain from Theorem \ref{threetermconstant} is such that

\vskip 0.1in

$\operatorname{Li_{2}}\left(\frac{a}{(a+1)^{2}}\right)+\operatorname{Li_{2}}\left(\frac{1}{a^{2}+a+1}\right)
-\operatorname{Li_{2}}\left(\frac{a+1}{a^{2}+a+1}\right)-\ln\left(\frac{a+1}{a}\right)\ln\left(\frac{(a+1)^{2}}{a^{2}+a+1}\right)=0.$

\vskip 0.1in

First, we formulate an equation in such a way that the arguments of the first and the second term become equal. Hence, we can write

\vskip 0.1in

$\frac{a}{(a+1)^{2}}=\frac{1}{a^{2}+a+1}\Rightarrow a^{3}-a-1=0 \Rightarrow a=\frac{\sqrt[3]{9+\sqrt{69}}+\sqrt[3]{9-\sqrt{69}}}{\sqrt[3]{18}}
\Rightarrow a\approx1.324718, $

\vskip 0.1in

which is the only real root. This root $a$ is also known as the \emph{plastic constant} $P$, with reference to the decimal expansion
for this constant given in the On-Line Encyclopedia of Integer Sequences (OEIS) \cite{oeis} as sequence {\tt A060006}. 
This constant $P$ also satisfies $P-1=\frac{1}{P^{4}}$ and $P^{2}-2=-\frac{1}{P^{5}}$. Next, we insert $P$ into the three-term identity
in Theorem \ref{threetermconstant}. By combining two first terms, whose argument values are the same, we get

\vskip 0.1in

$ 2\operatorname{Li_{2}}(\frac{P}{(P+1)^{2}})-\operatorname{Li_{2}}(\frac{P+1}{P^{2}+P+1})+\ln(\frac{P+1}{P})\ln(\frac{(P+1)^{2}}{P^{2}+P+1})=0
\Rightarrow$
\vskip 0.01in
$ 2\operatorname{Li_{2}}(\frac{1}{P^{5}})-\operatorname{Li_{2}}(\frac{P+1}{P^{2}+P+1})+\ln(\frac{P+1}{P})\ln(\frac{(P+1)^{2}}{P^{2}+P+1})=0
\Rightarrow$

\vskip 0.01in

$ 2\operatorname{Li_{2}}(\frac{p-1}{P})-\operatorname{Li_{2}}(\frac{P^{3}}{P^{2}+P^{3}})+\ln(\frac{P^{3}}{P})\ln(\frac{P^{6}}{P^{2}+P^{3}})=0
\Rightarrow$

\vskip 0.01in

$-2\operatorname{Li_{2}}(\frac{1}{P})+\frac{\pi^{2}}{3}-2\ln(P)\ln(P^{5})-\operatorname{Li_{2}}(\frac{1}{1+\frac{1}{P}})
+2\ln(P)\ln(\frac{P^{4}}{1+P})=0 \Rightarrow$
\vskip 0.01in
$-2\operatorname{Li_{2}}(\frac{1}{P})+\frac{\pi^{2}}{3}-10\ln^{2}(P)-\operatorname{Li_{2}}(\frac{1}{1+\frac{1}{P}})+2\ln^{2}(P)=0 \Rightarrow$
\vskip 0.01in
$-2\operatorname{Li_{2}}(\frac{1}{P})+\frac{\pi^{2}}{3}-8\ln^{2}(P)-\operatorname{Li_{2}}(\frac{1}{1+\frac{1}{P}})=0 \Rightarrow$
\vskip 0.01in
$-2\operatorname{Li_{2}}(\frac{1}{P})+\frac{\pi^{2}}{3}-8\ln^{2}(P)-\operatorname{Li_{2}}(-\frac{1}{P})-\frac{\pi^{2}}{6}
+\frac{1}{2}\ln(1+\frac{1}{P})\ln((1+\frac{1}{P})P^{2})=0 \Rightarrow$
\vskip 0.01in
$-2\operatorname{Li_{2}}(\frac{1}{P})+\frac{\pi^{2}}{6}-4\ln^{2}(P)-\operatorname{Li_{2}}(-\frac{1}{P})=0 \Rightarrow$

\begin{eqnarray}
\operatorname{Li_{2}}\left(\frac{1}{P}\right)=-\frac{1}{2}\operatorname{Li_{2}}\left(-\frac{1}{P}\right)+\frac{\pi^{2}}{12}-2\ln^{2}(P).
\end{eqnarray}

This result can also be called a semi-trivial, because it can be derived solely from the algebraic properties of the constant $P$
by means of basic identities. We can also derive a similar kind of two term identity, which connects $\operatorname{Li_{2}}(\frac{1}{P})$
and $\operatorname{Li_{2}}(\frac{1}{P^{2}})$ by applying the five-term gemini-identity in \eqref{finalfivegemini}. Let us choose the
integration limits in such a way that $x_{1}=\ln(P)$ and $x_{2}=\ln(P^{3})$. The next task is to calculate the corresponding shape factor
for the respective gemini function. It is given by

\vskip 0.1in

${\frac{P^{3}+a}{P^{3}-1}=P} \Rightarrow \frac{P^{3}+a}{P}=P \Rightarrow P^{3}+a=P^{2} \Rightarrow a=P^{2}-P^{3}=P^{2}(1-P)=-\frac{1}{P^{2}}$.

\vskip 0.1in

Now we can build a five-term identity with these initial values.

\vskip 0.1in

{$\operatorname{Li_{2}}(\frac{1}{P^{2}} \cdot 
\frac{1}{P})-\operatorname{Li_{2}}(\frac{1}{P})-\operatorname{Li_{2}}(\frac{1}{P^{2}})+\frac{\pi^{2}}{6}
-\ln(P)\ln(P^{3})=\operatorname{-Li_{2}}(\frac{1}{P^{2}}\cdot\frac{1}{P^{3}}) + 
\operatorname{Li_{2}}(\frac{1}{P^{3}}) \Rightarrow$}
\vskip 0.01in
${\operatorname{Li_{2}}(\frac{1}{P^{3}})-\operatorname{Li_{2}}(\frac{1}{P})-\operatorname{Li_{2}}(\frac{1}{P^{2}})+\frac{\pi^{2}}{6}-3\ln^{2}(P)=
-\operatorname{Li_{2}}(\frac{1}{P^{5}})+\operatorname{Li_{2}}(\frac{1}{P^{3}}) \Rightarrow}$
\vskip 0.01in
${-\operatorname{Li_{2}}(\frac{1}{P})-\operatorname{Li_{2}}(\frac{1}{P^{2}})+\frac{\pi^{2}}{6}-3\ln^{2}(P)=
-\operatorname{Li_{2}}(\frac{P-1}{P}) \Rightarrow}$
\vskip 0.01in
${-\operatorname{Li_{2}}(\frac{1}{P})-\operatorname{Li_{2}}(\frac{1}{P^{2}})+\frac{\pi^{2}}{6}-3\ln^{2}(P)=
\operatorname{Li_{2}}(\frac{1}{P})-\frac{\pi^{2}}{6}+\ln(\frac{1}{P})\ln(\frac{1}{P^{5}}) \Rightarrow}$
\vskip 0.01in
$-2\operatorname{Li_{2}}(\frac{1}{P})-\operatorname{Li_{2}}(\frac{1}{P^{2}})-3\ln^{2}(P)=-\frac{\pi^{2}}{3}+5\ln^{2}(P) \Rightarrow$

\begin{eqnarray}
\operatorname{Li_{2}}\left(\frac{1}{P}\right)=-\frac{1}{2}\operatorname{Li_{2}}\left(\frac{1}{P^{2}}\right)+\frac{\pi^{2}}{6}-4\ln^{2}(P)
\end{eqnarray}

Below is a list of all other two-term identities related to the plastic constant P, which can be trivially derived.

\vskip 0.1in

$\operatorname{Li_{2}}(P-1)=\operatorname{Li_{2}}(\frac{1}{P^{4}})$
\vskip 0.01in
$\operatorname{Li_{2}}(P^{2}-2)=\operatorname{Li_{2}}(-\frac{1}{P^{5}})$
\vskip 0.01in
$\operatorname{Li_{2}}(\frac{1}{P})=\frac{1}{2}\operatorname{Li_{2}}(\frac{1}{P^{3}})+\frac{\pi^{2}}{12}-\ln^{2}(P)$
\vskip 0.01in
$\operatorname{Li_{2}}(\frac{1}{P})=-\operatorname{Li_{2}}(\frac{1}{P^{5}})+\frac{\pi^{2}}{6}-5\ln^{2}(P)$
\vskip 0.01in
$\operatorname{Li_{2}}(\frac{1}{P})=\operatorname{Li_{2}}(-\frac{1}{P^{4}})+\frac{\pi^{2}}{6}-\frac{9}{2}\ln^{2}(P)$

\vskip 0.1in

From the above identities, we can also write the following formula, as shown below.

\vskip 0.1in

$\operatorname{Li_{2}}(\frac{1}{P^{5}})+\operatorname{Li_{2}}(-\frac{1}{P^{4}})=-\frac{1}{2}\ln^{2}(P)$

\vskip 0.1in

Numerous valid dilogarithm ladder relations can be derived for the plastic constant $P$, as in the following result. 
Results of a similar nature were recently introduced by Hakimoglu-Brown \cite{Hak_25}. A detailed study on the generation of
dilogarithm ladders can be found in Lewin's monograph on structural properties of polylogarithms \cite{Lewin1991Structural}.

\begin{theorem}
The valid ladder relation such that 
$$ 2\operatorname{Li_{2}}\left(\frac{1}{P^{3}}\right)+2\operatorname{Li_{2}}\left(\frac{1}{P^{4}}\right) + 
2\operatorname{Li_{2}}\left(\frac{1}{P^{5}}\right)-\operatorname{Li_{2}}\left(\frac{1}{P^{8}}\right)-\frac{\pi^{2}}{3}+15\ln^{2}(P) = 0$$
holds. 
\end{theorem}

\textit{Proof.} Our derivation is based on the five-term gemini identity, for the case whereby $a=+P^{4}$, letting the integration limits be
such that 
$x_{1}=\ln(P^{3})$ and $$x_{2}=\ln\left(\frac{P^{3}+P^{4}}{P^{3}-1}\right)=\ln\left(\frac{P^{3}(1+P)}{P}\right)=\ln(P^{2} \cdot P^{3})=\ln(P^{5}).$$
From the specified values, the five-term gemini identity gives us that 
\begin{multline*}
\operatorname{Li_{2}}\left(-\frac{P^{4}}{P^{3}}\right) - 
\operatorname{Li_{2}}\left(\frac{1}{P^{3}}\right) - 
\operatorname{Li_{2}}(-P^{4})+\frac{\pi^{2}}{6} - \\ \ln(P^{3})\ln(P^{5}) = 
-\operatorname{Li_{2}}\left(-\frac{P^{4}}{P^{5}}\right) + \operatorname{Li_{2}}\left(\frac{1}{P^{5}}\right). 
\end{multline*}
Consequently, we have that 
\begin{multline*} 
\operatorname{Li_{2}}(-P)-\operatorname{Li_{2}}\left(\frac{1}{P^{3}}\right) - 
\operatorname{Li_{2}}(-P^{4})+\frac{\pi^{2}}{6} - \\ 15\ln(P)\ln(P) = 
-\operatorname{Li_{2}}\left(-\frac{1}{P}\right) + \operatorname{Li_{2}}\left(\frac{1}{P^{5}}\right). 
\end{multline*}
As a consequence, we obtain that 
\begin{multline*}
-\operatorname{Li_{2}}\left(-\frac{1}{P} \right) - 
\frac{1}{2}\ln^{2}(P)-\operatorname{Li_{2}}\left(\frac{1}{P^{3}}\right) + 
\operatorname{Li_{2}}\left(-\frac{1}{P^{4}} \right)
+ \\ \frac{\pi^{2}}{6}-7\ln^{2}(P)=-\operatorname{Li_{2}}\left(-\frac{1}{P}\right) + 
\operatorname{Li_{2}}\left(\frac{1}{P^{5}}\right). 
\end{multline*}
Finally, we obtain that 
$$ -\operatorname{Li_{2}}\left(\frac{1}{P^{3}} \right)+\operatorname{Li_{2}}\left(-\frac{1}{P^{4}}\right)+ 
\frac{\pi^{2}}{6}-\frac{15}{2}\ln^{2}(P)=\operatorname{Li_{2}}\left(\frac{1}{P^{5}} \right), $$
giving us an equivalent version of the desired result. 

\subsection{The super golden ratio $\operatorname{\Psi}$ and a dilogarithm}

The \emph{super golden ratio} $\Psi$ is a constant that, informally, may be seen as especially well suited for dilogarithm identities. 
The decimal expansion for this constant is indexed in the OEIS as {\tt A092526}. It is the positive root of a cubic trinomial, which 
is given by $$ x^{3}-x^{2}-1=0 \Rightarrow x=\Psi=\frac{1}{3}\left(1+\sqrt[3]{\frac{29}{2}-\frac{3\sqrt{93}}{2}}+\sqrt[3]{\frac{29}{2} + 
\frac{3\sqrt{93}}{2}}\right). $$ Next, we derive a ladder by using the fixed-point identity displayed in Theorem \ref{firstfixedpoint}.
First, we define the relationship between the fixed point and the shape factor in such a way that the $x_{0}=\ln(x)$ and the shape factor
$a$ depend on each other according to $$x^{2}=1+\sqrt{1+\frac{1}{x^{3}}} \Rightarrow x^{7}-2x^{5}+1=0 \Rightarrow x=-P, +1 \vee \Psi.$$
Let us choose $\Psi$. Hence, the fixed point is such that $x_{0}=\ln(\Psi^{2})$ and the shape factor $a=+\frac{1}{\Psi^{3}}$. The three-term
fixed point identity is given by

$\operatorname{Li_{2}}(-\frac{1}{\Psi^{5}})-\operatorname{Li_{2}}(\frac{1}{\Psi^{2}})+\frac{\pi^{2}}{12}-\frac{1}{2}\ln^{2}(\Psi^{2})
-\frac{1}{2}\operatorname{Li_{2}}(-\frac{1}{\Psi^{3}})=0 \Rightarrow$

\vskip 0.02in

$\frac{1}{2}\operatorname{Li_{2}}(\frac{1}{\Psi^{10}})-\operatorname{Li_{2}}(\frac{1}{\Psi^{5}})-\operatorname{Li_{2}}(\frac{1}{\Psi^{2}})+
\frac{\pi^{2}}{12}-2\ln^{2}(\Psi)-\frac{1}{4}\operatorname{Li_{2}}(\frac{1}{\Psi^{6}})+\frac{1}{2}\operatorname{Li_{2}}(\frac{1}{\Psi^{3}})=0
\Rightarrow$

\begin{eqnarray}
4\operatorname{Li_{2}}\left(\frac{1}{\Psi^{2}}\right)-2\operatorname{Li_{2}}\left(\frac{1}{\Psi^{3}}\right)+
4\operatorname{Li_{2}}\left(\frac{1}{\Psi^{5}}\right)+\operatorname{Li_{2}}\left(\frac{1}{\Psi^{6}}\right)-
\end{eqnarray}
\begin{eqnarray}
2\operatorname{Li_{2}}\left(\frac{1}{\Psi^{10}}\right)-\frac{\pi^{2}}{3}+8\ln^{2}(\Psi)=0.
\nonumber
\end{eqnarray}

We can write the following two-term identities for $\Psi$.

\vskip 0.1in

$\operatorname{Li_{2}}(\Psi-1)=\operatorname{Li_{2}}(\frac{1}{\Psi^{2}})$
\vskip 0.01in
$\operatorname{Li_{2}}(\Psi^{2}-2)=\operatorname{Li_{2}}(\frac{1}{\Psi^{5}})$
\vskip 0.01in
$\operatorname{Li_{2}}(\frac{1}{\Psi})=\operatorname{Li_{2}}(-\frac{1}{\Psi^{2}})+\frac{\pi^{2}}{6}-\frac{5}{2}\ln^{2}(\Psi)$
\vskip 0.01in
$\operatorname{Li_{2}}(\frac{1}{\Psi})=-\operatorname{Li_{2}}(\frac{1}{\Psi^{3}})+\frac{\pi^{2}}{6}-3\ln^{2}(\Psi)$

\vskip 0.1in

From the above identities, we can also write two following formulae, as shown below.
\vskip 0.1in
$\operatorname{Li_{2}}(\frac{1}{\Psi^{3}})+\operatorname{Li_{2}}(-\frac{1}{\Psi^{2}})=-\frac{1}{2}\ln^{2}(\Psi)$
\vskip 0.01in
$\operatorname{Li_{2}}(3-2\Psi)+\operatorname{Li_{2}}(-\frac{1}{2\Psi^{5}})=\ln(\frac{1}{2}\Psi^{2})\ln(3\Psi^{7}-2\Psi^{8})-
\frac{1}{2}\ln^{2}(2)+2\ln(2)\ln(\Psi)$
\vskip 0.01in
$-2\ln^{2}(\Psi)$.

\subsection{The second smallest Pisot number $\operatorname{\theta_{1}}$ and the dilogarithm}

The constant $\theta_{1}$ is the second smallest Pisot number, as the positive root of the quartic trinomial $x^{4}-x^{3}-1=0$. 
We proceed to generate a six-term ladder with $\theta_{1}$ by applying the five-term gemini-identity. Let us set the integration
limits in such a way that $x_{1}=\ln(\theta)$ and $x_{2}=\ln(1+\theta^{2})$. Hence, the shape factor $a$ is given by

\vskip 0.1in

$\frac{\theta+a}{\theta-1}=\theta^{2}+1 \Rightarrow a=\theta^{3}-\theta^{2}-1=\theta^{2}(\theta-1)-1=\theta^{2}\cdot\frac{1}{\theta^{3}}-1=$
\vskip 0.01in
$\frac{1}{\theta}-1=\frac{1-\theta}{\theta}=-\frac{1}{\theta^{4}}$.

\vskip 0.1in

\begin{theorem}\label{Pisotladdertheorem}
The valid ladder relation such that 
\begin{multline*}
\operatorname{Li}_2\left(\frac{1}{\theta ^{14}}\right)-2 \operatorname{Li}_2\left(\frac{1}{\theta ^7}\right)-
2\operatorname{Li}_2\left(\frac{1}{\theta^5}\right)+\operatorname{Li}_2\left(\frac{1}{\theta ^4}\right) + \\
2 \operatorname{Li}_2\left(\frac{1}{\theta ^2}\right) + 
2 \operatorname{Li}_2\left(\frac{1}{\theta}\right)=\frac{\pi ^2}{3}-5 \log ^2\left(\frac{1}{\theta }\right) 
\end{multline*}
holds. 
\end{theorem}

\textit{Proof.} This evaluation is based on the $\gemini_{-\frac{1}{\theta^{4}}}(x)$-function. Let us insert the initial values in the five term
identity. For simplicity, let us denote $\theta_{1}=\theta$. Hence, we can write

\vskip 0.1in

$\operatorname{Li_{2}}(\frac{1}{\theta^{5}})-\operatorname{Li_{2}}(\frac{1}{\theta})-\operatorname{Li_{2}}(\frac{1}{\theta^{4}})+\frac{\pi^{2}}{6}
-\ln(\theta)\ln(\theta^{2}+1)+\operatorname{Li_{2}}(\frac{1}{\theta^{6}+\theta^{4}})-\operatorname{Li_{2}}(\frac{1}{1+\theta^{2}})$
\vskip 0.01in
=0$ \Rightarrow$
\vskip 0.01in
$\operatorname{Li_{2}}(\frac{1}{\theta^{5}})-\operatorname{Li_{2}}(\frac{1}{\theta})-\operatorname{Li_{2}}(\frac{1}{\theta^{4}})+\frac{\pi^{2}}{6}
-\ln(\theta)\ln(\theta^{2}+1)+\operatorname{Li_{2}}(\frac{1}{\theta^{7}+1})-\operatorname{Li_{2}}(-\theta^{2})-$
\vskip 0.01in
$\frac{\pi^{2}}{6}
+\frac{1}{2}\ln(\theta^{2}+1)\ln(\frac{\theta^{2}+1}{\theta^{4}})=0 \Rightarrow$
\vskip 0.01in
$\operatorname{Li_{2}}(\frac{1}{\theta^{5}})-\operatorname{Li_{2}}(\frac{1}{\theta})-\operatorname{Li_{2}}(\frac{1}{\theta^{4}})+\frac{\pi^{2}}{6}
-\ln(\theta)\ln(\theta^{2}+1)+\operatorname{Li_{2}}(\frac{1}{\theta^{7}+1})+\operatorname{Li_{2}}(-\frac{1}{\theta^{2}})+$
\vskip 0.01in
$2\ln^{2}(\theta)+\frac{1}{2}\ln(\theta^{2}+1)\ln(\frac{\theta^{2}+1}{\theta^{4}})=0 \Rightarrow$
\vskip 0.01in
$\operatorname{Li_{2}}(\frac{1}{\theta})+\operatorname{Li_{2}}(\frac{1}{\theta^{2}})+\frac{1}{2}\operatorname{Li_{2}}(\frac{1}{\theta^{4}})-
\operatorname{Li_{2}}(\frac{1}{\theta^{5}})-\operatorname{Li_{2}}(\frac{1}{\theta^{7}+1})-\frac{\pi^{2}}{6}+\ln(\theta)\ln(\theta^{2}+1)-$
\vskip 0.01in
$2\ln^{2}(\theta)-\frac{1}{2}\ln(\theta^{2}+1)\ln(\frac{\theta^{2}+1}{\theta^{4}})=0 \Rightarrow$
\vskip 0.01in
$\operatorname{Li_{2}}(\frac{1}{\theta})+\operatorname{Li_{2}}(\frac{1}{\theta^{2}})+\frac{1}{2}\operatorname{Li_{2}}(\frac{1}{\theta^{4}})-
\operatorname{Li_{2}}(\frac{1}{\theta^{5}})-\operatorname{Li_{2}}(-\theta^{7})-\frac{\pi^{2}}{3}+\ln(\theta)\ln(\theta^{2}+1)+$
\vskip 0.01in
$\frac{1}{2}\ln(1+\theta^{7})\ln(\frac{1+\theta^{7}}{\theta^{14}})-2\ln^{2}(\theta)-\frac{1}{2}\ln(\theta^{2}+1)\ln(\frac{\theta^{2}+1}{\theta^{4}})=0
\Rightarrow$
\vskip 0.01in
$\operatorname{Li_{2}}(\frac{1}{\theta})+\operatorname{Li_{2}}(\frac{1}{\theta^{2}})+\frac{1}{2}\operatorname{Li_{2}}(\frac{1}{\theta^{4}})
-\operatorname{Li_{2}}(\frac{1}{\theta^{5}})+\operatorname{Li_{2}}(-\frac{1}{\theta^{7}})-\frac{\pi^{2}}{6}+\ln(\theta)\ln(\theta^{2}+1)+$
\vskip 0.01in
$\frac{45}{2}\ln^{2}(\theta)+\frac{1}{2}\ln(1+\theta^{7})\ln(\frac{1+\theta^{7}}{\theta^{14}})-
\frac{1}{2}\ln(\theta^{2}+1)\ln(\frac{\theta^{2}+1}{\theta^{4}})=0$

\subsection{Two ladders related to a quartic equation $x^{4}-x-1=0$}

Let us denote the root of a quartic trinomial $x^{4}-x-1=0$ as $x=a_{4}$. The subscript stands for the highest exponents of this equation. 
For simplicity, let us again denote $a_{4}=a$. 

\begin{theorem}
The valid ladder relation such that 
\begin{multline*}
\operatorname{Li}_2\left(\frac{1}{a^{10}}\right)-2 \operatorname{Li}_2\left(\frac{1}{a^5}\right)-\operatorname{Li}_2\left(\frac{1}{a^4}\right)-2
\operatorname{Li}_2\left(\frac{1}{a^3}\right) + \\
4 \operatorname{Li}_2\left(\frac{1}{a^2}\right)+2 \operatorname{Li}_2\left(\frac{1}{a}\right) = \frac{\pi ^2}{3}-5 \log ^2(a)
\end{multline*}
holds. 
\end{theorem}

\textit{Proof.} Next, we derive a ladder by applying the $\gemini_{-\frac{1}{a^{2}}}(x)$-function. Let the integration limits be such that

\vskip 0.1in

$x_{1}=\ln(a)$ and $x_{2}=\ln(\frac{a-\frac{1}{a^{2}}}{a-1})=\ln(\frac{a^{3}-1}{a^{2}(a-1)})=\ln(\frac{a^{2}+a+1}{a^{2}})=\ln(\frac{a^{2}+a^{4}}{a^2})=$
\vskip 0.01in
$\ln(1+a^{2})$.

\vskip 0.05in

Hence, we get the equation shown below.

\vskip 0.1in

$\operatorname{Li_{2}}(\frac{1}{a^{3}})-\operatorname{Li_{2}}(\frac{1}{a})-\operatorname{Li_{2}}(\frac{1}{a^{2}})+\frac{\pi^{2}}{6}-\ln(a)\ln(a^{2}+1)+
\operatorname{Li_{2}}(\frac{1}{a^{4}+a^{2}})-\operatorname{Li_{2}}(\frac{1}{a^{2}+1})=$
\vskip 0.01in
$0 \Rightarrow$
\vskip 0.01in
$\operatorname{Li_{2}}(\frac{1}{a^{3}})-\operatorname{Li_{2}}(\frac{1}{a})-\operatorname{Li_{2}}(\frac{1}{a^{2}})-\ln(a)\ln(a^{2}+1)
+\operatorname{Li_{2}}(\frac{1}{a^{5}+1})-\operatorname{Li_{2}}(-a^{2})+$
\vskip 0.01in
$\frac{1}{2}\ln(a^{2}+1)\ln(\frac{a^{2}+1}{a^{4}})=0 \Rightarrow $
\vskip 0.01in
$\operatorname{Li_{2}}(\frac{1}{a^{3}})-\operatorname{Li_{2}}(\frac{1}{a})-\operatorname{Li_{2}}(\frac{1}{a^{2}})-\ln(a)\ln(a^{2}+1)
+\operatorname{Li_{2}}(\frac{1}{a^{5}+1})+\operatorname{Li_{2}}(-\frac{1}{a^{2}})+\frac{\pi^{2}}{6}+$
\vskip 0.1in
$2\ln^{2}(a)+\frac{1}{2}\ln(a^{2}+1)\ln(\frac{a^{2}+1}{a^{4}})=0\Rightarrow$
\vskip 0.1in
$\operatorname{Li_{2}}(\frac{1}{a^{3}})-\operatorname{Li_{2}}(\frac{1}{a})-2\operatorname{Li_{2}}(\frac{1}{a^{2}})
-\ln(a)\ln(a^{2}+1)+\operatorname{Li_{2}}(\frac{1}{a^{5}+1})+\frac{1}{2}\operatorname{Li_{2}}(\frac{1}{a^{4}})+\frac{\pi^{2}}{6}+$
\vskip 0.1in
$2\ln^{2}(a)+\frac{1}{2}\ln(a^{2}+1)\ln(\frac{a^{2}+1}{a^{4}})=0 \Rightarrow $
\vskip 0.1in
$\operatorname{Li_{2}}(\frac{1}{a})+2\operatorname{Li_{2}}(\frac{1}{a^{2}})-\operatorname{Li_{2}}(\frac{1}{a^{3}})-
\frac{1}{2}\operatorname{Li_{2}}(\frac{1}{a^{4}})
-\operatorname{Li_{2}}(-a^{5})-\frac{\pi^{2}}{3}+\frac{1}{2}\ln(a^{5}+1)\ln(\frac{a^{5}+1}{a^{10}})+$
\vskip 0.1in
$\ln(a)\ln(a^{2}+1)-2\ln^{2}(a)-\frac{1}{2}\ln(a^{2}+1)\ln(\frac{a^{2}+1}{a^{4}})=0 $


\begin{theorem}
The valid ladder relation such that 
\begin{multline*}
2\operatorname{Li_{2}}\left(\frac{1}{a}\right)+\operatorname{Li_{2}}\left(\frac{1}{a^{2}}\right)+2\operatorname{Li_{2}}\left(\frac{1}{a^{3}}\right)
-\operatorname{Li_{2}}\left(\frac{1}{a^{4}}\right)+2\operatorname{Li_{2}}\left(\frac{1}{a^{5}}\right) + \\ 
2\operatorname{Li_{2}}\left(\frac{1}{a^{7}}\right)-\operatorname{Li_{2}}\left(\frac{1}{a^{10}}\right)-\frac{2\pi^{2}}{3}+38\ln^{2}(a)=0
\end{multline*}
holds. 
\end{theorem}

\textit{Proof.} We apply the $\gemini_{a^{2}}(x)$-function with $x_{1}=\ln(a^{3})$ and $x_{2}=\ln(\frac{a^{3}+a^{2}}{a^{3}-1})=\ln(\frac{a^{2}(a+1)}{\frac{1}{a}})=
\ln(a \cdot a^{2} \cdot a^{4})=
\ln(a^{7})$. By inserting the initial values in the five-term identity, we get the formulae, as shown below.

\vskip 0.1in

$\operatorname{Li_{2}}(-\frac{1}{a})-\operatorname{Li_{2}}(\frac{1}{a^{3}})-\operatorname{Li_{2}}(-a^{2})+\frac{\pi^{2}}{6}-21\ln^{2}(a)+
\operatorname{Li_{2}}(-\frac{1}{a^{5}})-\operatorname{Li_{2}}(\frac{1}{a^{7}})=0 \Rightarrow $

\vskip 0.1in

$\frac{1}{2}\operatorname{Li_{2}}(\frac{1}{a^{2}})-\operatorname{Li_{2}}(\frac{1}{a})-
\operatorname{Li_{2}}(\frac{1}{a^{3}})+\operatorname{Li_{2}}(-\frac{1}{a^{2}})
+\frac{\pi^{2}}{3}-19\ln^{2}(a)+\operatorname{Li_{2}}(-\frac{1}{a^{5}})-\operatorname{Li_{2}}(\frac{1}{a^{7}})=$
\vskip 0.01in
$0$


\subsection{A general two-term identity related to a trinomial equation of $x^{n}-x^{m}-1=0$}

The preceding constants involved in our ladder constructions all are roots of the trinomial equation of the form $x^{n}-x^{m}-1=0$, where $n$ 
and $m$ are positive real numbers in such a way that $n>m$. Next, we generate a general identity formula based on this trinomial equation.
First, we apply the reflection formula as follows,

\vskip 0.1in

$\operatorname{Li_{2}}(\frac{1}{x^{n}})=-\operatorname{Li_{2}}(\frac{x^{n}-1}{x^{n}})+\frac{\pi^{2}}{6}-\ln(\frac{1}{x^{n}})\ln(\frac{x^{n}-1}{x^{n}})=
-\operatorname{Li_{2}}(\frac{x^{m}}{x^{n}})+\frac{\pi^{2}}{6}+\ln(x^{n})\ln(\frac{x^{m}}{x^{n}})=-\operatorname{Li_{2}}(\frac{1}{x^{n-m}})
+\frac{\pi^{2}}{6}-\ln(x^{n})\ln(x^{n-m})=-\operatorname{Li_{2}}(\frac{1}{x^{n-m}})+\frac{\pi^{2}}{6}-n(n-m)\ln^{2}(x)=$
\vskip 0.01in
$-\operatorname{Li_{2}}(\frac{1}{x^{n-m}})+\frac{\pi^{2}}{6}-n^{2}\ln^{2}(x)+nm\ln^{2}(x)$.
On the other hand, we can write a second formula by applying Landen's identity for representing the same thing with a different way, so that 
\begin{multline*}
\operatorname{Li_{2}}\left(\frac{1}{x^{n}}\right)=\operatorname{Li_{2}}\left(\frac{1}{x^{m}+1}\right)=\operatorname{Li_{2}}(-x^{m})+\frac{\pi^{2}}{6}
-\frac{1}{2}\ln(x^{m}+1)\ln\left(\frac{x^{m}+1}{x^{2m}}\right) = \\ 
\operatorname{Li_{2}}(-x^{m})+\frac{\pi^{2}}{6}-\frac{1}{2}\ln(x^{n})\ln(x^{n-2m})=
\operatorname{Li_{2}}(-x^{m})+\frac{\pi^{2}}{6}-\frac{n^{2}}{2}\ln^{2}(x)+nm\ln^{2}(x). 
\end{multline*}
Next, we combine the above equations and we can write
$$-\operatorname{Li_{2}}\left(\frac{1}{x^{n-m}}\right)+\frac{\pi^{2}}{6}-n^{2}\ln^{2}(x) + 
nm\ln^{2}(x)=\operatorname{Li_{2}}(-x^{m})+\frac{\pi^{2}}{6}-\frac{n^{2}}{2}\ln^{2}(x)+nm\ln^{2}(x).$$
Consequently, we have that 
\begin{eqnarray}\label{formula46}
\operatorname{Li_{2}}\left(\frac{1}{x^{n-m}}\right)+\operatorname{Li_{2}}\left(-x^{m}\right)=-\frac{1}{2}n^{2}\ln^{2}(x).
\end{eqnarray}

Below is a list of all the trinomial equation identities related to the constants we have dealt with in the previous sections.

\vskip 0.1in

$\phi^{2}-\phi-1=0 \Rightarrow \operatorname{Li_{2}}(\frac{1}{\phi})+\operatorname{Li_{2}}(-\phi)=-2\ln^{2}(\phi)$
\vskip 0.01in
$P^{3}-P-1=0 \Rightarrow \operatorname{Li_{2}}(\frac{1}{P^{2}})+\operatorname{Li_{2}}(-P)=-\frac{9}{2}\ln^{2}(P)$
\vskip 0.01in
$P^{5}-P^{4}-1=0 \Rightarrow \operatorname{Li_{2}}(\frac{1}{P})+\operatorname{Li_{2}}(-P^{4})=-\frac{25}{2}\ln^{2}(P)$
\vskip 0.01in
$\Psi^{3}-\Psi^{2}-1=0 \Rightarrow \operatorname{Li_{2}}(\frac{1}{\Psi})+\operatorname{Li_{2}}(-\Psi^{2})=-\frac{9}{2}\ln^{2}(\Psi))$
\vskip 0.01in
$\theta_{1}^{4}-\theta_{1}^{3}-1=0 \Rightarrow \operatorname{Li_{2}}(\frac{1}{\theta_{1}})+\operatorname{Li_{2}}(-\theta_{1}^{3})=
-8\ln^{2}(\theta_{1}^{4})$
\vskip 0.01in
$a_{4}^{4}-a_{4}-1=0 \Rightarrow \operatorname{Li_{2}}(\frac{1}{a_{4}^{3}})+\operatorname{Li_{2}}(-a_{4})=-8\ln^{2}(a_{4})$

\subsection{The transcendental intersection point of $\gemini_{0}(x)$ and $\gemini_{0}^{rot}(x)$}
 If an arbitrary gemini function $\gemini_{a}(x)$ is rotated by $45^{\circ}$ counterclockwise, then the common intersection point of 
the rotated version and the original function is on the line of $y=(\sqrt{2}+1)x=x\tan(\frac{3\pi}{8})$. The equation for this
intersection point can be written as $$y=(\sqrt{2}+1)x=\ln\left(\frac{1+ae^{-x}}{1-e^{-x}}\right), $$ giving us that
$e^{x(\sqrt{2}+2)}- e^{x(\sqrt{2}+1)}- e^{x}-a=0$. Next, we substitute such that $k_{a}=e^{x}$ and we can write
$$k_{a}^{\sqrt{2}+2}-k_{a}^{\sqrt{2}+1}-k_{a}-a=0.$$
Thus, the coordinates corresponding to the intersection point of the $\gemini_{a}(x)$-function are such that $x=\ln(k_{a})$ and
$y=\ln(\frac{k_{a}+a}{k_{a}-1})$. Let us consider the $\gemini_{0}(x)$-function, whose $x$-coordinate of the intersection point
is $\ln(k_{0})$. Here, $k_{0}$ is the root of the equation such that 
 $$k_{0}^{\sqrt{2}+1}-k_{0}^{\sqrt{2}}-1=0, $$
 and there is no known closed form for $k_{0}$. 
 The constant $k_{0}$ satisfies the trinomial equation
 identity \eqref{formula46}, and this gives us that 
\begin{eqnarray}\label{formula47}
\operatorname{Li_{2}}\left(\frac{1}{k_{0}}\right)+\operatorname{Li_{2}}\left(-k_{0}^{\sqrt{2}}\right)=-\frac{1}{2}(\sqrt{2}+1)^{2}\ln^{2}(k_{0}).
\end{eqnarray}

 We can derive a related identity in a manner suggested in Figure \ref{Figure6}. First, we introduce a special property of gemini functions. 
 Given two line segments from the origin to integration limit points of a curve of a gemini function, these line segments
 form a sector-like figure, 
 within the curve between the integration limits. The area $A_{s}$ of
this figure is equal to the area $A_{c}$, which lies between the integration limits $x_{1}$ and $x_{2}$. Proving this is a simple task.
We can write $A_{tot}=A_{a}+A_{r}+A_{c}+A_{a}$. Secondly, we can also write $A_{tot}=A_{a}+\frac{1}{2}A_{r}+A_{s}+\frac{1}{2}A_{r}+A_{a} \Rightarrow A_{s}=A_{c}$,
where $A_{a}$ is the apex area and $A_{r}=x_{1}x_{2}$.
 This derivation is based on the $\gemini_{0}(x)$- and $\gemini_{0}^{rot}(x)$-functions. 
 By applying the $\gemini_{0}(x)$-function, we can write 
\begin{multline*}
 A_{s}=A_{c}=\int_{\ln(k_{0})}^{\ln(\frac{k_{0}}{k_{0}-1})}\gemini_{0}(x) \, 
 dx=\int_{\ln(k_{0})}^{\ln(\frac{k_{0}}{k_{0}-1})}\ln\left( 
 \frac{1}{1-e^{-x}} \right) \, dx = \\ 
 \operatorname{Li_{2}}\left(
 \frac{1}{k_{0}} \right) - 
 \operatorname{Li_{2}}\left(
 \frac{k_{0}-1}{k_{0}} \right) = 
2\operatorname{Li_{2}}\left(
 \frac{1}{k_{0}} \right) - 
 \frac{\pi^{2}}{6}+\ln(k_{0}) \ln\left(
 \frac{k_{0}}{k_{0}-1} \right). 
\end{multline*}

\begin{figure}[htbp!]
\begin{center}
\resizebox{7.0cm}{!}
{\includegraphics{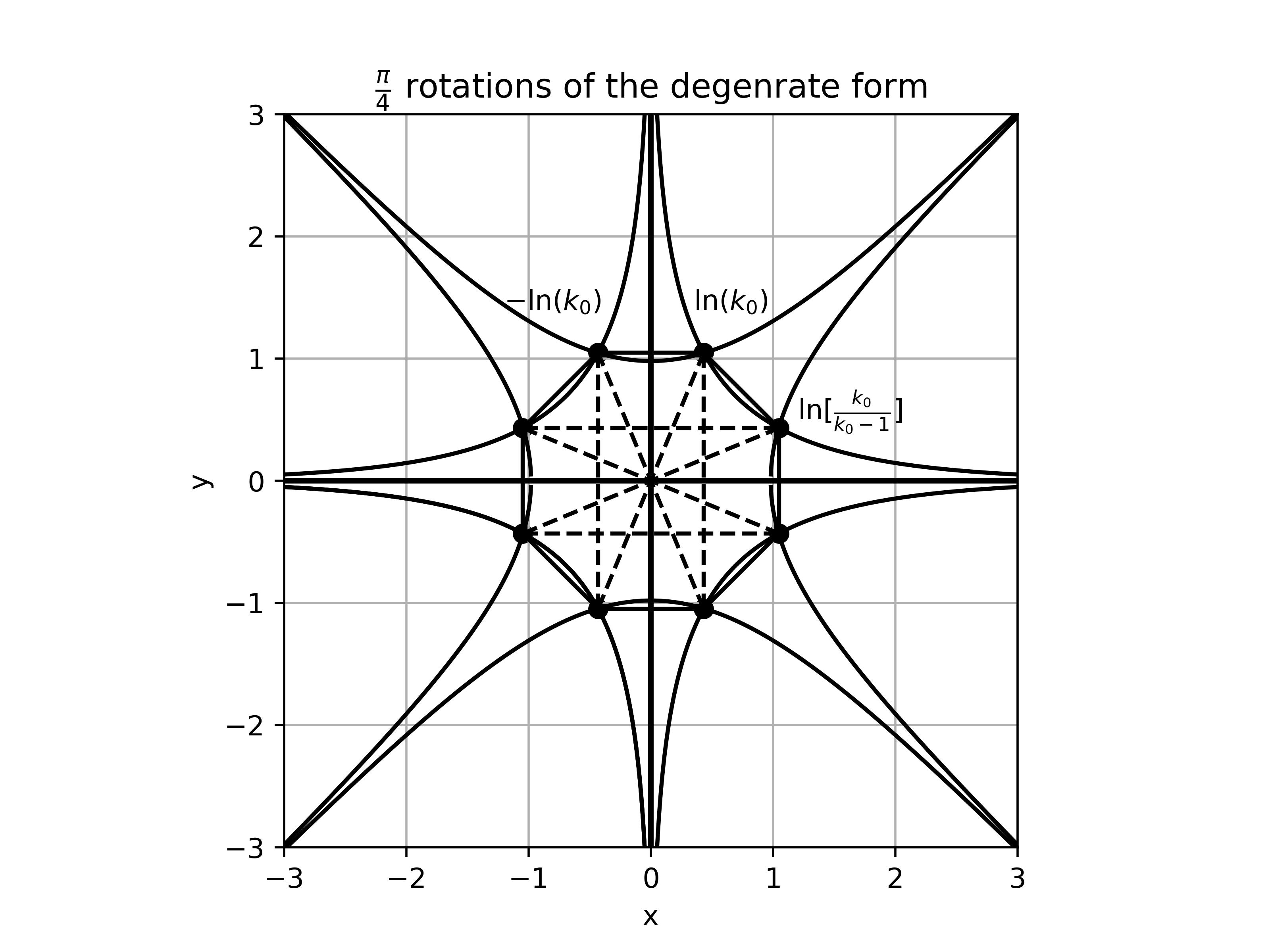}}
\end{center}
\caption{\label{Figure6} This star-like figure is obtained by rotating counter clockwise the $\gemini_{0}(x)$-function seven times
by an angle of $\frac{\pi}{4}$. The area of the figure is finite. If this figure is rotated about the $x$- or $y$-axis, the volume
of the resulting solid of revolution is also finite. }
\end{figure}

 We can obtain the same area from the rotated function $\gemini_{0}^{rot}(x)$. By subtracting two triangular areas, relative to the
 definite integral being used, we obtain
\begin{multline*}
 A_{s}=\int_{-\ln(k_{0})}^{+\ln(k_{0})}\gemini_{0}^{rot}(x) \, dx-\frac{1}{2}A_{r}-\frac{1}{2}A_{r} = \\ 
 \int_{-\ln(k_{0})}^{+\ln(k_{0})}\frac{1}{\sqrt{2}}\ln(2\cosh(x\sqrt{2})+2) \, dx-x_{1}x_{2}, 
\end{multline*}
 and this gives us that the same area $A_{s}$ may be expressed as 
\begin{multline*}
\operatorname{Li_{2}}(-\frac{1}{k_{0}^{\sqrt{2}}})-\operatorname{Li_{2}}(-k_{0}^{\sqrt{2}})-\ln(k_{0})\ln(\frac{k_{0}}{k_{0}-1}) = \\ 
-2\operatorname{Li_{2}}(-k_{0}^{\sqrt{2}})-\frac{\pi^{2}}{6}-\frac{1}{2}\ln^{2}(k_{0}^{\sqrt{2}})-\ln(k_{0})\ln(\frac{k_{0}}{k_{0}-1}).
\end{multline*}
Through a combination of the above formulas, we find that 
\begin{multline*}
2\operatorname{Li_{2}}(\frac{1}{k_{0}})-\frac{\pi^{2}}{6}+\ln(k_{0})\ln(\frac{k_{0}}{k_{0}-1}) = \\ 
-2\operatorname{Li_{2}}(-k_{0}^{\sqrt{2}})-\frac{\pi^{2}}{6}-\frac{1}{2}\ln^{2}(k_{0}^{\sqrt{2}})-\ln(k_{0})\ln(\frac{k_{0}}{k_{0}-1}).
\end{multline*}
Consequently, we have that 
\begin{eqnarray}\label{formula48} 
\operatorname{Li_{2}}\left(\frac{1}{k_{0}}\right)+\operatorname{Li_{2}}\left(-k_{0}^{\sqrt{2}}\right)=
-\ln(k_{0})\ln\left(\frac{k_{0}\sqrt{k_{0}}}{k_{0}-1}\right).
\end{eqnarray}

The identities shown in \eqref{formula47} and \eqref{formula48} are exactly the same despite the different representations of the
constant terms. By setting an equal sign between the constant terms, we get
 $$-\frac{1}{2}(\sqrt{2}+1)^{2}\ln^{2}(k_{0})=-\ln(k_{0})\ln(\frac{k_{0}\sqrt{k_{0}}}{k_{0}-1}), $$ 
 so that 
 $k_{0}^{\sqrt{2}+1}-k_{0}^{\sqrt{2}}-1=0$. 
 Next, we calculate the area of the star like figure shown in Figure \ref{Figure6} with the aid of the constant $k_{0}$. First, the area of the
octagon must be evaluated. Its area is given by
 $$A_{oct}=8\ln(k_{0})\ln(\frac{k_{0}}{k_{0}-1}).$$
 Observe that $k_{0}$ satisfies the relation 
 $$\frac{\ln(k_{0})}{\ln(\frac{k_{0}}{k_{0}-1})}=\sqrt{2}-1 \Rightarrow
\ln(\frac{k_{0}}{k_{0}-1})=\ln(k_{0}^{\sqrt{2}+1}).$$ The area of a single vertex of this star like figure is given by
 $$A_{v}=2\int_{\ln(\frac{k_{0}}{k_{0}-1})}^{\infty}\gemini_{0}(x) \, dx=2\operatorname{Li_{2}}(\frac{k_{0}-1}{k_{0}})=
2\operatorname{Li_{2}}(\frac{1}{k_{0}^{\sqrt{2}+1}}).$$ 
 Hence, the total area of this star like figure is given by
 $$A_{tot}=A_{oct}+8A_{v}=8\ln(k_{0})\ln\left(\frac{k_{0}}{k_{0}-1}\right)+16\operatorname{Li_{2}}\left(\frac{k_{0}-1}{k_{0}}\right).$$ 
We can derive a third two-term identity for the constant $k_{0}$ based on its arithmetical properties, so that 
\begin{multline*}
 \operatorname{Li_{2}}\left( 
 \frac{k_{0}-1}{k_{0}} \right) = 
 \operatorname{Li_{2}}\left( 
 \frac{1}{k_{0}^{\sqrt{2}+1}} \right) \Rightarrow \\ 
-\operatorname{Li_{2}}\left(\frac{1}{k_{0}} \right) + 
 \frac{\pi^{2}}{6}-\ln(k_{0})\ln(\frac{k_{0}}{k_{0}-1})=
\operatorname{Li_{2}}\left(\frac{1}{k_{0}^{\sqrt{2}+1}} \right). 
\end{multline*} 
 Consequently, we have that 
 $$-\operatorname{Li_{2}}\left(\frac{1}{k_{0}} \right) + 
 \frac{\pi^{2}}{6}-\ln(k_{0})\ln(k_{0}^{\sqrt{2}+1}) = 
\operatorname{Li_{2}}\left(\frac{1}{k_{0}^{\sqrt{2}+1}} \right), $$
so that 
 $$-\operatorname{Li_{2}}\left(\frac{1}{k_{0}} \right) + 
 \frac{\pi^{2}}{6}-(\sqrt{2}+1)\ln^{2}(k_{0}) = 
\operatorname{Li_{2}}\left(\frac{1}{k_{0}^{\sqrt{2}+1}} \right), $$
and continuing in a similar manner gives us that 

\begin{eqnarray}\label{formula50}
\operatorname{Li_{2}}\left(\frac{1}{k_{0}^{\sqrt{2}+1}}\right)+\operatorname{Li_{2}}\left(-\frac{1}{k_{0}^{\sqrt{2}}}\right)+\frac{1}{2}\ln^{2}(k_{0})=0.
\end{eqnarray}

\begin{figure}[!h]
\begin{center}
\resizebox{10.0cm}{!} 
{\includegraphics{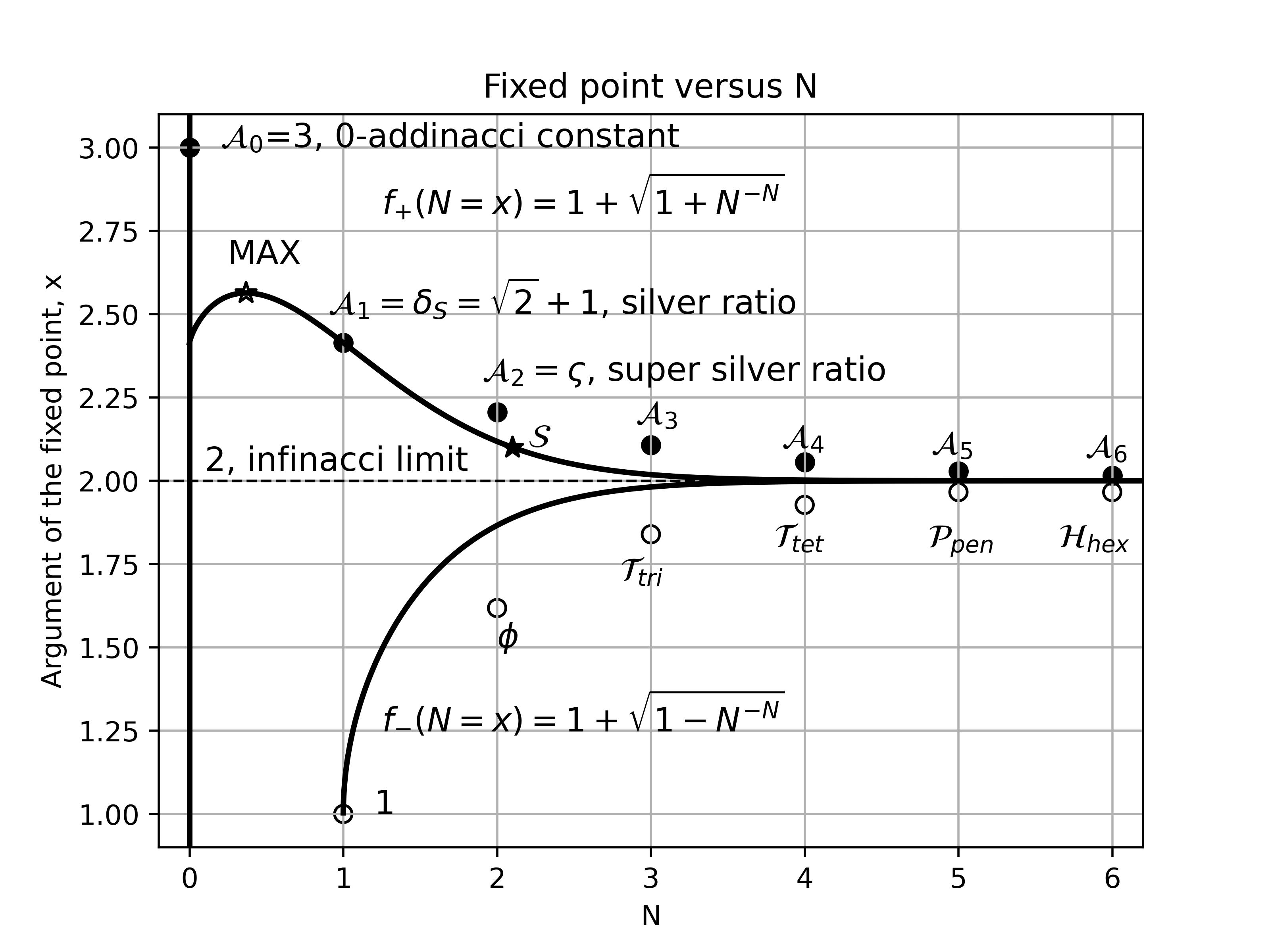}}
\end{center}
\caption{\label{Figure7} The first 7 N-addinacci and the 6 first N-bonacci constants together with two approximating functions are
illustrated in this plot. The star symbols correspond to the maximum value and the super fixed point of the approximated addinacci function.}
\end{figure}

\subsection{N-bonacci and N-addinacci constants related to the fixed point identity}

 The fixed point $x_{0}$ and the shape factor $a$ of the $\gemini_{a}(x)$-function has a simple relation with N-addinacci and
N-bonacci constants. These constants are the limiting ratios of the two successive terms in the sequences of Fibobacci n-step
numbers. These sequences are generalizations of the Fibonacci sequence where the limiting ratio is $\phi$. The $N$-bonacci
sequence starts with $N-1$ zeroes followed by 1, with subsequent terms being generated by taking the sum of the $N$ previous terms.
There exists another way to define a Fibonacci (or 2-bonacci) sequence. The sequence can also be written in a following way,
$F_{i}=2F_{i-1}-F_{i-3}$, where i is the index or the position of the number $F_{i}$ in the sequence. The $N$-addinacci sequence 
can be defined by adding these two terms instead of subtracting them. Hence, the recursive formula, e.g. for the 3-addinacci
sequence is given by $F_{i}=2F_{i-1}+F_{i-4}$. In general, the N-bonacci and N-addinacci sequences are defined by
$F_{i}=2F_{i-1}\pm F_{i-(N+1)}$. We can write the following equations for N-bonacciand and N-addinacci constants, as shown below.
The plus sign corresponds to N-addinacci constants and minus sign respectively to N-bonacci constants. 

\vskip 0.1in

$x=1+\sqrt{1\pm\frac{1}{x^{N-1}}} \Rightarrow x-2=\pm\frac{1}{x^{N}} \Rightarrow x^{N+1}-2x^{N}\mp1=0 \Rightarrow$
\vskip 0.01in
$x_{0}=\ln(x)=\ln(1+\sqrt{1\pm\frac{1}{x^{N-1}}})$ and $a=\pm\frac{1}{x^{N-1}}$.

\vskip 0.1in

We can approximate the equations $x=1+\sqrt{1\pm\frac{1}{x^{N-1}}}$ with the functions $f_{\pm}(N)=1+\sqrt{1\pm\frac{1}{N^{N}}}$,
as shown in Figure \ref{Figure7}. Both of the functions approach to the value 2, which is the infinacci constant. The infinacci constant is
the argument of a fixed point of the degenerate form of a gemini function, i.e., $x_{0}=\ln(2)$. The approximate addinacci-function
$f_{+}(N)=1+\sqrt{1+\frac{1}{N^{N}}}$ has two special points marked with star symbols in Figure \ref{Figure7}. The maximum is at $N=x=\frac{1}{e}$,
whose corresponding fixed point is such that $N_{0}=x_{0}=\ln(1+\sqrt{1+\sqrt[e]{e}})$. The other star denoted point with the symbol
$\mathcal{S}$ corresponds to the addinacci super fixed point. At that point $N_{s}=x_{s}=1+\sqrt{1+\frac{1}{x_{s}^{x_{s}}}}\approx2.100211$,
where $x_{0}=\ln(x_{s})$ and $a=+\frac{1}{x_{s}^{x_{s}}}$.

\begin{theorem}\label{theoremaddin}
The valid ladder relation 
\begin{multline*}
4\operatorname{Li_{2}}\left(\frac{1}{x}\right)-2\operatorname{Li_{2}}\left(\frac{1}{x^{N-1}}\right)+4\operatorname{Li_{2}}\left(\frac{1}{x^{N}}\right)
+ \\ 
\operatorname{Li_{2}}\left(\frac{1}{x^{2N-2}}\right)-2\operatorname{Li_{2}}\left(\frac{1}{x^{2N}}\right)-\frac{\pi^{2}}{3}+2\ln^{2}(x)=0
\end{multline*}
holds. 
\end{theorem}

\textit{Proof.} We employ the fixed-point identity in Theorem \ref{firstfixedpoint}, giving us that the N-bonacci constant formula contains only three
dilogarithm terms, with 
\begin{eqnarray}\label{bonacci}
2\operatorname{Li_{2}}\left(\frac{1}{x}\right)+\operatorname{Li_{2}}\left(\frac{1}{x^{N-1}}\right)-2\operatorname{Li_{2}}\left(\frac{1}{x^{N}}\right)
-\frac{\pi^{2}}{6}+\ln^{2}(x)=0. 
\end{eqnarray}

The ladder formula for N-addinacci constants includes five dilogarithm terms. In this case, the shape factor and the fixed point are
defined in such a way that $a=+\frac{1}{x^{N-1}}$ and $x_{0}=\ln(x)$. Hence, we get
$\operatorname{Li_{2}}\left(-\frac{1}{x^{N-1}} \cdot \frac{1}{x}\right)-\operatorname{Li_{2}}\left(\frac{1}{x}\right)-
\frac{1}{2}\operatorname{Li_{2}}\left(-\frac{1}{x^{N-1}}\right)+\frac{\pi^{2}}{12}-\frac{1}{2}\ln^{2}(x)=0$, 
and this gives us an equivalent version of the desired result. 

\begin{example}
Let us construct one addinacci ladder as an example. By setting $N=4$ the respective minimal polynomial is such that 
$x^{5}-2x^{4}-1=0$, which has one real root, corresponding to the 4-addinacci constant. The corresponding ladder is given by
\begin{multline*}
4\operatorname{Li_{2}}\left(\frac{1}{\mathcal{A}_{4}}\right)-2\operatorname{Li_{2}}\left(\frac{1}{\mathcal{A}_{4}^{3}}\right)+
4\operatorname{Li_{2}}\left(\frac{1}{\mathcal{A}_{4}^{4}}\right)+\operatorname{Li_{2}}\left(\frac{1}{\mathcal{A}_{4}^{6}}\right) - \\ 
2\operatorname{Li_{2}}\left(\frac{1}{\mathcal{A}_{4}^{8}}\right)+2\ln^{2}(\mathcal{A}_{4})-\frac{\pi^{2}}{3}=0. 
\end{multline*}
\end{example}

\begin{example}
It is a simple task to build a three-term 3-bonacci constant ladder by setting $N=3$. This ladder is given by
\begin{eqnarray}\label{tri1}
2\operatorname{Li_{2}}\left(\frac{1}{\mathcal{T}_{tri}}\right)+\operatorname{Li_{2}}\left(\frac{1}{\mathcal{T}_{tri}^{2}}\right)-
2\operatorname{Li_{2}}\left(\frac{1}{\mathcal{T}_{tri}^{3}}\right)
-\frac{\pi^{2}}{6}+\ln^{2}(\mathcal{T}_{tri})=0.
\end{eqnarray}
\end{example}

\begin{example}
We can build another tribonacci ladder by applying the five-term gemini-identity based on the $\gemini_{\mathcal{T}_{tri}^2}(x)$-function.
Let us set the lower integration limit so that $x_{1}=\ln(\mathcal{T}_{tri})$. Hence, the upper integration limit becomes such that
$x_{2}=\ln(\frac{\mathcal{T}_{tri}^{2}+\mathcal{T}_{tri}}{\mathcal{T}_{tri}-1})=\ln(\frac{\mathcal{T}_{tri}^{3}-1}{\mathcal{T}_{tri}-1})=
\ln(\mathcal{T}_{tri}^{2}+\mathcal{T}_{tri}+1)=\ln(\mathcal{T}_{tri}^{3})$. The outcome formula is
shown below.
\begin{multline*}
2\operatorname{Li_{2}}\left(\frac{1}{\mathcal{T}_{tri}}\right)+2\operatorname{Li_{2}}\left(\frac{1}{\mathcal{T}_{tri}^{2}}\right) + \\ 
2\operatorname{Li_{2}}\left(\frac{1}{\mathcal{T}_{tri}^{3}}\right)-\operatorname{Li_{2}}\left(\frac{1}{\mathcal{T}_{tri}^{4}}\right)-
\frac{\pi^{2}}{3}+3\ln^{2}\left(\mathcal{T}_{tri}\right)=0.
\end{multline*}
\end{example}

\begin{example}
Using \eqref{tri1}, we can obtain another four-term identity for the tribonacci constant with the same argument values, but with the
different multiplicative coefficients. This can be written by
\begin{multline*}
6\operatorname{Li_{2}}\left(\frac{1}{\mathcal{T}_{tri}}\right)+5\operatorname{Li_{2}}\left(\frac{1}{\mathcal{T}_{tri}^{2}}\right)+ \\ 
2\operatorname{Li_{2}}\left(\frac{1}{\mathcal{T}_{tri}^{3}}\right)-2\operatorname{Li_{2}}\left(\frac{1}{\mathcal{T}_{tri}^{4}}\right)-
\frac{5\pi^{2}}{6}+7\ln^{2}(\mathcal{T}_{tri})=0. 
\end{multline*}
\end{example}

\section{Calculation in the complex domain with gemini identities}\label{sectioncomplex}
 Gemini functions were initially defined for real arguments, but can be extended to complex arguments. The graphs shown in 
 Figure \ref{Figure5} give a hint that the identities allow the use of complex numbers. The maximum value of the arguments of the 
 first term
 in Theorem \ref{mainthreeterm} is only $\frac{1}{4}$, so as to achieve greater values extending up to the range between
$\frac{1}{4}$ and $1$, then the variable $a$, i.e., the shape factor in the argument must be a complex number. We have to keep in mind
that generally the valid domain for the shape factor in the five-term gemini-identities is such that $a\ge-1$. Later on, we will see
that we can break this restriction, at least in some cases.

\subsection{Derivation of a generalized identity in the complex domain}
 Next, we derive an identity, which enables us the determine at least three exact values for dilogarithms, where the argument is greater than one.
It is well known fact, when the argument of the dilogarithm is greater than one then the value of a dilogarithm includes a negative complex term
$-i \pi \ln(x)$ such that $\operatorname{Li_{2}}(x)=\mathfrak{Re}\biggl\{\operatorname{Li_{2}}(x)\bigg\}-i\pi\ln(x)$ for $x>1$ and
$x \in \mathbb{R}$. This identity is given in Lewin's text on dilogarithms
 \cite{Lew_58}. Our derivation is based on an application of 
 Theorem \ref{firstfixedpoint}. The idea is as follows. We make these two terms equal such that
$\operatorname{Li_{2}}(-\frac{a}{1+\sqrt{1+a}})=\operatorname{Li_{2}}(2-x)$. We know that the argument of the fixed point is related to the shape
factor in such a way that $x=1+\sqrt{1+a}$. Next, we manipulate the first term of the first fixed point identity with \eqref{Landensformula}.
Hence, we can write
\begin{multline*}
 \operatorname{Li_{2}}\left(-\frac{a}{1+\sqrt{1+a}}\right)=\operatorname{Li_{2}}\left(1-\sqrt{1+a}\right)=\\
\operatorname{Li_{2}}\left(\frac{1}{\sqrt{1+a}}\right)-\frac{\pi^{2}}{6}+ 
 \frac{1}{2}\ln\sqrt{a+1}\ln\left[\frac{\sqrt{a+1}}{(\sqrt{a+1}-1)^{2}}\right].
\end{multline*}
 Next, we combine these two first terms
 $$\operatorname{Li_{2}}\left(\frac{1}{\sqrt{1+a}}\right)-\frac{\pi^{2}}{6}+
\frac{1}{2}\ln\sqrt{a+1}\ln\left[\frac{\sqrt{a+1}}{(\sqrt{a+1}-1)^{2}}\right]=\operatorname{Li_{2}}\left(2-x\right).$$
 By applying the \eqref{reflectionformula} to the RHS, we get

\vskip 0.1in

$\operatorname{Li_{2}}\left(\frac{1}{\sqrt{1+a}}\right)-\frac{\pi^{2}}{6}+
\frac{1}{2}\ln\sqrt{a+1}\ln\left[\frac{\sqrt{a+1}}{(\sqrt{a+1}-1)^{2}}\right]=$
\vskip 0.01in
$-\operatorname{Li_{2}}\left(x-1\right)+\frac{\pi^{2}}{6}-
\ln(2-x)\ln(x-1).$

\vskip 0.1in

Next, we apply the above shown relation $x-1=\sqrt{1+a}.$ Hence, we can write

\vskip 0.1in

$\operatorname{Li_{2}}\left(\frac{1}{x-1}\right)+\operatorname{Li_{2}}\left(x-1\right)-\frac{\pi^{2}}{3}+
\frac{1}{2}\ln(x-1)\ln\left(\frac{x-1}{(2-x)^{2}}\right)+\ln(2-x)\ln(x-1)=0$.

\vskip 0.1in

Substituting $z=x-1$, we get

\vskip 0.1in

$\operatorname{Li_{2}}(z)+\operatorname{Li_{2}}\left(\frac{1}{z}\right)-\frac{\pi^{2}}{3}+\ln(z)\ln(-\sqrt{z})=0, \quad z<0 \vee z \geq 1$.

\vskip 0.1in

This result above can be rewritten in an other way. Hence, we get the more familiar complex domain identity, as shown below.

\begin{eqnarray}\label{complex}
\operatorname{Li_{2}}(z)+\operatorname{Li_{2}}\left(\frac{1}{z}\right)-\frac{\pi^{2}}{3}+\frac{1}{2}\ln^{2}(z)+i\pi\ln(z)=0,\quad z > 1
\end{eqnarray}

With the aid of this identity, the exact values for $\operatorname{Li_{2}}(2)=\frac{\pi^{2}}{4}-i\pi \ln(2)$, $\operatorname{Li_{2}}(\phi)=
\frac{7\pi^{2}}{30}+\frac{1}{2}\ln^{2}(\phi)-i\pi \ln(\phi)$ and $\operatorname{Li_{2}}(\phi^{2})=\frac{4\pi^{2}}{15}-\ln^{2}(\phi)-2i\pi \ln(\phi)$
can be evaluated based on the known values of $\operatorname{Li_{2}}(\frac{1}{2})$, $\operatorname{Li_{2}}(\frac{1}{\phi})$ and
$\operatorname{Li_{2}}(\frac{1}{\phi^{2}})$.

\subsection{Derivation of the exact value for $\operatorname{Li_{2}}(\frac{1-i}{2})$}

Let us begin this Section by calculating the exact value for $\operatorname{Li_{2}}(\frac{1-i}{2})$, which is already known, but we
want to show, how easily it can be derived with the aid of the fixed-point identity in 
 Theorem \ref{firstfixedpoint}. Next, we formulate
the familiar equation, which connects the argument of a fixed point and the shape factor in such a way that
$x=1+\sqrt{1-\frac{x^{2}}{x-1}} \Rightarrow x=1+i \Rightarrow x_{0}=\ln(x)$ and $a=-\frac{(1+i)^{2}}{1+i-1}=-2$. Hence, we can obtain
the following, letting $\operatorname{C}$ denote \emph{Catalan's constant} $\frac{1}{1^2}-\frac{1}{3^2}+\frac{1}{5^2}-\cdots$. 

\vskip 0.1in

$\operatorname{Li_{2}}(\frac{2}{1+i})-\operatorname{Li_{2}}(\frac{1}{1+i})-\frac{1}{2}\operatorname{Li_{2}}(2)+\frac{\pi^{2}}{12}-
\frac{1}{2}\ln^{2}(1+i)=0 \Rightarrow $

\vskip 0.1in

$-\operatorname{Li_{2}}(\frac{1-i}{2})+\operatorname{Li_{2}}(1-i)-\frac{\pi^2}{8}+\frac{i\pi \ln(2)}{2}+\frac{\pi^{2}}{12}+
\frac{\pi^{2}}{32}-\frac{\ln^{2}(2)}{8}-\frac{i\pi \ln(2)}{8}=0 \Rightarrow $

\vskip 0.1in

$-\operatorname{Li_{2}}(\frac{1-i}{2})+\operatorname{Li_{2}}(1-i)-\frac{\pi^{2}}{96}-\frac{\ln^{2}(2)}{8}+\frac{3i\pi \ln(2)}{8}=0 \Rightarrow $

\vskip 0.1in

$-\operatorname{Li_{2}}(\frac{1-i}{2})-\operatorname{Li_{2}}(i)+\frac{\pi^{2}}{6}-\ln(i)\ln(1-i)-\frac{\pi^{2}}{96}-
\frac{\ln^{2}(2)}{8}+\frac{3i\pi \ln(2)}{8}=0 \Rightarrow $
 Consequently, we have that 
$$\operatorname{Li_{2}}\left(\frac{1-i}{2}\right)=\frac{5\pi^{2}}{96}-\frac{\ln^{2}(2)}{8}+i\left[\frac{\pi \ln(2)}{8}-\operatorname{C}\right]. 
$$
 The above fixed point identity can also be shown to hold 
 in the complex domain. 
 Observe that the evaluation 
 $\operatorname{Li_{2}}(i)=-\frac{\pi^{2}}{48}+i\operatorname{C}$ is required, as above. 

\subsection{Rederiving two recent results}
 
The most recent results related to the complex valued two-term identities can be found in the work of Campbell \cite{Campbell202122}.
Campbell's method to prove two-term dilogarithm evaluations is based on series transform and Legendre polynomial expansions. We can
reproduce almost all their results and infinitely of similar kind of identities with the aid of the five-term gemini-identity obtained
 from the fundamental form. Our method is based on the well known fact. If the absolute value of the complex argument is equal to one
 then the dilogarithm can be evaluated analytically with the aid of the Clausen-function $\operatorname{Cl}_{2}(\theta)$ and trigamma
 function $\psi_{1}(x)$. See for example Lewin's monograph on dilogarithms
 \cite{Lew_58}. In this case, the respective angles of the complex argument value must be rational
 multiples of $\pi$ radians, which also means a rational fraction of a circle. Campbell has proved the following identity, 
 which was initially discovered by Ramanujan in an equivalent way, 
 whereas Campbell obtained a new Fourier--Legendre-based derivation. 

\begin{theorem} 
 (Ramanujan, 1915) The evaluation 
 $$\operatorname{Li_{2}}\left(i(2-\sqrt{3})\right)-\operatorname{Li_{2}}\left(-i(2-\sqrt{3})\right)=
\frac{2i\sqrt{7-4\sqrt{3}}\left[8C-\pi \ln(2+\sqrt{3})\right]}{3(8-4\sqrt{3})}$$ 
 holds true \cite{Ramanujan1915}.
\end{theorem}

\textit{Proof.} Consider the argument value $\ln(i(2-\sqrt{3}))$. Let the radical conjugate of this argument value be the integration limit
$x_{1}=\ln(i(2+\sqrt{3}))$ for the fundamental form of the gemini-function, where $a=+1$. Next, we calculate the corresponding
other integration limit such that $x_{2}=\ln(\frac{i(2+\sqrt{3})+1}{i(2+\sqrt{3})+1})=\ln(\frac{\sqrt{3}-1}{2})=e^{-\frac{i\pi}{6}}$.
The absolute value of the argument is such that $|e^{-\frac{i\pi}{6}}|=1$ and the multiplier in the exponent is rational, i.e.,
$-\frac{\pi}{6}$. Hence, we can easily prove the above identity. Next, we put all the initials in the five-term gemini-identity
\eqref{finalfivegemini}, where the shape factor $a=+1$. Thus, we obtain that

\begin{multline*}
 \operatorname{Li_{2}}\left(-e^{\frac{i\pi}{6}}\right)-\operatorname{Li_{2}}\left(e^{\frac{i\pi}{6}}\right)+\frac{\pi^{2}}{4} - \\ 
\ln\left(e^{-\frac{i\pi}{6}}\right)\ln\left(\frac{1}{i(2+\sqrt{3})}\right)=-\operatorname{Li_{2}}\left(-\frac{1}{i(2+\sqrt{3})}\right)+ 
\operatorname{Li_{2}}\left(\frac{1}{i(2+\sqrt{3})}\right). 
\end{multline*}

Consequently, we have that 
$$\operatorname{Li_{2}}(-e^{\frac{i\pi}{6}})-\operatorname{Li_{2}}(e^{\frac{i\pi}{6}})+\frac{\pi^{2}}{6}+\frac{1}{6}i\pi \ln(2+\sqrt{3}) = 
-\operatorname{Li_{2}}(i(2-\sqrt{3}))+\operatorname{Li_{2}}(i(\sqrt{3}-2)). $$
Therefore, we obtain that 
$$\operatorname{Li_{2}}(i(2-\sqrt{3}))-\operatorname{Li_{2}}(i(\sqrt{3}-2))=-\frac{\pi^{2}}{6}-\frac{1}{6}i\pi \ln(2+\sqrt{3})+
\operatorname{Li_{2}}(e^{\frac{i\pi}{6}})-\operatorname{Li_{2}}(-e^{\frac{i\pi}{6}}). $$
The right-hand side may be expressed in terms of $\psi_1$ in a routine way, giving us an equivalent formulation of the desired result. 

We may also derive another result from the work of Campbell \cite{Campbell202122}, as below. Again, this was given in an equivalent
way by Ramanujan, whereas Campbell provided a new Fourier--Legendre-based derivation. 

\begin{theorem}
(Ramanujan, 1915) The evaluation 
$$\operatorname{Li_{2}}(i(\sqrt{2}-1))-\operatorname{Li_{2}}(i(1-\sqrt{2}))=
\frac{i[\sqrt{2}(\psi_{1}(\frac{1}{8})+\psi_{1}(\frac{3}{8}))+8\pi\ln(\sqrt{2}-1)-4\sqrt{2}\pi^{2}]}{32}$$
holds. 
\end{theorem}

\textit{Proof.} This argument value $\ln(i(\sqrt{2}-1))$ plays the key role now. Let the other integration limit be such that
$x_{1}=\ln(i(\sqrt{2}+1))$, which is the radical conjugate of the argument value of the identity under investigation.
We apply again the five-term identity obtained from the fundamental form, i.e., $a=+1$. The other integration limit is such that
$x_{2}=\ln(\frac{i(\sqrt{2}+1)+1}{i(\sqrt{2}+1)-1})=\ln(e^{-\frac{i\pi}{4}})$. The absolute value of the argument of the
$x_{2}$-term is again one and the multiplicative coefficient $\frac{1}{4}$ in the exponent is rational. Thus, we can build the
following five-term identity with these integration limits. We then obtain that 

\begin{multline*}
\operatorname{Li_{2}}(-e^{\frac{i\pi}{4}})-\operatorname{Li_{2}}(e^{\frac{i\pi}{4}}) + 
\frac{\pi^{2}}{4}-\ln(e^{-\frac{i\pi}{4}})\ln(i(\sqrt{2}+1)) = \\ 
-\operatorname{Li_{2}}\left(-\frac{1}{i(\sqrt{2}+1)}\right)+\operatorname{Li_{2}}\left(\frac{1}{i(\sqrt{2}+1)}\right).
\end{multline*}
Consequently, we have that 
\begin{multline*}
\operatorname{Li_{2}}(-e^{\frac{i\pi}{4}})-\operatorname{Li_{2}}(e^{\frac{i\pi}{4}})+\frac{\pi^{2}}{4}-\frac{\pi^{2}}{8}+\frac{1}{4}i\pi \ln(\sqrt{2}+1)= \\ 
-\operatorname{Li_{2}}(i(\sqrt{2}-1))+\operatorname{Li_{2}}(-i(1-\sqrt{2})). 
\end{multline*}
So, we find that
\begin{multline*}
\operatorname{Li_{2}}(i(\sqrt{2}-1))-\operatorname{Li_{2}}(i(1-\sqrt{2})) = \\ 
\operatorname{Li_{2}}(e^{\frac{i\pi}{4}})-
\operatorname{Li_{2}}(-e^{\frac{i\pi}{4}})-\frac{\pi^{2}}{8}-\frac{1}{4}i\pi \ln(\sqrt{2}+1). 
\end{multline*}
By expressing the right-hand side in terms of $\psi_1$, we obtain an equivalent version of the desired result. 
 
\subsection{On the fundamental form of a gemini function}

It is possible to choose an arbitrary complex number to be the initial value to generate a two-term identity. The only requirement
is that the absolute value of the argument must be one and the angle must be a rational multiple of $\pi$ number, as earlier stated.
Let us set, e.g. the lower integration limit such that $x_{1}=\ln(e^{\frac{i\pi}{5}})$. Hence, the upper limit is given by
$x_{2}=\ln\left(\frac{e^{\frac{i\pi}{5}}+1}{e^{\frac{i\pi}{5}}-1}\right)=\ln(-i\sqrt[4]{5}\phi^{\frac{3}{2}})$. Next, we set these initials
again in the gemini five-term identity \eqref{finalfivegemini}, with the scale factor $a=+1$. This gives us that 
\begin{eqnarray}
\operatorname{Li_{2}}\left(\frac{i}{\sqrt[4]{5}\phi^{\frac{3}{2}}}\right) - 
 \operatorname{Li_{2}}\left(-\frac{i}{\sqrt[4]{5}\phi^{\frac{3}{2}}}\right) 
\nonumber
\end{eqnarray}
may be expressed as 
\begin{multline*}
\frac{i\sqrt{\phi^{2}+1}}{200}\biggl\{\frac{1}{\phi}\psi_{1}\left(\frac{1}{10}\right)+\left(\frac{4}{\phi}+1\right)\psi_{1}\left(\frac{1}{5}\right)+
\psi_{1}\left(\frac{3}{10}\right)+\left(\frac{1}{\phi}-4\right)\psi_{1}\left(\frac{2}{5}\right) + \\ 
 \left(4-\frac{1}{\phi}\right)\psi_{1}\left(\frac{3}{5}\right)-
\psi_{1}\left(\frac{7}{10}\right)-\left(\frac{4}{\phi}+1\right)\psi_{1}\left(\frac{4}{5}\right)-
\frac{1}{\phi}\psi_{1}\left(\frac{9}{10}\right)\biggl\}-\frac{i\pi}{20}\ln(5\phi^{6}). 
\end{multline*}

\subsection{Applying Theorem \ref{threetermconstant} with the third term fixed by $\phi$}

We can derive also single and two-term identities in the complex domain by applying three-term cancellation identities as in Theorem
\ref{mainthreeterm} and (19). Next, we calculate an exact complex valued single-term identity with the aid of Theorem \ref{threetermconstant}
by setting the third term to be equal to $\phi$. Hence, we can write

\vskip 0.1in

$\operatorname{Li_{2}}(\frac{a+1}{a^{2}+a+1})=\operatorname{Li_{2}}(\phi) \Rightarrow \frac{a+1}{a^{2}+a+1}=\phi \Rightarrow
a=-\frac{1}{2\phi^{2}}+\frac{i}{2\phi}\sqrt{\phi^{2}+1}$.

\vskip 0.1in

By inserting this obtained formula related to $a$ into the other terms, we get

\vskip 0.1in

$\operatorname{Li_{2}}(\frac{1}{2}+\frac{i}{2\phi^{2}}\sqrt{\phi^{2}+1})+\operatorname{Li_{2}}(\frac{1}{2}\phi^{2}-\frac{i}{2}\sqrt{\phi^{2}+1})
-\operatorname{Li_{2}}(\phi)+$
\vskip 0.01in
$\ln(\frac{1}{2}-\frac{i\phi}{2}\sqrt{\phi^{2}+1})\ln(\frac{1}{2}\phi^{2}+\frac{i}{2}\sqrt{\phi^{2}+1})=0$.

\vskip 0.1in

We know the exact value of $\operatorname{Li_{2}}(\phi)=\frac{7\pi^{2}}{30}+\frac{1}{2}\ln^{2}(\phi)-i\pi\ln(\phi)$, which is inserted into the
formula above. By simplifying the constant terms, we get

\vskip 0.1in

$\operatorname{Li_{2}}(\frac{1}{2}+\frac{i}{2\phi^{2}}\sqrt{\phi^{2}+1})+\operatorname{Li_{2}}(\frac{1}{2}\phi^{2}-\frac{i}{2}\sqrt{\phi^{2}+1})
-\frac{7\pi^{2}}{30}-\frac{1}{2}\ln^{2}(\phi)+i\pi\ln(\phi)+\frac{2\pi^{2}}{25}+$
\vskip 0.01in
$\ln^{2}(\phi)-\frac{1}{5}i\pi\ln(\phi)=0 \Rightarrow$

\vskip 0.1in

$\operatorname{Li_{2}}(\frac{1}{2}+\frac{i}{2\phi^{2}}\sqrt{\phi^{2}+1})+\operatorname{Li_{2}}(\frac{1}{2}\phi^{2}-\frac{i}{2}\sqrt{\phi^{2}+1})
-\frac{23\pi^{2}}{150}+\frac{1}{2}\ln^{2}(\phi)+\frac{4}{5}i\pi\ln(\phi)=0$.

\vskip 0.1in

Next, we apply the reflection identity \eqref{reflectionformula} to the second term to get a new argument, whose absolute value is one.
Hence, we can write

\vskip 0.1in

$\operatorname{Li_{2}}(\frac{1}{2}\phi^{2}-\frac{i}{2}\sqrt{\phi^{2}+1})=$
\vskip 0.01in
$-\operatorname{Li_{2}}(-\frac{1}{2\phi}+\frac{i}{2}\sqrt{\phi^{2}+1})+\frac{\pi^{2}}{6}-\ln(\frac{1}{2}\phi^{2}-
\frac{i}{2}\sqrt{\phi^{2}+1})\ln(-\frac{1}{2\phi}+\frac{i}{2}\sqrt{\phi^{2}+1})=$
\vskip 0.01in
$-\operatorname{Li_{2}}(e^{\frac{3i\pi}{5}})+\frac{\pi^{2}}{6}-\frac{3\pi^{2}}{25}-\frac{3}{5}i\pi\ln(\phi)=
-\operatorname{Li_{2}}(e^{\frac{3i\pi}{5}})+\frac{7\pi^{2}}{150}-\frac{3}{5}i\pi\ln(\phi)$.

\vskip 0.1in

By inserting the above formula into the original identity, we get

\vskip 0.1in

$\operatorname{Li_{2}}(\frac{1}{2}+\frac{i}{2\phi^{2}}\sqrt{\phi^{2}+1})-\operatorname{Li_{2}}(e^{\frac{3i\pi}{5}})-
\frac{8\pi^{2}}{75}+\frac{1}{2}\ln^{2}(\phi)+\frac{1}{5}i\pi\ln(\phi)=0$.

\vskip 0.1in

The term $-\operatorname{Li_{2}}(e^{\frac{3i\pi}{5}})$ can be expressed with the aid of trigamma functions, because the absolute
value of the argument is one and the multiplicative coefficient in the exponent is rational. After a workable manipulation,
we get the final single value representation, as shown in (36).

\begin{eqnarray}
\operatorname{Li_{2}}\left(\frac{1}{2}+i\frac{\sqrt{\phi^{2}+1}}{2\phi^{2}}\right)=\operatorname{Li_{2}}\left(\frac{e^{\frac{i\pi}{5}}}{\phi}\right)=
\frac{19\pi^{2}}{300}-\frac{1}{2}\ln^{2}(\phi)+
\nonumber
\end{eqnarray}
\begin{eqnarray}
i\biggl\{\frac{\sqrt{\phi^{2}+1}}{200}\left[\psi_{1}\left(\frac{1}{10}\right)+\psi_{1}\left(\frac{2}{5}\right)-
\psi_{1}\left(\frac{3}{5}\right)-\psi_{1}\left(\frac{9}{10}\right) \right]
\nonumber
\end{eqnarray}
\begin{eqnarray}
+\frac{\sqrt{\phi^{2}+1}}{200\phi}\left[\psi_{1}\left(\frac{4}{5}\right)+\psi_{1}\left(\frac{7}{10}\right)-
\psi_{1}\left(\frac{3}{10}\right)-\psi_{1}\left(\frac{1}{5}\right)\right]-\frac{1}{5}\pi\ln(\phi)\biggl\} 
\end{eqnarray}

\subsection{Case I: Applying Theorem \ref{mainthreeterm} for deriving a real part for a complex valued dilogarithm}

Next, we derive an exact value for a real part of a complex dilogarithm. 

\begin{theorem}
The closed-form evaluation
\begin{multline*}
\mathfrak{Re}\biggl\{\operatorname{Li_{2}}\left(\frac{3+i\sqrt{7}}{2}\right)\biggl\}=-\frac{\pi^{2}}{24}-\frac{1}{4}\ln^{2}(2)
+\frac{1}{2}\pi\operatorname{arctan}\left(\frac{\sqrt{7}}{5}\right) + \\ 
\frac{1}{2}\operatorname{arctan}\left(\frac{\sqrt{7}}{3}\right)\operatorname{arctan}(\sqrt{7}). 
\end{multline*}
holds. 
\end{theorem}

\textit{Proof.} We apply Theorem \ref{mainthreeterm} in the following manner. Let the argument of the first term be 2 in such a way that
$\operatorname{Li_{2}}(\frac{a-1}{a^{2}})=\operatorname{Li_{2}}(2) \Rightarrow \frac{a-1}{a^{2}}=2 \Rightarrow
a=\frac{1\pm i\sqrt{7}}{4}$. Next, we insert the root with the positive imaginary part of $a$ into Theorem \ref{mainthreeterm} 
and we can write

\vskip 0.1in

$\operatorname{Li_{2}}(2)+\operatorname{Li_{2}}(\frac{3+i\sqrt{7}}{2})-\operatorname{Li_{2}}(\frac{-1+i\sqrt{7}}{2})+\left[\frac{1}{2}\ln(2)
 +i\operatorname{arctan}(\sqrt{7})\right]i\pi=0 \Rightarrow$

\vskip 0.1in

$\frac{\pi^{2}}{4}-i\pi\ln(2)+\operatorname{Li_{2}}(\frac{3+i\sqrt{7}}{2})-\operatorname{Li_{2}}(\frac{-1+i\sqrt{7}}{2})+\frac{1}{2}i\pi\ln(2)
-\pi\operatorname{arctan}(\sqrt{7})=0 \Rightarrow$

\vskip 0.1in

$\operatorname{Li_{2}}(\frac{3+i\sqrt{7}}{2})-\operatorname{Li_{2}}(\frac{-1+i\sqrt{7}}{2})+\frac{\pi^{2}}{4}-\frac{1}{2}i\pi\ln(2)
-\pi\operatorname{arctan}(\sqrt{7})=0.$
Consequently, we have that 
\begin{multline*}
\operatorname{Li_{2}}\left(\frac{3+i\sqrt{7}}{2}\right)+\operatorname{Li_{2}}\left(\frac{3-i\sqrt{7}}{2}\right)+\frac{\pi^{2}}{12} - \\ 
\frac{1}{2}i\pi\ln(2)
-\pi\operatorname{arctan}(\sqrt{7})+\ln\left(\frac{-1+i\sqrt{7}}{2}\right)\ln\left(\frac{3-i\sqrt{7}}{2}\right)=0.
\end{multline*}
This gives us an equivalent version of the desired result.

\subsection{Applying Theorem \ref{mainthreeterm} in the derivation of the 
 real part for a complex-valued dilogarithm}
 We can also derive a further 
 exact expression for the 
 real part for a complex-valued dilogarithm in a similar manner, 
 by applying Theorem \ref{mainthreeterm}.
 In this direction, 
 we set the third term in Theorem \ref{mainthreeterm} to be equal to 2 such that
$\operatorname{Li_{2}}(\frac{a}{a^{2}-a+1})=\operatorname{Li_{2}}(2) \Rightarrow \frac{a}{a^{2}-a+1}=2 \Rightarrow a=\frac{3+i\sqrt{7}}{4}$.
By inserting this obtained root $a$ and the value 2 into Theorem \ref{mainthreeterm}, and hence
\begin{multline*}
\operatorname{Li_{2}}\left( \frac{3-i\sqrt{7}}{2} \right) + 
 \operatorname{Li_{2}}\left( 
 \frac{5+i\sqrt{7}}{8} \right) 
 - \frac{\pi^{2}}{4}+\frac{1}{2}\ln^{2}(2)+\frac{1}{2}i\pi\ln(2) \\ 
 +\operatorname{arctan}\left( 
 \frac{\sqrt{7}}{3} \right) 
 \operatorname{arctan}(\sqrt{7})+
i\ln(2) 
 \operatorname{arctan}\left(
 \frac{\sqrt{7}}{3} \right) = 0. 
\end{multline*}
 A simplification then gives us that 
\begin{multline*}
 \mathfrak{Re}\biggl\{\operatorname{Li_{2}}\left(\frac{5+i\sqrt{7}}{8}\right)\biggl\}=\frac{\pi^{2}}{4}-\frac{1}{2}\ln^{2}(2) - \\ 
 \operatorname{arctan}\left(\frac{\sqrt{7}}{3}\right)\operatorname{arctan}(\sqrt{7})
-\mathfrak{Re}\biggl\{\operatorname{Li_{2}}\left(\frac{3-i\sqrt{7}}{2}\right)\biggl\}, 
\end{multline*}
 and we obtain 
 the final representation for the exact real part term by substituting the previous result into the formula above.
Hence, we can write
\begin{multline*}
\mathfrak{Re}\biggl\{\operatorname{Li_{2}}\left(\frac{5+i\sqrt{7}}{8}\right)\biggl\}=\frac{7\pi^{2}}{24}-\frac{1}{4}\ln^{2}(2)-
\frac{1}{2}\pi\operatorname{arctan}\left(\frac{\sqrt{7}}{5}\right) - \\ 
\frac{3}{2}\operatorname{arctan}\left(\frac{\sqrt{7}}{3}\right)\operatorname{arctan}\left(\sqrt{7}\right).
\end{multline*}
Similar kinds of exact real part values are also introduced in the paper of Hakimoglu-Brown \cite{Hak_25}.

\subsection{Dilogarithm and the imaginary golden ratio \emph{$\phi_{i}$}}

The imaginary golden ratio is given by $\phi_{i}=\frac{1+i\sqrt{3}}{2}=e^{\frac{i\pi}{3}}$ and $|\phi_{i}|=1$. It is the root of the equation
$x^{2}-x+1=0$. It has analogous algebraic properties as like the real golden ratio $\phi$ has. Among other things, it satisfies the following
formulae, $\phi_{i}=1-\frac{1}{\phi_{i}}$ and $\phi_{i}^{n}=\phi_{i}^{n-1}-\phi_{i}^{n-2}$.The nested radical representation for the imaginary golden
ratio is such that $\phi_{i}=\sqrt{-1+\sqrt{-1+\sqrt{-1+...\quad}}}$. It also gives particularly simple results in these two following formulae,
$\sin(i\ln(\phi_{i}))=-\frac{1}{2}\sqrt{3}$ and $\sin[\frac{\pi}{2}-i\ln(\phi_{i})]=\frac{1}{2}$. Since the absolute value of $\phi_{i}$ is equal
to 1 and arg$(\phi_{i})$ is rational multiple of $\pi$, it is a trivial task to evaluate its dilogarithm, which is given by
$\operatorname{Li_{2}}(\phi_{i})=\frac{\pi^{2}}{36}-i\operatorname{Cl_{2}}(\frac{\pi}{3})$. The term $\operatorname{Cl_{2}}(\frac{\pi}{3})$
stands for the Clausen function at $\theta=\frac{\pi}{3}$, whose value is also referred to as \emph{Gieseking's constant}
\cite{Adams1998} \cite[pp.\ 232--233]{Finch2003}. Gieseking's constant can also be expressed in terms of the trigamma function as 
$\operatorname{\mathcal{G}_{GI}}=\frac{9-\psi_{1}(\frac{2}{3})+\psi_{1}(\frac{4}{3})}{4\sqrt{3}}\approx1.014943$.

\vskip 0.1in

Let us once more return to the results derived in the paper of Campbell \cite{Campbell202122}. The following identity is interesting,
since it is related to the imaginary golden ratio $\phi_{i}$.

\vskip 0.1in

$\operatorname{Li_{2}}(\frac{i}{\sqrt{3}})-\operatorname{Li_{2}}(-\frac{i}{\sqrt{3}})=
i\left[\frac{3\psi_{1}(\frac{1}{6})+15\psi_{1}(\frac{1}{3})-6\sqrt{3}\pi\ln(3)-16\pi^{2}}{36\sqrt{3}}\right]$

\vskip 0.1in

By setting, e.g. the integration limit in such a way that $x_{1}=i\sqrt{3}$. Hence, the other limit is given by
$x_{2}=\ln(\frac{i\sqrt{3}+1}{i\sqrt{3}-1})=\ln(\frac{1-i\sqrt{3}}{2})=\ln(e^{-\frac{i\pi}{3}})=\ln(\bar{\phi_{i}})$.
Here, we apply again the five-term gemini-identity \eqref{finalfivegemini} with $a=+1$. Hence, we get

\vskip 0.1in

$\operatorname{Li_{2}}(-\frac{1}{i\sqrt{3}})-\operatorname{Li_{2}}(\frac{1}{i\sqrt{3}})+\frac{\pi^{2}}{4}-\ln(i\sqrt{3})\ln(\bar{\phi_{i}})
+\operatorname{Li_{2}}(-\frac{1}{\bar{\phi_{i}}})-\operatorname{Li_{2}}(\frac{1}{\bar{\phi_{i}}})=0 \Rightarrow $

\vskip 0.1in

For clarity: $\frac{1}{\bar{\phi_{i}}}=\phi_{i}=e^{\frac{i\pi}{3}}, \quad |e^{\frac{i\pi}{3}}|=1 \Rightarrow $

\vskip 0.1in

$\operatorname{Li_{2}}(\frac{i}{\sqrt{3}})-\operatorname{Li_{2}}(-\frac{i}{\sqrt{3}})+\frac{\pi^{2}}{12}+\frac{i\pi\ln(3)}{6}+
\operatorname{Li_{2}}(-e^{\frac{i\pi}{3}})-\operatorname{Li_{2}}(e^{\frac{i\pi}{3}})=0\Rightarrow $

\vskip 0.1in

$\operatorname{Li_{2}}(\frac{i}{\sqrt{3}})-\operatorname{Li_{2}}(-\frac{i}{\sqrt{3}})=-\frac{\pi^{2}}{12}-\frac{i\pi\ln(3)}{6}+
\operatorname{Li_{2}}(-\phi_{i})-\operatorname{Li_{2}}(\phi_{i})=0 \Rightarrow $

\vskip 0.1in

$\operatorname{Li_{2}}(\frac{i}{\sqrt{3}})-\operatorname{Li_{2}}(-\frac{i}{\sqrt{3}})=i[\frac{\psi_{1}(\frac{1}{6})+5\psi_{1}(\frac{1}{3})
-5\psi_{1}(\frac{2}{3})-\psi_{1}(\frac{5}{6})}{24\sqrt{3}}-\frac{\pi\ln(3)}{6}] \Rightarrow $

\begin{eqnarray}
\operatorname{Li_{2}}\left(\frac{i}{\sqrt{3}}\right)-\operatorname{Li_{2}}\left(-\frac{i}{\sqrt{3}}\right)=
i\left[\frac{5}{3}\operatorname{\mathcal{G}_{GI}}-\frac{1}{6}\pi\ln(3)\right].
\end{eqnarray}

The constant term in our formula differs again from Campbell's result \cite{Campbell202122}, but those are numerically exactly
the same. By setting an equal sign between these two constant terms, we get a nice trigamma-identity of the form

\vskip 0.1in

$i\left[\frac{3\psi_{1}(\frac{1}{6})+15\psi_{1}(\frac{1}{3})-6\sqrt{3}\pi\ln(3)-16\pi^{2}}{36\sqrt{3}}\right]
=i\left[\frac{\psi_{1}(\frac{1}{6})+5\psi_{1}(\frac{1}{3})-5\psi_{1}(\frac{2}{3})-\psi_{1}(\frac{5}{6})}{24\sqrt{3}}-
\frac{\pi\ln(3)}{6}\right] \Rightarrow$

\vskip 0.1in

$\psi_{1}(\frac{1}{6})+5\psi_{1}(\frac{1}{3})+5\psi_{1}(\frac{2}{3})+\psi_{1}(\frac{5}{6})=\frac{32\pi^{2}}{3}.$

\vskip 0.1in

Next, we derive the imaginary part for $\operatorname{Li_{2}}(\frac{1}{2}\phi_{i})$. First, we have to build two separate identities,
where the terms $\operatorname{Li_{2}}(\frac{1}{2}\phi_{i})$ and $\operatorname{Li_{2}}(\frac{1}{2}\bar{\phi_{i}})$ are connected to the
$\operatorname{Li_{2}}(-i\sqrt{3})$-term. The calculation related to the first identity goes as follows:

\vskip 0.1in

$\operatorname{Li_{2}}(-i\sqrt{3})=\operatorname{Li_{2}}(\frac{1}{1+i\sqrt{3}})-
\frac{\pi^{2}}{6}+\frac{1}{2}\ln(1+i\sqrt{3})\ln[\frac{1+i\sqrt{3}}{(i\sqrt{3})^{2}}]=$
\vskip 0.01in
$\operatorname{Li_{2}}(\frac{1-i\sqrt{3}}{4})-\frac{\pi^{2}}{18}-\frac{1}{2}\ln(2)\ln(\frac{3}{2})-\frac{1}{6}i\pi\ln(6)=$
\vskip 0.01in
$\operatorname{Li_{2}}(\frac{1}{2}\bar{\phi_{i}})-\frac{\pi^{2}}{18}-\frac{1}{2}\ln(2)\ln(\frac{3}{2})-\frac{1}{6}i\pi\ln(6) \Rightarrow$
\vskip 0.01in
$\operatorname{Li_{2}}(-i\sqrt{3})=\operatorname{Li_{2}}(\frac{1}{2}\bar{\phi_{i}})-\frac{\pi^{2}}{18}-
\frac{1}{2}\ln(2)\ln(\frac{3}{2})-\frac{1}{6}i\pi\ln(6). $

\vskip 0.1in

The other connection is obtained from the five-term gemini-identity \eqref{finalfivegemini} in a following way. Let us define the
integration limits in such a way that $x_{1}=\ln(\frac{1+i\sqrt{3}}{2})=\ln(\phi_{i})$ and $x_{2}=\ln(2)$. Hence, the formula for
the shape factor is such that $\frac{2+a}{2-1}=\phi_{i} \Rightarrow a=\phi_{i}-2$. By setting the initial values into the five-term
identity \eqref{finalfivegemini}, we get

\vskip 0.1in

$\operatorname{Li_{2}}(\frac{2-\phi_{i}}{\phi_{i}})-\operatorname{Li_{2}}(\frac{1}{\phi_{i}})-\operatorname{Li_{2}}(2-\phi_{i})+\frac{\pi^{2}}{6}-
\ln(\phi_{i})\ln(2)=-\operatorname{Li_{2}}(\frac{2-\phi_{i}}{2})+\operatorname{Li_{2}}(\frac{1}{2}) \Rightarrow$

\vskip 0.1in

$\operatorname{Li_{2}}(-i\sqrt{3})-\operatorname{Li_{2}}(\bar{\phi_{i}})-\operatorname{Li_{2}}(1+\bar{\phi_{i}})+\frac{\pi^{2}}{6}-\ln(\phi_{i})\ln(2)=
-\operatorname{Li_{2}}(\frac{1}{2}+\frac{1}{2}\bar{\phi_{i}})+\operatorname{Li_{2}}(\frac{1}{2}) \Rightarrow$

\vskip 0.1in

$\operatorname{Li_{2}}(-i\sqrt{3})-\operatorname{Li_{2}}(\bar{\phi_{i}})+\operatorname{Li_{2}}(-\bar{\phi_{i}})+
\ln(-\bar{\phi_{i}})\ln(1+\bar{\phi_{i}})-\ln(\phi_{i})\ln(2)=$
\vskip 0.01in
$\operatorname{Li_{2}}(\frac{1}{2}\phi_{i})-\frac{\pi^{2}}{6}+\ln(\frac{1}{2}\phi_{i})\ln(\frac{1}{2}+\frac{1}{2}\bar{\phi_{i}})+
\operatorname{Li_{2}}(\frac{1}{2}) \Rightarrow$

\vskip 0.1in

$\operatorname{Li_{2}}(-i\sqrt{3})=\operatorname{Li_{2}}(\frac{1}{2}\phi_{i})+\operatorname{Li_{2}}(\bar{\phi_{i}})-
\operatorname{Li_{2}}(-\bar{\phi_{i}})-\ln(-\bar{\phi_{i}})\ln(1+\bar{\phi_{i}})+\ln(\phi_{i})\ln(2)$
\vskip 0.01in
$-\frac{\pi^{2}}{6}+\ln(\frac{1}{2}\phi_{i})\ln(\frac{1}{2}+\frac{1}{2}\bar{\phi_{i}})+\operatorname{Li_{2}}(\frac{1}{2}).$

\vskip 0.1in

Our next task is to combine these two auxiliary equations as follows:

\vskip 0.1in

$\operatorname{Li_{2}}(\frac{1}{2}\bar{\phi_{i}})-\frac{\pi^{2}}{18}-\frac{1}{2}\ln(2)\ln(\frac{3}{2})-\frac{1}{6}i\pi\ln(6)=
\operatorname{Li_{2}}(\bar{\phi_{i}})-\operatorname{Li_{2}}(-\bar{\phi_{i}})$
\vskip 0.01in
$-\ln(-\bar{\phi_{i}})\ln(1+\bar{\phi_{i}})+\ln(\phi_{i})\ln(2)+\operatorname{Li_{2}}(\frac{1}{2}\phi_{i})-
\frac{\pi^{2}}{6}+\ln(\frac{1}{2}\phi_{i})\ln(\frac{1}{2})+\frac{1}{2}\bar{\phi_{i}}+\operatorname{Li_{2}}(\frac{1}{2})$

\vskip 0.1in

 After a simplification, we get

\vskip 0.1in

$\operatorname{Li_{2}}(\frac{1}{2}\phi_{i})-\operatorname{Li_{2}}(\frac{1}{2}\bar{\phi_{i}})=\frac{\pi^{2}}{12}-\frac{1}{3}i\pi\ln(2)
+\operatorname{Li_{2}}(-\bar{\phi_{i}})-\operatorname{Li_{2}}(\bar{\phi_{i}})\Rightarrow$

\vskip 0.1in

$\operatorname{Li_{2}}(\frac{1}{2}\phi_{i})-\operatorname{Li_{2}}(\frac{1}{2}\bar{\phi_{i}})=-\frac{1}{3}i\pi\ln(2)+
\frac{i}{24\sqrt{3}}[\psi_{1}(\frac{1}{6})+5\psi_{1}(\frac{1}{3})-5\psi_{1}(\frac{2}{3})-\psi_{1}(\frac{5}{6})] \Rightarrow$

\vskip 0.1in

$\operatorname{Li_{2}}(\frac{1}{2}\phi_{i})-\operatorname{Li_{2}}(\frac{1}{2}\bar{\phi_{i}})=i[\frac{5}{3}\operatorname{\mathcal{G}_{GI}}-
\frac{1}{3}\pi\ln(2)] \Rightarrow$

\begin{eqnarray}\label{38}
\mathfrak{Im}\biggl\{\operatorname{Li_{2}}\left(\frac{1}{2}\phi_{i}\right)\biggl\}=\frac{5}{6}\operatorname{\mathcal{G}_{GI}}-\frac{1}{6}\pi\ln(2).
\end{eqnarray}

We can exploit the above result to derive also the exact value for the $\mathfrak{Im}\bigl\{\operatorname{Li_{2}}(2\phi_{i})\bigl\}$. We begin
this evaluation by applying the relation below.

\vskip 0.1in

$\operatorname{Li_{2}}(2\phi_{i})=\operatorname{Li_{2}}(1+i\sqrt{3})=-\operatorname{Li_{2}}(-i\sqrt{3})+\frac{\pi^{2}}{6}-
\ln(1+i\sqrt{3})\ln(-i\sqrt{3}).$

\vskip 0.1in

On the other hand, we can write
$$
 \operatorname{Li_{2}}(-i\sqrt{3})=\operatorname{Li_{2}}(\frac{1}{2}\bar{\phi_{i}})-\frac{\pi^{2}}{18}-\frac{1}{2}\ln(2)\ln(\frac{3}{2})-
\frac{1}{6}i\pi\ln(6)$$ 
 and $\mathfrak{Im}\bigl\{\operatorname{Li_{2}}(\frac{1}{2}\bar{\phi_{i}})\bigl\}=
-\frac{5}{6}\operatorname{\mathcal{G}_{GI}}+\frac{1}{6}\pi\ln(2).$

\vskip 0.1in

By putting all together, we get

\begin{eqnarray}\label{39}
\mathfrak{Im}\biggl\{\operatorname{Li_{2}}\left(2\phi_{i}\right)\biggl\}=\frac{5}{6}\operatorname{\mathcal{G}_{GI}}+\frac{1}{2}\pi\ln(2).
\end{eqnarray}

\section{The unresolved $\operatorname{Li_{2}}\left(-\frac{1}{2}\right)$}
 The exact value of $\operatorname{Li_{2}}(-\frac{1}{2})$ is a great mystery, since this term does 
 reveal almost nothing about itself. It is known that it has a close connection to the number three in the sets of natural and rational 
 numbers. The term $\operatorname{Li_{2}}(-\frac{1}{2})$ is kind of a self-destructive entity. When it appears in the equation, 
 the very next step the same amount of these terms appear in the equation
with opposite signs and they cancel each other out. This can be seen, for example, in Ramanujan's identities introduced earlier. By carrying out
the simplification to the end with these identities, final outcomes are that all the $\operatorname{Li_{2}}(-\frac{1}{2})$ terms cancel each
other out. This term seems to be included in the calculation only in a supporting role. Its purpose is simply to make the equations computationally
true. The term $\operatorname{Li_{2}}(-\frac{1}{2})$ behaves somewhat like a catalyst in a chemical reaction without participating in the end
result itself. It always disappears from the stage before the performance itself ends, preserving its mystery. In fact, the term
$\operatorname{Li_{2}}(-\frac{1}{2})$ behaves in a completely different way with irrational and complex numbers. So next we will examine its
connections with these numbers. It is an easy task to generate identity formulae for $\operatorname{Li_{2}}(-\frac{1}{2})$ by applying five-term
or three-term cancellation gemini-identities in such a way that the representation of $\operatorname{Li_{2}}(-\frac{1}{2}$) contains two other
dilogarithm terms. We can derive a couple of three-term identities in the real domain for $\operatorname{Li_{2}}(-\frac{1}{2})$, which are listed
below with the initial values needed to build the particular identity.

\vskip 0.1in

\textbf{1.} $\gemini_{+\frac{1}{\sqrt{2}}}(x)$; $a=+\frac{1}{\sqrt{2}}$, $x_{1}=\ln(\sqrt{2})$ and $x_{2}=\ln\left(\frac{6+3\sqrt{2}}{2}\right)$:

\begin{eqnarray}
\operatorname{Li_{2}}\left(-\frac{1}{2}\right)=\frac{\pi^{2}}{24}-\operatorname{Li_{2}}\left(\frac{1+\sqrt{2}}{3}\right)-
\operatorname{Li_{2}}\left(\frac{1-\sqrt{2}}{3}\right)
\end{eqnarray}
\begin{eqnarray}
-\ln\left(\frac{2-\sqrt{2}}{3}\right)\ln\left(\frac{1+\sqrt{2}}{3}\right)-\frac{1}{2}\ln(2)\ln\left(\frac{2\sqrt{2}-2}{3}\right)
\nonumber
\end{eqnarray}

\textbf{2.} $\gemini_{+\frac{1}{\sqrt{3}}}(x)$; $a=+\frac{1}{\sqrt{3}}$, $x_{1}=\ln(\sqrt{3})$ and $x_{2}=\ln\left(\frac{6+2\sqrt{3}}{3}\right)$:

\begin{eqnarray}
\operatorname{Li_{2}}\left(-\frac{1}{2}\right)=\frac{2}{3}\operatorname{Li_{2}}\left(\frac{1-\sqrt{3}}{4}\right)
+\frac{2}{3}\operatorname{Li_{2}}\left(\frac{1+\sqrt{3}}{4}\right)-\frac{\pi^{2}}{9}+\frac{5}{6}\ln^{2}(2)-\frac{1}{3}\ln^{2}(3) \quad
\end{eqnarray}
\begin{eqnarray} 
 +\frac{1}{6}\ln(3)\ln\left(\frac{3}{4}\right)+\frac{4}{3}\ln(2)\ln\left(\frac{3}{2}\right)+
 \frac{1}{3}\ln\left(\frac{16-8\sqrt{3}}{3}\right)\ln\left(\frac{6+2\sqrt{3}}{3}\right) 
\nonumber
\end{eqnarray}

\textbf{3.} $\gemini_{+2}(x)$; $a=+2$, $x_{1}=\ln(\sqrt{2})$ and $x_{2}=\ln(4+3\sqrt{2})$: 
\begin{multline*}
\operatorname{Li_{2}}\left(-\frac{1}{2}\right)=\operatorname{Li_{2}}\left(\frac{3\sqrt{2}-4}{2}\right)- \\ 
 \operatorname{Li_{2}}\left(4-3\sqrt{2}\right)-\frac{\pi^{2}}{8}-\frac{5}{8}\ln^{2}(2)+\frac{1}{2}\ln(2)\ln(4+3\sqrt{2}).
\end{multline*}

\textbf{4.} $\gemini_{+3}(x)$; $a=+3$, $x_{1}=\ln(\sqrt{3})$ and $x_{2}=\ln(3+2\sqrt{3})$:

\begin{eqnarray}
\operatorname{Li_{2}}\left(-\frac{1}{2}\right)=\frac{2}{3}\operatorname{Li_{2}}(3-2\sqrt{3})-
\frac{2}{3}\operatorname{Li_{2}}\left(\frac{2\sqrt{3}-3}{3}\right)
\end{eqnarray}
\begin{eqnarray}
+\frac{1}{12}\ln^{2}(3)-\frac{1}{2}\ln^{2}(2)+\ln(2)\ln(3)-\frac{1}{3}\ln(3)\ln(3+2\sqrt{3})
\nonumber
\end{eqnarray}

We found two three-term identities, which are related to $\phi^{4}$.

\vskip 0.1in

\textbf{5.} $\gemini_{-\frac{2}{\phi^{4}}}(x)$; $a=-\frac{2}{\phi^{4}}$, $x_{1}=\ln\left(\frac{2}{\phi}\right)$ and $x_{2}=\ln(4)$:
\begin{multline*}
 \operatorname{Li_{2}}\left(-\frac{1}{2}\right)=\frac{1}{2}\operatorname{Li_{2}}\left(\frac{1}{2\phi^{4}}\right) - \\ 
 \frac{1}{2}\operatorname{Li_{2}}\left(\frac{2}{\phi^{4}}\right)-\frac{\pi^{2}}{24}-\frac{1}{4}\ln^{2}(2)+2\ln(2)\ln(\phi) - 
 2\ln^{2}(\phi).
\end{multline*}
 
\textbf{6.} $\gemini_{-\frac{1}{3\phi^{2}}}(x)$; $a=-\frac{1}{3\phi^{2}}$, $x_{1}=\ln\left(\frac{3}{\phi^{2}}\right)$ and
$x_{2}=\ln\left(\frac{8\phi^{2}}{3}\right)$:

\begin{eqnarray}
 \operatorname{Li_{2}}\left(-\frac{1}{2}\right)=\frac{1}{6}\operatorname{Li_{2}}\left(\frac{1}{8\phi^{4}}\right)
 +\frac{1}{6}\operatorname{Li_{2}}\left(\frac{\phi^{4}}{8}\right)-\frac{\pi^{2}}{12}-\frac{4}{3}\ln^{2}(\phi)+\ln^{2}(2)
\end{eqnarray}

Next, we perform an unorthodox maneuver with the five-term identity by setting the shape factor in such a way that $a=-2$. According to the
original definition, the shape factor must be greater or equal to -1. Despite that, we set a following relation between the integration limits
$x_{2}=x_{1}^{2} \Rightarrow \frac{x_{1}+2}{x_{1}-1}=x_{1}^{2} \Rightarrow $. Hence, the roots are $2$ and $\pm e^{-\frac{2i\pi}{3}}$. Next, we select
in such a way that $x_{1}=\ln(e^{-\frac{2i\pi}{3}})$ and respectively $x_{2}=\ln(e^{-\frac{4i\pi}{3}})$. By putting the initial values in the five-term
identity, we get the final three-term formula, which connects $\operatorname{Li_{2}}(-\frac{1}{2})$ to the imaginary
 golden ratio, i.e., $\phi_{i}=e^{\frac{i\pi}{3}}=\frac{1+i\sqrt{3}}{2}$.

\vskip 0.1in

\textbf{7.} $\gemini_{-2}(x)$; $a=-2$, $x_{1}=\ln(e^{-\frac{2i\pi}{3}})$ and $x_{2}=\ln(e^{-\frac{4i\pi}{3}})$:

\begin{eqnarray}\label{73}
 \operatorname{Li_{2}}\left(-\frac{1}{2}\right)=-2\mathfrak{Re}\left\{\operatorname{Li_{2}}(2\phi_{i})\right\}-\frac{1}{2}\ln^{2}(2)
\end{eqnarray} 

We proceed to apply our method to evaluate another relation between $\operatorname{Li_{2}}(-\frac{1}{2})$ and $\phi_{i}$ related
dilogrithm. Let us apply the identity shown in Theorem \ref{threetermconstant} by defining the first term to be equal to
$\frac{1}{3}$. Hence, we get the following equation for solving the variable $a$ in the argument of the first term of
Theorem \ref{threetermconstant}.

\vskip 0.1in

\textbf{8.} $\frac{a}{(a+1)^{2}}=\frac{1}{3} \Rightarrow a=\frac{1+i\sqrt{3}}{2}=\phi_{i}$. Next, we substitute this value in the
all other terms.

\vskip 0.1in

$\operatorname{Li_{2}}\left(\frac{1}{3}\right)+\operatorname{Li_{2}}\left(\frac{1-i\sqrt{3}}{4}\right)
-\operatorname{Li_{2}}\left(\frac{3-i\sqrt{3}}{4}\right)+\ln\left(\frac{3}{2}\right)\ln\left(\frac{3-i\sqrt{3}}{2}\right)=0 \Rightarrow$
\vskip 0.01in
$\operatorname{Li_{2}}\left(-\frac{1}{2}\right)=-\frac{\pi^{2}}{9}+\frac{1}{2}\ln^{2}(2)+\operatorname{Li_{2}}\left(\frac{1+i\sqrt{3}}{4}\right)
+\operatorname{Li_{2}}\left(\frac{1-i\sqrt{3}}{4}\right) \Rightarrow$
\vskip 0.01in
$\operatorname{Li_{2}}\left(-\frac{1}{2}\right)=-\frac{\pi^{2}}{9}+\frac{1}{2}\ln^{2}(2)+\operatorname{Li_{2}}(\frac{1}{2}\phi_{i})+
\operatorname{Li_{2}}(\frac{1}{2}\bar{\phi_{i}}) \Rightarrow$

\begin{eqnarray}\label{74}
 \operatorname{Li_{2}}\left(-\frac{1}{2}\right)=2\mathfrak{Re}\left\{\operatorname{Li_{2}}\left(\frac{1}{2}\phi_{i}\right)\right\}-
 \frac{\pi^{2}}{9}+\frac{1}{2}\ln^{2}(2)
\end{eqnarray} 

Let us do one more evaluation similarly by applying the three-term cancellation identity represented in Theorem \ref{mainthreeterm}. By
setting the first argument equal to $\frac{3}{4}$, then the equation becomes as follows:

\vskip 0.1in

\textbf{9.} $\frac{a-1}{a^{2}}=\frac{3}{4} \Rightarrow$
\vskip 0.01in
$\operatorname{Li_{2}}(\frac{3}{4})+\operatorname{Li_{2}}(-1-2\sqrt{2})-\operatorname{Li_{2}}(2-2\sqrt{2})+
\ln(4)\ln\left(\frac{2-2i\sqrt{2}}{3}\right)=0 \Rightarrow$

\begin{eqnarray}\label{magic}
\operatorname{Li_{2}}\left(-\frac{1}{2}\right)=-\mathfrak{Re}\left\{\operatorname{Li_{2}}\left(2+2i\sqrt{2}\right)\right\}+\frac{\pi^{2}}{12}
-\theta_{m}^{2}-\frac{1}{2}\ln^{2}(2)-\frac{1}{4}\ln^{2}(3).
\end{eqnarray} 

The term $\theta_{m}=\arctan(\sqrt{2})$ in \eqref{magic} is a constant sometimes referred to as the \emph{magic angle}. Equations \eqref{73} and
\eqref{74} do not reveal much about $\operatorname{Li_{2}}(-\frac{1}{2})$, but it seems to have a strong connection to the imaginary golden
ratio. The term $\operatorname{Li_{2}}(-\frac{1}{2})$ can be represented by a single other dilogarithm term with constant terms. A more
detailed study might reveal a pattern between $\operatorname{Li_{2}}(-\frac{1}{2})$ and real parts of particular complex numbers, which
depends on the initial value set for one argument out of three. We have observed that finding the exact value for
$\operatorname{Li_{2}}(-\frac{1}{2})$ has also been dealt by some others, e.g. \cite{Boy_22} . Anyway, we can write a following two-term
identity related to $\phi_{i}$ by applying \eqref{73} and \eqref{74}.

\begin{eqnarray}\label{47}
 \mathfrak{Re}\left\{\operatorname{Li_{2}}\left(\frac{1}{2}\phi_{i}\right)\right\}+
 \mathfrak{Re}\left\{\operatorname{Li_{2}}(2\phi_{i})\right\}=\frac{\pi^{2}}{18}-\frac{1}{2}\ln^{2}(2)
\end{eqnarray} 

By combining \eqref{38}, \eqref{39} and \eqref{47}, we get a nice two-term identity including the complex golden ratio $\phi_{i}$, as
shown below.

\begin{eqnarray}
\operatorname{Li_{2}}\left(\frac{1}{2}\phi_{i}\right)+\operatorname{Li_{2}}\left(2\phi_{i}\right)
 =\frac{\pi^{2}}{18}-\frac{1}{2}\ln^{2}(2)+i\left[\frac{5}{3}\operatorname{\mathcal{G}_{GI}}+\frac{1}{3}\pi\ln(2)\right]
\end{eqnarray}

\section{Geometric properties of gemini functions versus the representation of a dilogarithm}

\begin{figure}[htbp!]
\begin{center}
 \resizebox{10.0cm}{!}
 {\includegraphics{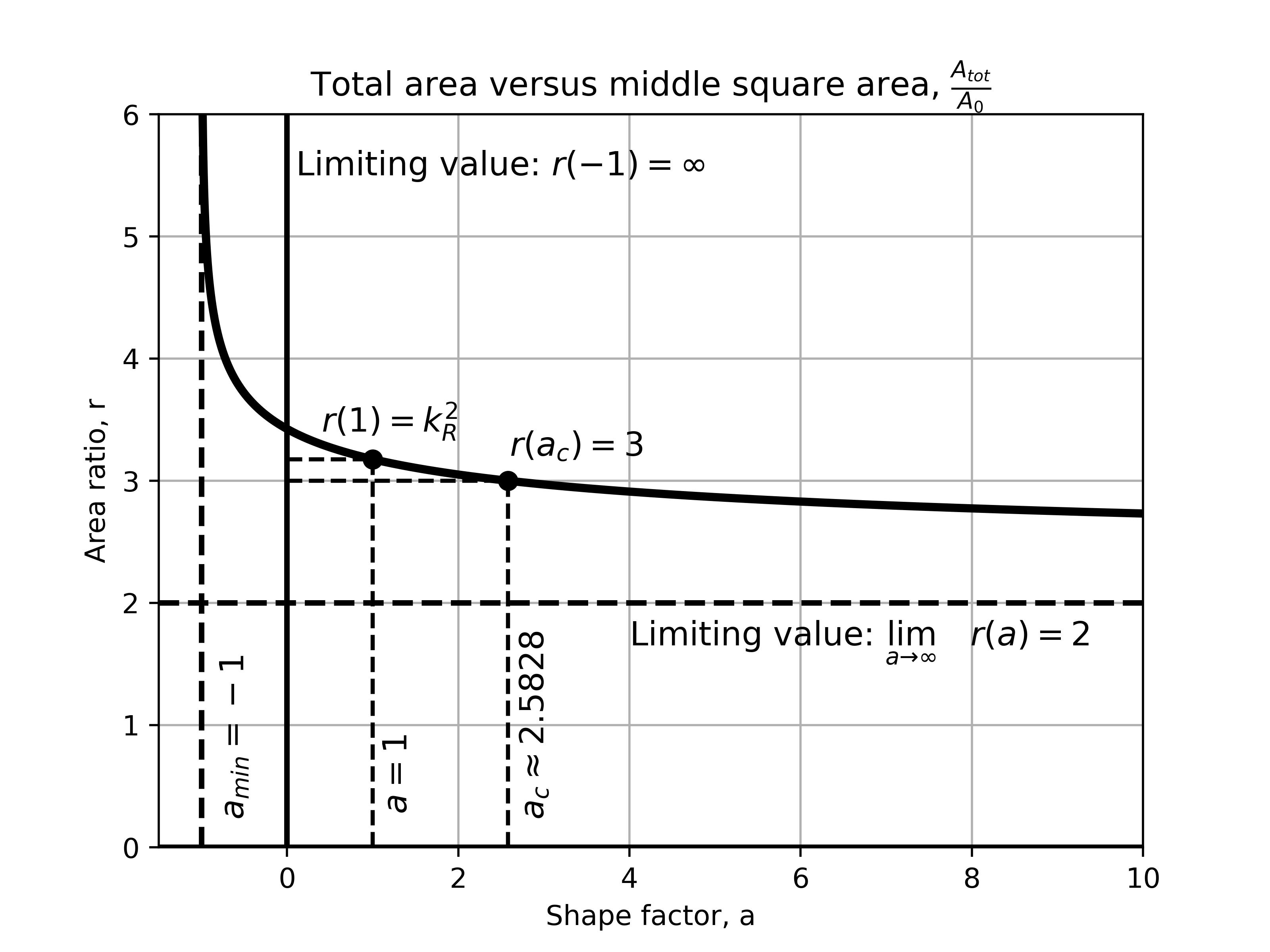}}
\end{center}
\caption{\label{Figure8} The area ratio $r$ versus the shape factor $a$.}
\end{figure}

This section discusses the effect of geometric properties of gemini functions on the representation of a dilogarithm. In other words,
we study how the shape of different area sections appear in the expressions of a dilogarithm function. All terms in a functional
identity always correspond to a certain plane area. Thus, gemini functions can be used to illustrate the formation of terms in these
particular identities. The valid domain for the shape factor of a gemini function is such that $a\in[-1,\infty)$. The graph in
Figure \ref{Figure8} illustrates the ratio of the total area $A_{tot}$ of a gemini function to the area of a middle square $A_{0}$.
The formula for this ratio is given by

\begin{figure}[htbp!]
\begin{center}
 \resizebox{10.0cm}{!}
 {\includegraphics{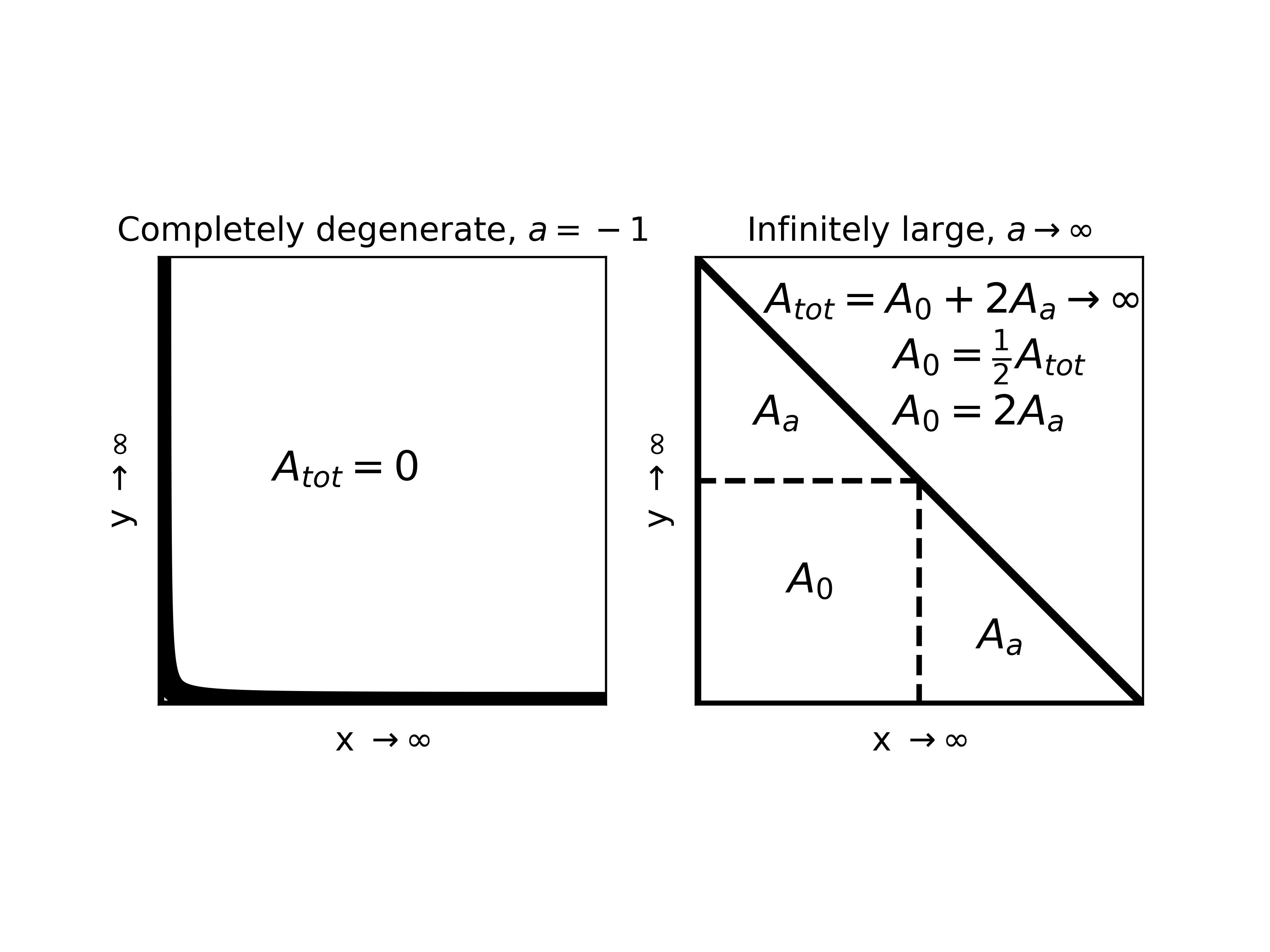}}
\end{center}
\caption{\label{Figure09} Schematic graphs of the $\gemini_{-1}(x)$- and $\gemini_{\infty}(x)$-functions.}
\end{figure}

\begin{eqnarray}
r(a)=\frac{\frac{\pi^{2}}{6}-\operatorname{Li_{2}}(-a)}{\ln^{2}(1+\sqrt{1+a})}.
\end{eqnarray}

The limiting values for this ratio $r(a)$ are such that $r(-1)=\infty$ and $\lim\limits_{a \to \infty}r(a)=2$. When the shape factor
$a$ (see Figure \ref{Figure09}) is equal to 1, i.e., $r(1)=\frac{\pi^{2}}{4\ln^{2}(1+\sqrt{2})}=k_{R}^{2}$. This is the area ratio
obtained from the fundamental form of the gemini function, i.e. , $\gemini_{1}(x)$. This value $k_{R}$ may be referred to as the
\emph{Grothendieck--Krivine constant} and is involved in the work of Pain \cite{Pai_23} on dilogarithm identities. Our construction
gives a geometric way of interpreting the Grothendieck--Krivine constant. 

\subsection{Representations related to a middle square}

In this section, we introduce identities related to the expression for a dilogarithm with a middle square area $A_{0}$ as a measure.

\begin{example}
Let us define the critical shape factor $a_{c}$ for a gemini function, whose middle square area $A_{o}$ is equal to the apex area $A_{a}$.
   We thus obtain the equality $\frac{1}{2}(A_{tot}-A_{0})= A_{a}$, and this implies that 
 $$ \frac{1}{2}\left[\frac{\pi^{2}}{6}-\operatorname{Li_{2}}(-a_{c})-\ln^{2}\left(
 1+\sqrt{1+a_{c}} \right) \right] = 
\ln^{2}\left(
 1+\sqrt{1+a_{c}} \right).$$
Consequently, we have that
$$\operatorname{Li_{2}}(-a_{c})=\frac{\pi^{2}}{6}-3\ln^{2}\left(
  1+\sqrt{1+a_{c}} \right)  \Rightarrow r(a_{c})=3. 
$$
The value $a_{c}$ is a critical value related to the middle square area $A_{0}$ versus the apex area $A_{a}$, because both areas
are equal in this case. When $a>a_{c} \Rightarrow A_{0}>A_{a}$ and vice versa.
\end{example}

\begin{example}
In Section 2.3, the reflection identity is derived with the aid of the degenerate gemini function. The respective areas for
$A_{0}$ and $A_{a}$ are also introduced in this context. Now, we derive on a general level, which gemini function satisfies a
condition such that $A_{a}=\operatorname{Li_{2}}(\frac{1}{x})$ and $A_{0}=\ln^{2}(x)$. Here, $\ln(x)$ denotes the fixed point of
the corresponding gemini function and the respective shape factor is such that $a=x^{2}-2x$. See for example Figure \ref{Figure5}
and \eqref{threetermfromfixed}. Hence, we can write

\vskip 0.1in

$A_{0}+2A_{a}=A_{tot} \Rightarrow \ln^{2}(x)+2\operatorname{Li_{2}}(\frac{1}{x})=\frac{\pi^{2}}{6}-\operatorname{Li_{2}}(2x-x^{2}) \Rightarrow$

\vskip 0.1in

$2\operatorname{Li_{2}}(\frac{1}{x})-\operatorname{Li_{2}}(1-2x+x^{2})-\ln(1-2x+x^{2})\ln(2x-x^{2})+\ln^{2}(x)=0 \Rightarrow$

\vskip 0.1in

$2\operatorname{Li_{2}}(\frac{1}{x})-\operatorname{Li_{2}}[(x-1)^{2}]-2\ln(x-1)\ln(2x-x^{2})+\ln^{2}(x)=0 \Rightarrow$

\vskip 0.1in

$2\operatorname{Li_{2}}\left(1-x\right)+\frac{\pi^{2}}{3}-\ln(x)\ln[\frac{x}{(x-1)^{2}}]-\operatorname{Li_{2}}[(x-1)^{2}]-$
\vskip 0.01in
$2\ln(x-1)\ln(2x-x^{2})+\ln^{2}(x)=0 \Rightarrow$

\vskip 0.1in

$2\operatorname{Li_{2}}(1-x)-\operatorname{Li_{2}}[(x-1)^{2}]+\frac{\pi^{2}}{3}-2\ln(x-1)\ln(2-x)=0 \Rightarrow$

\vskip 0.1in

$-\operatorname{Li_{2}}(x-1)+\frac{\pi^{2}}{6}-\ln(x-1)\ln(2-x)=0 \Rightarrow$

\vskip 0.1in

$\operatorname{Li_{2}}(2-x)+\ln(2-x)\ln(x-1)-\ln(x-1)\ln(2-x)=0 \Rightarrow$

\vskip 0.1in

$\operatorname{Li_{2}}(2-x)=0 \Rightarrow x=2 \Rightarrow x_{0}=\ln(2) \Rightarrow a=x^{2}-2x=2^{2}-2\cdot 2=0.$

\vskip 0.1in

We can determine the initial conditions for the degenerate form $\gemini_{0}(x)=\ln(\frac{1}{1-e^{-x}})$
of the gemini function, with $A_{0}=\ln^{2}(2)$ and $A_{a}=\operatorname{Li_{2}}(\frac{1}{2})$. If we assume that the
expression for the shape factor is greater than zero, then we can apply Landen's identity \eqref{Landensformula} to the term
$\operatorname{Li_{2}}(-(x^{2}-2x))$. In this case, we obtain that

\vskip 0.1in

$ \operatorname{Li_{2}}\left(\frac{1}{x-1}\right)=\frac{\pi^{2}}{6}-\ln(x-1)\ln\left(\frac{\sqrt{x-1}}{x-2}\right) \Rightarrow x=2. $

\end{example}

\begin{example}
Next, we investigate the second fixed point identity \eqref{threetermfromfixed}, where the area of the middle square is such that
${A_{0}=\frac{\pi^2}{6}}$. In this case, the terms $\frac{1}{2}A_{0}$ and $\frac{\pi^2}{12}$
vanish, as shown next. Hence, we can write

\vskip0.1in

$\frac{1}{2}A_{0}=\frac{1}{2}\ln^{2}(k)=\frac{\pi^2}{12} \Rightarrow k=e^{\pm \frac{\pi}{\sqrt{6}}}.$

\vskip0.1in

By inserting $k=e^{\frac{\pi}{\sqrt{6}}}$ into \eqref{threetermfromfixed}, we get the following identity without the constant
terms, as shown below $(k>1)$. We thus find that

\begin{eqnarray}\label{80}
\operatorname{Li_{2}}(2-k)-\operatorname{Li_{2}}\left(\frac{1}{k}\right)-\frac{1}{2}\operatorname{Li_{2}}(2k-k^{2})=0, \quad k=e^{\frac{\pi}{\sqrt{6}}}. 
\end{eqnarray}
\end{example}

It is mentioned earlier in Section 2.1 that the five-term gemini-identities \eqref{finalfivegemini} obtained from
$\gemini_{1}(a)$ and $\gemini_{a}(a)$ yield always to one and the same identity for $x_{1}=\ln(a)$. Next, we deal with this issue,
because the obtained result is linked to the previous examination of the second fixed point identity (19). Let us first manipulate
the identity obtained from the $\gemini_{1}(x)$-function at $x_{1}=\ln(a)$ and $x_{2}=\ln(\frac{a+1}{a-1})$. 
 This gives us that 
\begin{multline*}
 \operatorname{Li_{2}}\left(
 -\frac{1}{a} \right) - 
 \operatorname{Li_{2}}\left(
 \frac{1}{a} \right) + 
 \frac{\pi^{2}}{4}-\ln(a)\ln\left(
 \frac{a+1}{a-1} \right) = \\
 -\operatorname{Li_{2}}\left(
 -\frac{a-1}{a+1} \right) + 
 \operatorname{Li_{2}}\left(
 \frac{a-1}{a+1} \right).
\end{multline*} 
 By applying the reflection \eqref{reflectionformula} and Landen's \eqref{Landensformula} identities to the RHS terms, we get
 $$\operatorname{Li_{2}}\left(
 -\frac{1}{a} \right) - 
 \operatorname{Li_{2}}\left(
 \frac{1}{a} \right) + 
 \operatorname{Li_{2}}\left(
 \frac{2}{a+1} \right) + 
\operatorname{Li_{2}}\left(
 \frac{a+1}{2a} \right) = 
 \frac{\pi^{2}}{12}-\frac{1}{2}\ln^{2}\left(
 \frac{2a}{a+1} \right).$$ 
 By then manipulating the identity obtained from the $\gemini_{a}(x)$-function, where $x_{1}=\ln(a)$ and $x_{2}=\ln(\frac{2a}{a-1})$, 
  we find that 
\begin{multline*}
 \operatorname{Li_{2}}\left(
 -\frac{a}{a} \right) - 
 \operatorname{Li_{2}}\left(
 \frac{1}{a} \right) - 
 \operatorname{Li_{2}}(-a) + 
 \frac{\pi^{2}}{6}-\ln(a)\ln\left(
 \frac{2a}{a-1} \right) = \\ 
 -\operatorname{Li_{2}}\left(
 -a \cdot \frac{a-1}{2a} \right) + 
 \operatorname{Li_{2}}\left(
 \frac{a-1}{2a} \right).
\end{multline*} 
First, we apply the reflection identity \eqref{reflectionformula} to the third term of the LHS. Then we apply reflection
\eqref{reflectionformula} and Landen's \eqref{Landensformula} identities to the RHS terms. So, we can write

\vskip 0.1in

$\operatorname{Li_{2}}(-\frac{1}{a})-\operatorname{Li_{2}}(\frac{1}{a})+\operatorname{Li_{2}}(\frac{2}{a+1})+\operatorname{Li_{2}}(\frac{a+1}{2a})
=\frac{\pi^{2}}{12}-\frac{1}{2}\ln^{2}(a)-\ln(\frac{a+1}{2})\ln(\frac{\sqrt{2a+2}}{2a})$.

\vskip 0.1in

These two identities are the same, although the representations of the constant terms differ from each other. The constant terms of the first
identity enable the analytic calculation of the respective RHS root, which must also be the root for the whole identity. It simply means that
${\frac{\pi^{2}}{12}-\frac{1}{2}\ln(\frac{2a}{a+1})=0 \Rightarrow a=\frac{1}{2e^{\frac{\pi}{\sqrt{6}}}-1}\approx0.160988}$ or
${a=-\frac{e^{\frac{\pi}{\sqrt{6}}}}{2e^{\frac{\pi}{\sqrt{6}}}-2}\approx-2.245468}$. Next, we insert the positive root into the obtained identity, which
is also the common root for the both identities above. Hence, we get

\begin{eqnarray}\label{81}
\operatorname{Li_{2}}(k)+\operatorname{Li_{2}}\left(2-\frac{1}{k}\right)+\operatorname{Li_{2}}(1-2k)-\operatorname{Li_{2}}(2k-1)=0, \quad
k=e^{\pm\frac{\pi}{\sqrt{6}}}.
\end{eqnarray}

We can derive another four-term identity by combining these two constant term free identities \eqref{80} and \eqref{81}. The outcome
is shown below, and it is true, when $k>0$, although the identities \eqref{80} and \eqref{81} are only true at $k=e^{\frac{\pi}{\sqrt{6}}}$,
when $k\in\mathbb{R}$. (Identity \eqref{80} is also true at $k=e^{-\frac{\pi}{\sqrt{6}}}$).
\begin{multline*}
 \operatorname{Li_{2}}\left(\frac{1}{k}-1\right)-\operatorname{Li_{2}}\left(1-\frac{1}{k}\right)-\operatorname{Li_{2}}(1-2k) + \\ 
 \operatorname{Li_{2}}(2k-1)-\frac{\pi^{2}}{4}+\ln\left(\frac{1}{k}-1\right)\ln(2k-1)=0.
\end{multline*}

\subsection{On the median of a gemini function}
 A median can be defined for all gemini functions, as they are all monotonically decreasing functions and the area bounded by 
them with the positive coordinate axes is always finite, except when the shape factor $a$ tends to infinity. Whether one can ever calculate 
an analytic value for the median of the gemini function is another question, though. Anyway, it is a simple task to derive the general 
formula of a median for a $\gemini_{a}(x)$-function. Let $\ln(m)$ denote the median. Hence, the formula for a median is given by
 $$\int_{\ln(m)}^{\infty}\gemini_{a}(x) \, dx = 
 -\operatorname{Li_{2}}\left(-\frac{a}{m} \right) + 
 \operatorname{Li_{2}}\left(\frac{1}{m} \right) = 
 \frac{1}{2}\int_{0}^{\infty}\gemini_{a}(x) \, dx=\frac{\pi^{2}}{12}-\frac{1}{2}\operatorname{Li_{2}}(-a), $$
 and this gives us that 
\begin{eqnarray}\label{median}
\operatorname{Li_{2}}\left(\frac{1}{m}\right)-\operatorname{Li_{2}}\left(-\frac{a}{m}\right)=\frac{\pi^{2}}{12} - 
 \frac{1}{2}\operatorname{Li_{2}}(-a).
\end{eqnarray}
 This leads us toward the following properties concerning the medians for gemini functions.

 \ 

\noindent \textbf{Rule 1:} If the median $\ln(m)$ corresponds to the lower integration limit $x_{1}$ and the symmetric 
 upper integration limit
is such that $x_{2}=\ln(\frac{m+a}{m-1})$. Hence, the area $A_{c}$ between the integration limits is equal to the rectangle area
$A_{r}$, which is the product of $x_{1}$ and $x_{2}$, i.e., $A_{r}=\ln(m)\ln(\frac{m+a}{m-1})$. The formula for the first rule is given by 
\begin{multline*}
 A_{c}=\int_{\ln(m)}^{\ln(\frac{m+a}{m-1})}\gemini_{a}(x) \, dx = 
 \operatorname{Li_{2}}\left(-a\cdot\frac{m-1}{m + 
 a}\right)- \operatorname{Li_{2}}\left(\frac{m-1}{m+a}\right) - \\ 
 \operatorname{Li_{2}}\left(-\frac{a}{m}\right)+\operatorname{Li_{2}}\left(\frac{1}{m}\right)=
\ln(m)\ln\left(\frac{m+a}{m-1}\right).
\end{multline*} 
 This follows from the relations such that 
 $A_{tot}-(A_{c}+A_{r})=2A_{a}$ and such that 
 $\frac{1}{2}A_{tot}=A_{r}+A_{a}$, which imply that $A_{c}=A_{r}$.

 \ 

\noindent \textbf{Rule 2:} The area 
 $A_{\frac{1}{2}}$ between the median $\ln(m)$ and the fixed point $x_{0}=\ln(1+\sqrt{1+a})$ is always half of the area
of the middle square area $A_{0}$, i.e., $A_{\frac{1}{2}}=\frac{1}{2}\ln^{2}(1+\sqrt{1+a})$. This rule can be given by

\begin{eqnarray}
A_{\frac{1}{2}}=\int_{\ln(m)}^{\ln(1+\sqrt{1+a})}\gemini_{a}(x) \, dx=\operatorname{Li_{2}}\left(-\frac{a}{1+\sqrt{1+a}}\right)-
\nonumber
\end{eqnarray}
\begin{eqnarray}
\operatorname{Li_{2}}\left(\frac{1}{1+\sqrt{1+a}}\right)-
\operatorname{Li_{2}}\left(-\frac{a}{m}\right)+\operatorname{Li_{2}}\left(\frac{1}{m}\right)=\frac{1}{2}\ln^{2}(1+\sqrt{1+a}). 
\end{eqnarray}

This above formula can be simply obtained by combining the fixed-point identity in
 Theorem \ref{firstfixedpoint} and the median formula in \eqref{median}.
 
Next, we generate some two-term identities based on properties of a median of gemini functions without knowing exact values of
dilogarithm arguments. The primary aim is just to visualize representations of a dilogarithm function with different kind of arguments.
Let us study a following function, $\gemini_{m^{n}}(x)=\ln(\frac{1+m^{n}e^{-x}}{1-e{-x}})$. Now, the shape factor $a$ is
naturally $m^{n}$ for $n\in(-\infty,\sqrt{2}+2]$. Let the median be such that $x_{1}=\ln(m)$ for $m\in[m_{mim},\infty)$.
The maximum value for $n$ and the minimum value for $m$ will be clarified further ahead. The median equation is given by
 $$-\operatorname{Li_{2}}\left(
 -\frac{m^{n}}{m} \right) + 
 \operatorname{Li_{2}}\left(\frac{1}{m}\right)=\frac{\pi^{2}}{12}
-\frac{1}{2}\operatorname{Li_{2}}\left(
 -m^{n} \right), $$
 and this implies that 
 $$-\operatorname{Li_{2}}(-m^{n-1})+\operatorname{Li_{2}}\left(\frac{1}{m}\right)+
\frac{1}{2}\operatorname{Li_{2}}(-m^{n})-\frac{\pi^{2}}{12}=0. $$

\begin{example}
 Let us first define the liniting case for $m=m_{min}$, when $n$ tends to minus infinity. Thus, we get
$$\lim\limits_{n\to -\infty} -\operatorname{Li_{2}}\left(
 -m_{min}^{n-1} \right) + 
 \operatorname{Li_{2}}\left(
 \frac{1}{m_{min}} \right) + 
 \frac{1}{2}\operatorname{Li_{2}}\left(
 -m_{min}^{n} \right) - \frac{\pi^{2}}{12}=0, $$
 and this implies that $$ \operatorname{Li_{2}}\left(\frac{1}{m_{min}}\right)=\frac{\pi^{2}}{12}. $$
\end{example}

The obtained value of $\ln(m_{min})$ corresponds to the median of the degenerate form of a gemini function, i.e., $\gemini_{0}(x)$.

\begin{example}
Let us also evaluate the special case, when $n=0$. Hence, we can write
 $$-\operatorname{Li_{2}}\left(-m^{0-1} \right) + 
 \operatorname{Li_{2}}\left(\frac{1}{m} \right) + 
 \frac{1}{2}\operatorname{Li_{2}}\left(-m^{0} \right) - \frac{\pi^{2}}{12} = 0, $$
 and this implies that $$
\operatorname{Li_{2}}\left(\frac{1}{m}\right)-\operatorname{Li_{2}}\left(-\frac{1}{m}\right)=\frac{\pi^{2}}{8} \Rightarrow
\operatorname{\chi_{2}}\left(\frac{1}{m}\right)=\frac{\pi^{2}}{16}. $$
 The obtained value of $\ln(m)$ is 
 the median of the fundamental form of a gemini function, i.e., $\gemini_{1}(x)$.
\end{example}

\begin{example}
If $n=1$ then $x_{1}=\ln(m)=\ln(a)$ and the identity becomes extremely simple. We can write
 $$-\operatorname{Li_{2}}(-1)+\operatorname{Li_{2}}\left(\frac{1}{a}\right)+\frac{1}{2}\operatorname{Li_{2}}(-a)-\frac{\pi^{2}}{12}=0, $$
 and this implies that 
 $$\operatorname{Li_{2}}\left(\frac{1}{a}\right)+\frac{1}{2}\operatorname{Li_{2}}(-a)=0. $$
\end{example}

\begin{example}
 If $n=2$ then the median equation is given by $$-\operatorname{Li_{2}}\left( -\frac{m^{2}}{m} \right) + \operatorname{Li_{2}}\left( 
 \frac{1}{m} \right) - \frac{\pi^{2}}{12}+
\frac{1}{2}\operatorname{Li_{2}}(-m^{2}), $$
 and this implies that 
\begin{multline*}
 \operatorname{Li_{2}}\left(
 -\frac{1}{m} \right) + 
 \frac{\pi^{2}}{6}+\frac{1}{2}\ln^{2}(m)+\operatorname{Li_{2}}\left(
 \frac{1}{m} \right) - \\ 
 \frac{1}{2}\operatorname{Li_{2}}\left(
 -\frac{1}{m^{2}} \right) - 
 \frac{\pi^{2}}{12}-\frac{1}{4}\ln^{2}(m^{2})-\frac{\pi^{2}}{12} = 0, 
\end{multline*}
 so that 
 $$\frac{1}{2}\operatorname{Li_{2}}\left(
 \frac{1}{m^{2}} \right) - 
 \frac{1}{2}\operatorname{Li_{2}}\left(
 -\frac{1}{m^{2}} \right) + 
 \frac{1}{2}\ln^{2}(m)-\ln^{2}(m)=0. $$ 
 As a consequence, we find that $$ \operatorname{Li_{2}}\left(\frac{1}{m^{2}}\right) - 
 \operatorname{Li_{2}}\left(-\frac{1}{m^{2}}\right)=\ln^{2}(m)
\Rightarrow \operatorname{\chi_{2}}\left(\frac{1}{m^{2}}\right)=\frac{1}{2}\ln^{2}(m). 
$$
\end{example}

\begin{example}
 We apply the same $\gemini_{m^{2}}(x)$-function, where the lower integration
limit and the shape factor are related in such a way that $x_{1}=\ln(a)=\ln(m^{2})$. Here, the lower integration limit $x_{1}$ is not
 a median.
 Now, the corresponding upper integration limit is such that $x_{2}=\ln(\frac{2m^{2}}{m^{2}-1})$. Hence, this five-term identity is given by
\begin{multline*}
 \operatorname{Li_{2}}\left(
 -\frac{m^{2}}{m^{2}} \right) - 
 \operatorname{Li_{2}}\left(
 \frac{1}{m^{2}} \right) 
 - \operatorname{Li_{2}}(-m^{2}) + 
 \frac{\pi^{2}}{6} - \\
 \ln(m^{2})\ln\left(
 \frac{2m^{2}}{m^{2}-1} \right) + 
 \operatorname{Li_{2}}\left( 
 -m^{2}\frac{m^{2}-1}{2m^{2}}\right) - 
 \operatorname{Li_{2}}\left(
 \frac{m^{2}-1}{2m^{2}} \right) = 0. 
\end{multline*}
 Consequently, we have that

\vskip 0.1in

$\operatorname{Li_{2}}(-1)-\operatorname{Li_{2}}\left(\frac{1}{m^{2}}\right) + 
\operatorname{Li_{2}}\left(-\frac{1}{m^{2}}\right)+\frac{\pi^{2}}{3}+\frac{1}{2}\ln^{2}(m^{2})-$\\
\vskip 0.01in
$2\ln(m)\ln\left(\frac{2m^{2}}{m^{2}-1} \right) + \operatorname{Li_{2}}\left(-\frac{m^{2}-1}{2}\right)- 
\operatorname{Li_{2}}\left(\frac{m^{2}-1}{2m^{2}}\right)=0.$


\vskip 0.1in

$\frac{\pi^{2}}{4}+\ln^{2}(m)-2\ln(m)\ln\left(\frac{2m^{2}}{m^{2}-1}\right)+\operatorname{Li_{2}}\left(-\frac{m^{2}-1}{2}\right)-
\operatorname{Li_{2}}\left(\frac{m^{2}-1}{2m^{2}}\right)=0 \Rightarrow$

\begin{eqnarray}
\operatorname{Li_{2}}\left(\frac{2}{m^{2}+1}\right)-\operatorname{Li_{2}}\left(\frac{m^{2}-1}{2m^{2}}\right)+\frac{\pi^{2}}{12}+\ln^{2}(m)-
\nonumber
\end{eqnarray}
\begin{eqnarray}
2\ln(m)\ln\left(\frac{2m^{2}}{m^{2}-1}\right)+\frac{1}{2}\ln\left(\frac{m^{2}+1}{2}\right)\ln\left(\frac{2m^{2}+2}{m^{4}-2m^{2}+1}\right)=0
\end{eqnarray}

\end{example}

\begin{example}
If $n=3$, we get $-\operatorname{Li_{2}}(-m^{2})+\operatorname{Li_{2}}(\frac{1}{m})+\frac{1}{2}\operatorname{Li_{2}}(-m^{3})-\frac{\pi^{2}}{12}=0$. 
We then obtain the ladder relation for m  such that 
\begin{multline*}
4\operatorname{Li_{2}}\left(\frac{1}{m}\right)-4\operatorname{Li_{2}}\left(\frac{1}{m^{2}}\right)+
2\operatorname{Li_{2}}\left(\frac{1}{m^{3}}\right)+ \\ 
2\operatorname{Li_{2}}\left(\frac{1}{m^{4}}\right)-\operatorname{Li_{2}}\left(\frac{1}{m^{6}}\right)-\ln^{2}(m)=0. 
\end{multline*}
\end{example}

\begin{example}
 For $n=4$, the median equation is given by
$$-\operatorname{Li_{2}}(-m^{3})+\operatorname{Li_{2}}\left(\frac{1}{m}\right)+\frac{1}{2}\operatorname{Li_{2}}(-m^{4})-\frac{\pi^{2}}{12}=0.$$
 This equation has no real roots, i.e., it does not intersect the horizontal axis at all, when $n\ge4$ and $m>m_{min}$. This implies that
there must exist a critical or a limiting value for $n$, which defines whether the equation has a root or not at infinity, and this
value lies between 3 and 4. To find this critical value for the parameter $n$, we have to define the expression below in such a way that
 $$\lim\limits_{m\to \infty}-\operatorname{Li_{2}}(-m^{n-1}) + 
 \operatorname{Li_{2}}\left( 
 \frac{1}{m} \right) + 
 \frac{1}{2}\operatorname{Li_{2}}(-m^{n})-
\frac{\pi^{2}}{12}=0.$$
 Next, we have to convert the first and the third term equipped with the negative arguments by applying an inversion formula. We thus
obtain that
\begin{multline*}
 \lim\limits_{m\to \infty}\operatorname{Li_{2}}\left( 
 -\frac{1}{m^{n-1}} \right) + 
 \operatorname{Li_{2}}\left(
 \frac{1}{m} \right) + 
 \frac{1}{2}\operatorname{Li_{2}}\left(
 -\frac{1}{m^{n}} \right) + \\ 
 \left[\frac{1}{2}(n-1)^{2}-\frac{1}{4}n^{2}\right]\ln^{2}(m)=0. 
\end{multline*} 
 Hence, the equation for the parameter $n$ is given by
 $$\left[\frac{1}{2}(n-1)^{2}-\frac{1}{4}n^{2}\right]=0 \Rightarrow n=2\pm\sqrt{2}$$ for $n\in(3,4)$,
 and thus $n=2+\sqrt{2}$. 
 Let us substitute the parameter $n=2+\sqrt{2}$ back into the original median equation, and let $m$ tend to infinity. 
 We thus obtain the vanishing relation 
$$\lim\limits_{m\to \infty}-\operatorname{Li_{2}}(-m^{\sqrt{2}+1})+\operatorname{Li_{2}}\left(\frac{1}{m}\right)+
\frac{1}{2}\operatorname{Li_{2}}(-m^{\sqrt{2}+2})-\frac{\pi^{2}}{12}=0.$$
The middle term can be discarded, and hence the relation 
 $$ \lim_{m \to \infty} \frac{1}{2}\operatorname{Li_{2}}\left(-m^{\sqrt{2}+2}\right)-\operatorname{Li_{2}}\left(-m^{\sqrt{2} + 
 1}\right)=\frac{\pi^{2}}{12}. $$
\end{example}

\begin{figure}[htbp!]
\begin{center}
 \resizebox{10.0cm}{!}
 {\includegraphics{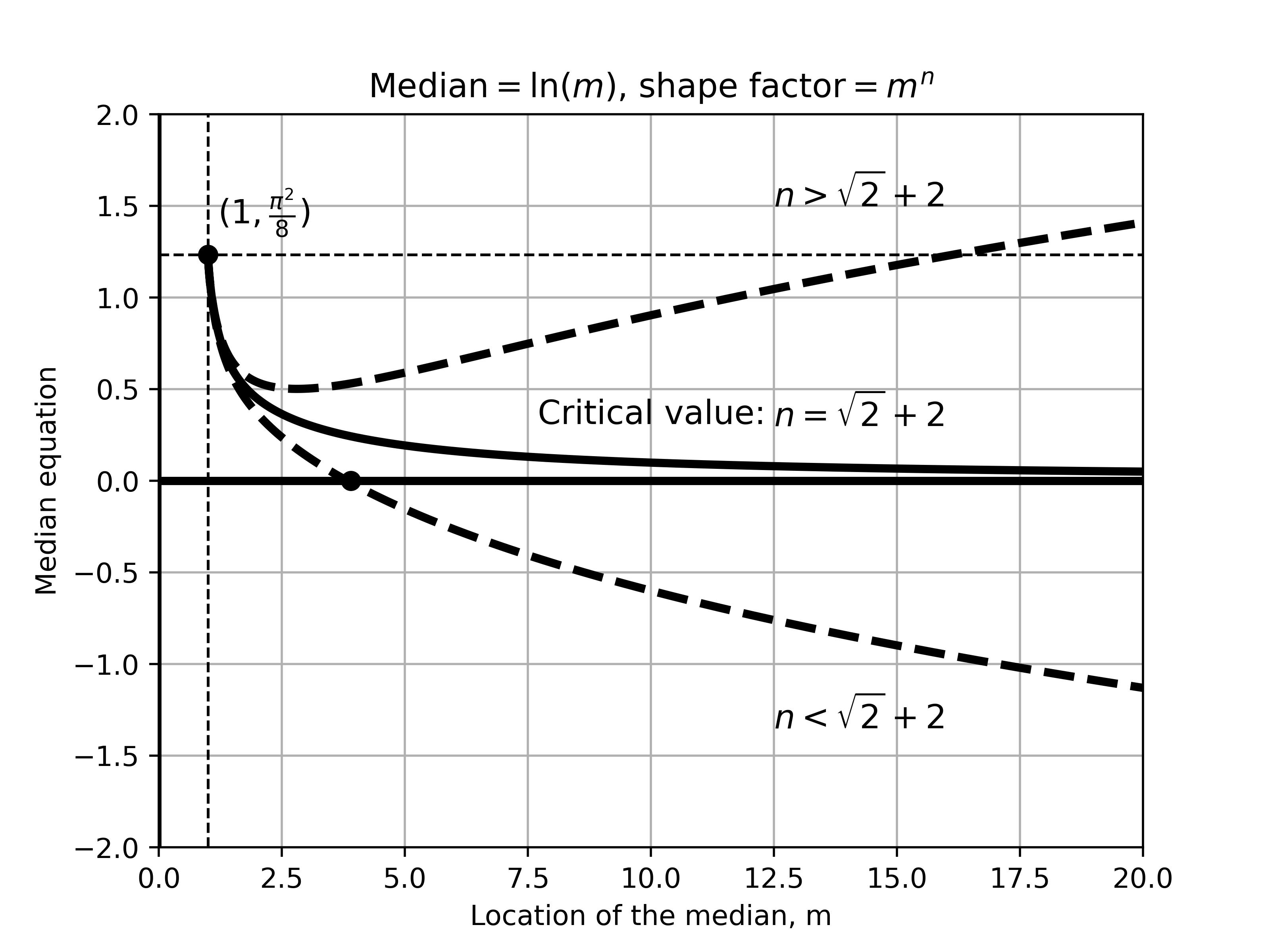}}
\end{center}
\caption{\label{Figure10} The solid asymptotic graph depicts the median at infinity.}
\end{figure}

The three graphs in Figure \ref{Figure10} illustrate the situation related to the behavior of the median. We defined the relation
between the shape factor and the median in such a way that the median is at $\ln(m)$ and the shape factor $a=m^{n}$ for
$n\in(-\infty,\sqrt{2}+2]$ and $m\in[m_{min},\infty)$. What does this mean in practice? We have a limiting gemini function
$\gemini_{\infty}(x)=\lim\limits_{m\to \infty}\ln(\frac{1+m^{\sqrt{2}+2}e^{-x}}{1-e^{-x}})$, whose median $\ln(m)$ is located at infinity.
The corresponding median function approaches asymptotically to the horizontal axis at infinity, when $n$ is critical, i.e., $n=\sqrt{2}+2$.
Immediately, after an infinitesimal increase of $n$, the function no longer touches the horizontal axis, and the median for the
limiting gemini function can not be determined. The formula above representing the asymptotic median function can also be interpreted
as a two-term single value identity for the infinity. It is worth to emphasize that this is not a unique case. We can derive at least one
corresponding asymptotic median equation by applying the function $\gemini_{m^{2n}}(x)$ in such a way that the median $x_{1}=\ln(m^{2n+1})$. The obtained formula is given by
 $$\lim\limits_{m\to \infty}-\operatorname{Li_{2}}(-m)+\operatorname{Li_{2}}\left(\frac{1}{m^{\sqrt{2}-1}}\right)+
\frac{1}{2}\operatorname{Li_{2}}(-m^{\sqrt{2}})-\frac{\pi^{2}}{12}=0, $$
 and this gives us that 
 $$ \lim\limits_{m\to \infty} \frac{1}{2}\operatorname{Li_{2}}\left(-m^{\sqrt{2}}\right)-\operatorname{Li_{2}}(-m)=\frac{\pi^{2}}{12}.$$

 \begin{figure}[htbp!]
\begin{center}
\resizebox{10.cm}{!}
{\includegraphics{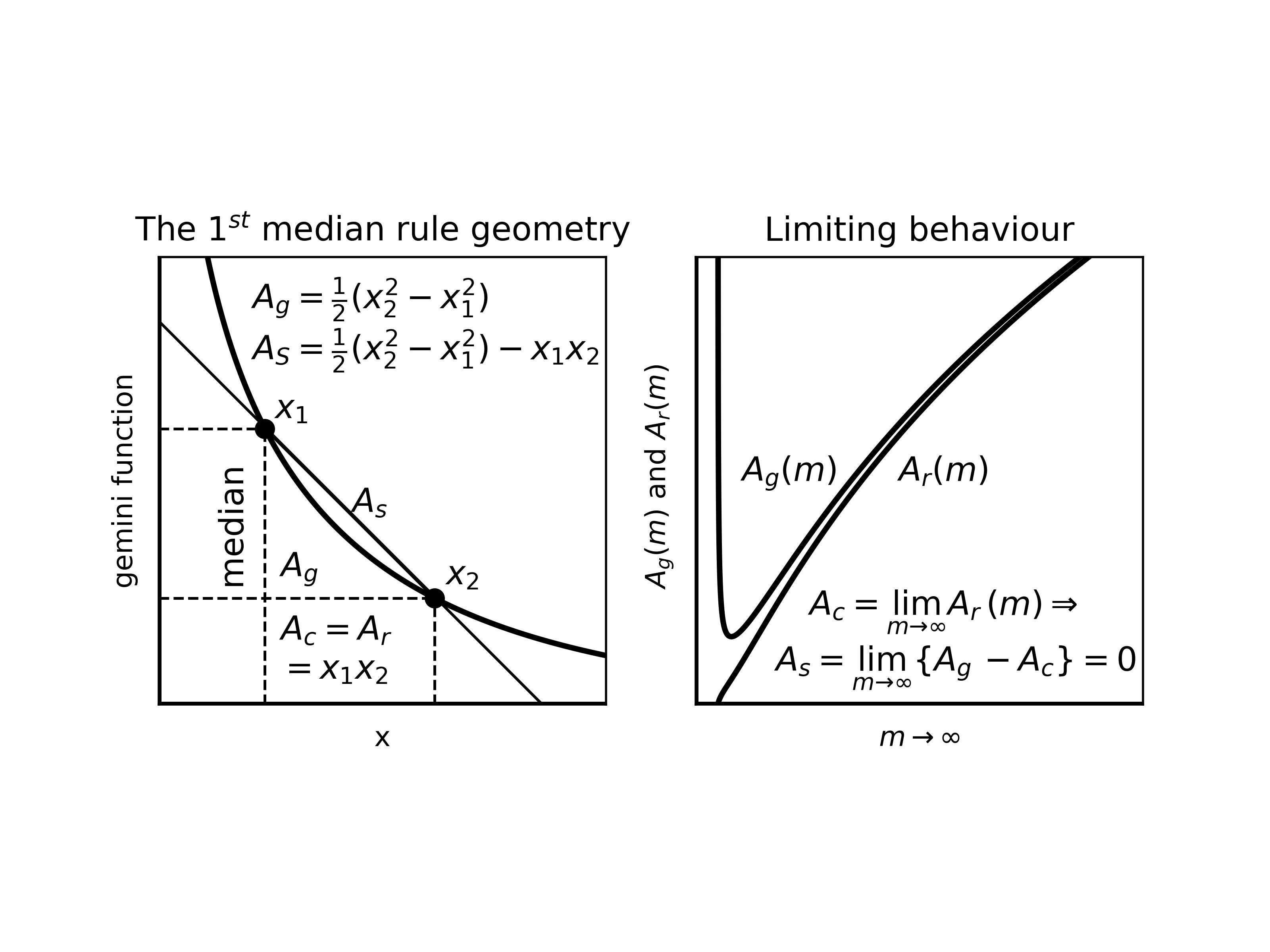}}
\end{center}
\caption{\label{Figure11} The left-hand
 plot illustrates the geometry related to the first median rule. The right-hand
 plot depicts how
the functions $A_{g}(m)$ and $A_{r}(m)$ asymptotically approach one another. }
\end{figure}
 
 The gemini function $$\gemini_{\infty}(x)=\lim\limits_{m\to \infty}\ln\left(\frac{1+m^{\sqrt{2}+2}e^{-x}}{1-e^{-x}}\right),$$ gives us a way, 
 informally, of showing how the graphs of gemini functions can straighten depending on the behavior of the associated parameters. It 
 is worth to point out that this proof is not based on the scale factor $b$. Here, the shape factor $a$ tends to infinity. Let us define such 
 that the shape factor $a=\lim\limits_{m\to \infty} m^{2+\sqrt{2}}$ and the median $x_{1}=\lim\limits_{m\to \infty} \ln(m)$ as earlier. 
 Hence, the corresponding upper integration limit or the symmetric point locates also at infinity in such a way that $x_{2}=\lim\limits_{m 
 \to \infty} \ln(\frac{m+m^{2+\sqrt{2}}}{m-1})$. This proof uses the first median rule, which states that $$A_{c}=\lim\limits_{m\to \infty} 
 \int_{\ln(m)}^{\ln(\frac{m+m^{2+\sqrt{2}}}{m-1})}\gemini_{m^{2+\sqrt{2}}}(x) \, dx = A_{r}=\lim\limits_{m\to \infty} \ln(m) 
 \ln\left(\frac{m+m^{2+\sqrt{2}}}{m-1}\right).$$ 

We proceed to calculate the segment area $A_{s}$ between $x_{1}$ and $x_{2}$, which is equal to $A_{g}-A_{r}$, and this is illustrated in 
 Figure \ref{Figure11}. In this case, the area $A_{g}$ may be expressed as $$ \lim\limits_{m\to \infty} 
 \left\{\ln(m)\left[\ln\left(\frac{m+m^{2+\sqrt{2}}}{m-1}\right)-\ln(m)\right] + \frac{1}{2}\left[\ln\left(\frac{m+m^{2+\sqrt{2}}}{m - 
 1}\right)-\ln(m)\right]^{2}\right\}, $$ and we rewrite this as $$ \lim\limits_{m\to \infty} \left\{\ln(m)\ln\left(\frac{1 + m^{1+\sqrt{2}}}{m - 
 1}\right)+\frac{1}{2}\ln^{2}\left(\frac{1+m^{1+\sqrt{2}}}{m-1}\right)\right\}. $$ The segment area may be expressed as $A_{s} = A_{g} 
 - A_{c}=A_{g}-A_{r}$, and we rewrite this as $$\lim\limits_{m\to \infty} \left\{\ln(m)\ln\left(\frac{1 + m^{1+\sqrt{2}}}{m-1}\right) + 
 \frac{1}{2}\ln^{2}\left(\frac{1+m^{1+\sqrt{2}}}{m-1}\right) - \ln(m)\ln\left(\frac{m+m^{2+\sqrt{2}}}{m-1}\right)\right\}, $$ and this implies 
 that $$A_{s}=\lim\limits_{m\to \infty} \left\{\frac{1}{2}\ln^{2}\left(\frac{1+m^{1+\sqrt{2}}}{m-1}\right) + \ln(m)\ln\left(\frac{1+m^{1 + 
 \sqrt{2}}}{m+m^{2+\sqrt{2}}}\right)\right\}, $$ so that $$\lim\limits_{m\to \infty} \left\{\frac{1}{2}\ln^{2}\left(\frac{1+m^{1 + 
 \sqrt{2}}}{m-1}\right)-\ln^2(m)\right\}=0.$$ Since the segment area $A_{s}$ vanishes between the median and its corresponding 
 symmetric point, the graph of an infinite gemini function is a straight line, i.e., linear in this domain. Thus, the graphs of infinite gemini 
 functions can be formulated in such a way that $\gemini_{\infty}(x)=\lim\limits_{C\to \infty}-x+C$ for $x \in (x_{1}, x_{2})$. Unfortunately, 
 this proof can not reveal the behavior of infinite gemini functions outside this domain, although the interval between $x_{1}$ and 
 $x_{2}$ is infinite.

\subsection{Representations related to areas determined by integration limits}
 Next, we deal with another area ratio $r$ related to the $\gemini_{a}(x)$-function. In this case, we define the ratio of the total area to the 
 area between the integration limits, i.e., $r=\frac{A_{tot}}{A_{c}}$. We can formulate the following equation for this purpose in such a way 
 that $2A_{a} = A_{tot} - A_{r} - A_{c}$, where $A_{a}$ denotes the apex area and $A_{r}$ denotes the product of the integration limits, i.e., 
 $A_{r}=x_{1}x_{2} = \ln(x) \ln(\frac{x+a}{x-1})$. Here, the valid domain for the lower integration limit is such that 
$x \in [0,1+\sqrt{1+a}]$. By combining all the terms together, we obtain the vanishing relation 
$$ \operatorname{Li_{2}}\left(-\frac{a}{x}\right)-\operatorname{Li_{2}}\left(\frac{1}{x}\right)+
 \left[\frac{r+1}{2r}\right]\left[\frac{\pi^{2}}{6}-\operatorname{Li_{2}}(-a)\right]
 -\frac{1}{2}\ln(x)\ln\left(\frac{x+a}{x-1}\right)=0, $$ which implies that 
$$ r\left(x,a\right)=\frac{\operatorname{Li_{2}}(-a)-\frac{\pi^{2}}{6}}{\frac{\pi^{2}}{6}-\ln(x)\ln\left(\frac{x+a}{x-1}\right)-
 \operatorname{Li_{2}}(-a)-2\operatorname{Li_{2}}(\frac{1}{x})+2\operatorname{Li_{2}}\left(-\frac{a}{x}\right)}. $$
 The limiting values for the ratio $r(x,a)$ are given by
 $$\lim\limits_{x\to 0}r(x,a)=1 \ \ \ \text{and} \ \ \ \lim\limits_{x\to 1+\sqrt{1+a}}r(x,a)=\infty.$$ 
 There exist only four gemini functions for which an exact value can be assigned to the ratio $r$. Two integer values can be defined for the
ratio $r$. These two functions are $\gemini_{0}(x)$ at $x=\phi$ and $\gemini_{1}(x)$ also at $x=\phi$. For the degenerate form we can write
 $$r(\phi,0)=\frac{\operatorname{Li_{2}}(0)-\frac{\pi^{2}}{6}}{\frac{\pi^{2}}{6}-\ln(\phi)\ln\left(\frac{\phi+0}{\phi-1}\right)-
\operatorname{Li_{2}}(0)-2\operatorname{Li_{2}}(\frac{1}{\phi})+2\operatorname{Li_{2}}\left(-\frac{0}{\phi}\right)}=5.$$ 

 Similarly, we have that 
 $$r(\phi,1)=\frac{\operatorname{Li_{2}}(-1)-\frac{\pi^{2}}{6}}{\frac{\pi^{2}}{6}-\ln(\phi)\ln\left(\frac{\phi+1}{\phi-1}\right)-
\operatorname{Li_{2}}(-1)-2\operatorname{Li_{2}}(\frac{1}{\phi})+2\operatorname{Li_{2}}\left(-\frac{1}{\phi}\right)}=3.$$ 
 It is possible to evaluate two 
 exact values for the ratio $r$ related to the same $\gemini_{+\phi}(x)$-function, whose lower integration
 limits are such that $x=\ln(\phi)$ or $x=\ln(\sqrt{\phi})$. The ratio for $x=\phi$ is given by
 $$r(\phi,\phi)=\frac{8\pi^{2}+30\ln^{2}(\phi)}{3\pi^{2}-60\ln^{2}(2)-30\ln(\frac{1}{4}\phi)\ln(2\phi)}.$$ 
 The ratio for $x_{1}=\sqrt{\phi}$ is given by
 $$r(\sqrt{\phi},\phi)=\frac{16\pi^{2}+60\ln^{2}(\phi)}{10\pi^{2}-105\ln^{2}(\phi)+30\ln(\phi)\ln(\frac{\sqrt{\phi}+\phi}{\sqrt{\phi} - 
 1})}.$$ 

\begin{example}
Next, we study gemini function pairs in such a way that the shape factors are reciprocal with respect to each other, e.g.
these two functions $\gemini_{a}(x)$ and $\gemini_{\frac{1}{a}}(x)$ form an inversion function pair. This simply means that $a_{1}=a$ and
$a_{2}=\frac{1}{a}$ for $a>1$. In addition to this, we define two rectangles, whose widths are equal in such a way that the abscissa value of
the lower right corner is the same for both rectangles, i.e., the lower integration limit is common for both functions. The heights of these two
rectangles are defined such that $y_{1}=\gemini_{a}(x_{1})$ and $y_{2}=\gemini_{\frac{1}{a}}(x_{1})$. The area ratio $n$ is defined in such a way that
$\frac{A_{r1}}{A_{r2}}=\frac{A_{tot1}}{A_{tot2}}$. See Figure \ref{Figure12}, which clarifies this configuration. By using the inversion identity,
we can write

\vskip 0.1in

$n=\frac{A_{tot1}}{A_{r1}}=\frac{A_{tot2}}{A_{r2}}=\frac{\frac{\pi^{2}}{6}-\operatorname{Li_{2}}(-a)}{\ln(x)\ln(\frac{x+a}{x-1})}=
\frac{\frac{\pi^{2}}{6}-\operatorname{Li_{2}}(-\frac{1}{a})}{\ln(x)\ln(\frac{x+\frac{1}{a}}{x-1})} \Rightarrow \frac{\frac{\pi^{2}}{6}
-\operatorname{Li_{2}}(-a)}{\ln(\frac{x+a}{x-1})}=\frac{\frac{\pi^{2}}{6}-\operatorname{Li_{2}}(-\frac{1}{a})}{\ln(\frac{x+\frac{1}{a}}{x-1})}
\Rightarrow $

\vskip 0.1in

$\ln[\frac{(a+x)(ax+1)}{a(x-1)^{2}}]\operatorname{Li_{2}}(-a)+\ln[\frac{a(a+x)^{2}}{(ax+1)(x-1)}]\frac{\pi^{2}}{6}+
\frac{1}{2}\ln^{2}(a)\ln(\frac{x+a}{x-1})=0
\Rightarrow$

\begin{eqnarray}\label{64}
\operatorname{Li_{2}}(-a)=-\biggl\{\frac{\ln[\frac{a(a+x)^{2}}{(ax+1)(x-1)}]}{\ln[\frac{(a+x)(ax+1)}{a(x-1)^{2}}]}\biggl\}\frac{\pi^{2}}{6}
-\frac{1}{2}\biggl\{\frac{\ln(\frac{x+a}{x-1})}{\ln[\frac{(a+x)(ax+1)}{a(x-1)^{2}}]}\biggl\}\ln^{2}(a).
\end{eqnarray}

\begin{figure}[htbp!]
\begin{center}
 \resizebox{8.0cm}{!}
 {\includegraphics{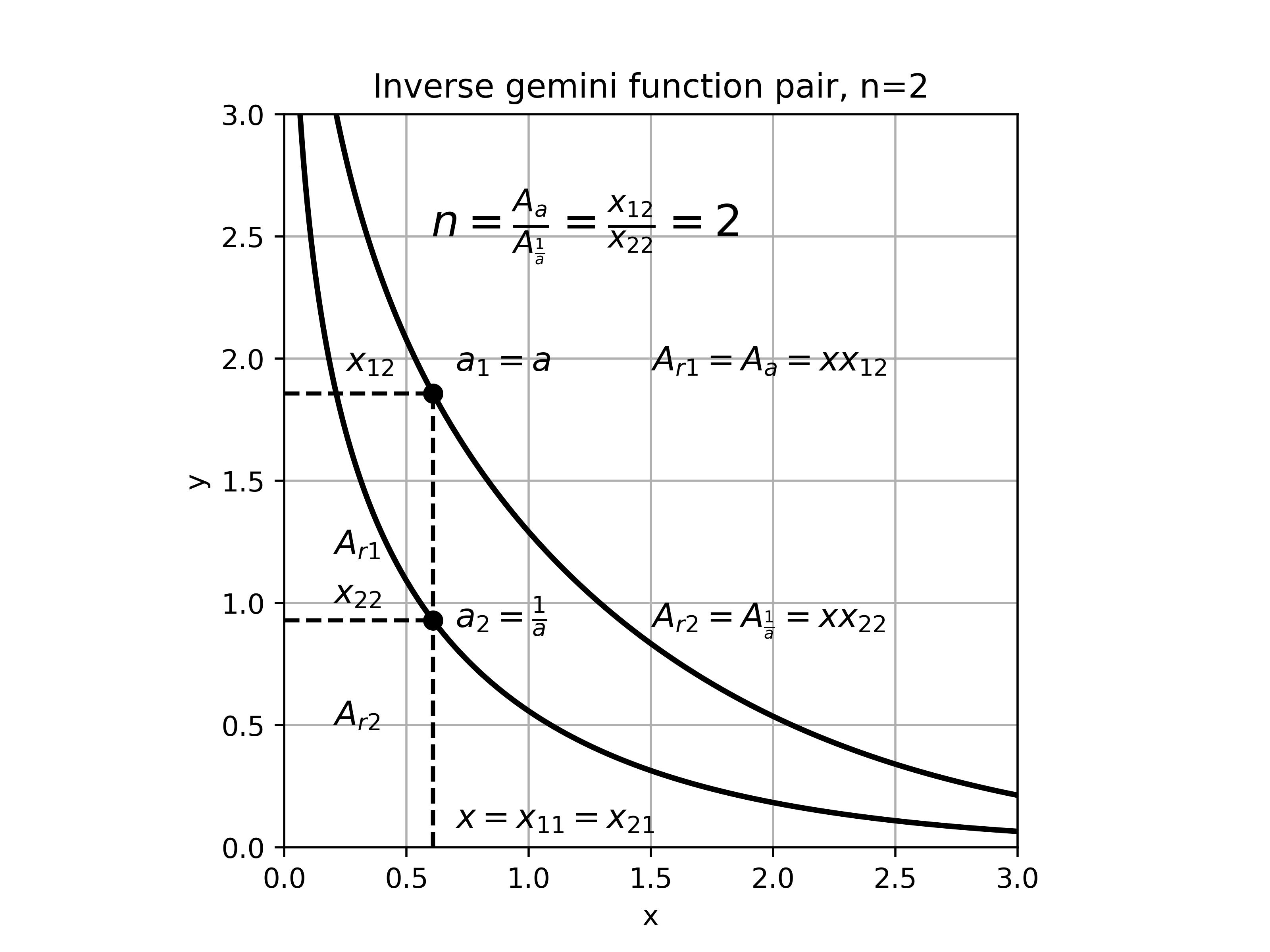}}
\end{center}
\caption{\label{Figure12} Illustration of the inverse gemini function pair, where the ratio $n=2$.}
\end{figure}
 
 Next, we simplify the coefficients of \eqref{64}. We denote these coefficients so that 
 $$C_{1}=\biggl\{\frac{\ln[\frac{a(a+x)^{2}}{(ax+1)(x-1)}]}{\ln[\frac{(a+x)(ax+1)}{a(x-1)^{2}}]}\biggl\}$$ 
 and so that 
$$C_{2}=\biggl\{\frac{\ln(\frac{x+a}{x-1})}{\ln[\frac{(a+x)(ax+1)}{a(x-1)^{2}}]}\biggl\}, $$ 
 and we simplify $C_{1}$ to obtain that 
\begin{multline*}
 C_{1}=\biggl\{\frac{\ln[\frac{a(x+a)(x+a)}{(ax+1)(x-1)}]}{\ln[\frac{(a+x)(ax+1)}{a(x-1)(x-1)}]}\biggl\}=
\biggl\{\frac{\ln(\frac{x+a}{x-1})+\ln(\frac{x+a}{x+\frac{1}{a}})}{\ln(\frac{x+a}{x-1})+\ln(\frac{x+\frac{1}{a}}{x-1})}\biggl\} = \\ 
\biggl\{\frac{\ln(\frac{x+a}{x-1})+\ln\left[\frac{(\frac{x+a}{x-1})}{(\frac{a+\frac{1}{a}}{x-1})}\right]}{\ln(\frac{x+a}{x-1})+
\ln(\frac{x+\frac{1}{a}}{x-1})}\biggl\} = 
\biggl\{\frac{2\ln(\frac{x+a}{x-1})-\ln(\frac{x+\frac{1}{a}}{x-1})}{\ln(\frac{x+a}{x-1})+
\ln(\frac{x+\frac{1}{a}}{x-1})}\biggl\}.
\end{multline*}
 According to Figure \ref{Figure12}, we can write $x_{12}=\ln(\frac{x+a}{x-1})$ and $x_{22}=\ln(\frac{x+\frac{1}{a}}{x-1})$. The area ratio
$n=\frac{A_{r1}}{A_{r2}}=\frac{x \cdot x_{21}}{x \cdot x_{22}}=\frac{x_{21}}{x_{22}} \Rightarrow x_{21}=nx_{22}$. Hence, we can write
 $$C_{1}=\biggl\{\frac{2x_{12}-x_{22}}{x_{12}+x_{22}}\biggl\} = 
 \biggl\{\frac{2nx_{22}-x_{22}}{nx_{22}+x_{22}}\biggl\}=\biggl\{\frac{2n-1}{n+1}\biggl\}. $$ 
Next, we simlify the coefficient $C_{2}$ in a similar manner. Hence, it is given by
\begin{multline*}
C_{2}=\biggl\{\frac{\ln(\frac{x+a}{x-1})}{\ln[\frac{(a+x)(ax+1)}{a(x-1)^{2}}]}\biggl\}=
\biggl\{\frac{\ln(\frac{x+a}{x-1})}{\ln[(\frac{x+a}{x-1})(\frac{x+\frac{1}{a}}{x-1})]}\biggl\} = \\ 
\biggl\{\frac{\ln(\frac{x+a}{x-1})}{\ln(\frac{x+a}{x-1})+\ln(\frac{x+\frac{1}{a}}{x-1})}\biggl\}
=\biggl\{\frac{x_{21}}{x_{21}+x_{22}}\biggl\}=\biggl\{\frac{nx_{22}}{nx_{22}+x_{22}}\biggl\}= 
\biggl\{\frac{n}{n+1}\biggl\}. 
\end{multline*}
 We thus obtain a formula for an inverse gemini function pair, namely 
\begin{eqnarray}\label{65}
\operatorname{Li_{2}}(-a)=-\biggl\{\frac{2n-1}{n+1}\biggl\}\frac{\pi^{2}}{6}-\frac{1}{2}\biggl\{\frac{n}{n+1}\biggl\}\ln^{2}(a)
\end{eqnarray}
 By applying the inversion identity \eqref{inversionformula} to \eqref{65}, we obtain 
\begin{eqnarray}\label{66}
\operatorname{Li_{2}}\left(-\frac{1}{a}\right)=\biggl\{\frac{n-2}{n+1}\biggl\}\frac{\pi^{2}}{6}-\frac{1}{2}\biggl\{\frac{1}{n+1}\biggl\}\ln^{2}(a)
\end{eqnarray}
 We next examine the case whereby $n=2$, and this gives us that 
 $$\operatorname{Li_{2}}(-a)=-\biggl\{\frac{2\cdot 2-1}{2+1}\biggl\}\frac{\pi^{2}}{6}-\frac{1}{2}\biggl\{\frac{2}{2+1}\biggl\}\ln^{2}(a), $$ so 
 that 
 $$\operatorname{Li_{2}}(-a)=-\frac{\pi^{2}}{6}-\frac{1}{3}\ln^{2}(a) \ \ \ \text{or} 
 \ \ \ \operatorname{Li_{2}}\left(-\frac{1}{a} \right) = 
 -\frac{1}{6}\ln^{2}(a).$$ 
 
The functions $\gemini_{\frac{1}{\phi}}(x)$ and $\gemini_{\phi}(x)$ form also an inversion function pair and the exact area
 ratio $n$ can be evaluated
 so that 
\begin{multline*}
 \operatorname{Li_{2}}(-\phi)=-\biggl\{\frac{2n-1}{n+1}\biggl\}\frac{\pi^{2}}{6}-\frac{1}{2}\biggl\{\frac{n}{n+1}\biggl\}\ln^{2}(\phi) \Rightarrow \\ 
 n=\frac{\pi^2 - 6\operatorname{Li_{2}}(-\phi)}{2pi^2 + 
 3\ln^2(\phi)+6\operatorname{Li_{2}}(-\phi)}=\frac{22\pi^{2}}{7\pi^{2}-15\ln^{2}(\phi)}- 2. 
\end{multline*}
\end{example}

 We demonstrate the application of \eqref{64} and \eqref{65}. 
If the total area of the gemini function is $\frac{\pi^{2}}{3}$, then the corresponding shape factor must be such that
$\operatorname{Li_{2}}(-a)=-\frac{\pi^{2}}{6}$ and $\operatorname{Li_{2}}(-\frac{1}{a})=-\frac{1}{2}\ln^{2}(a) \Rightarrow a\approx2.393308$.
The respective area ratio of the $\gemini_{a}(x)$- and $\gemini_{\frac{1}{a}}(x)$-function is 
 such that 
 $$n=\frac{\frac{\pi^{2}}{6}-\operatorname{Li_{2}}(-a)}{\frac{\pi^{2}}{6}-\operatorname{Li_{2}}(-\frac{1}{a})}=
\frac{\frac{\pi^{2}}{6}+\frac{\pi^{2}}{6}}{\frac{\pi^{2}}{6}+\frac{1}{2}\ln^{2}(a)}=
\frac{\frac{\pi^{2}}{3}}{\frac{\pi^{2}}{6}+\frac{1}{2}\ln^{2}(a)}. $$
Next, we insert this area ratio formula $n$ into the \eqref{64}, and we get
 $$\operatorname{Li_{2}}(-a)=-\left\{\frac{2\cdot\frac{\frac{\pi^{2}}{3}}{\frac{\pi^{2}}{6}+
\frac{1}{2}\ln^{2}(a)}-1}{\frac{\frac{\pi^{2}}{3}}{\frac{\pi^{2}}{6}+\frac{1}{2}\ln^{2}(a)}+1}\right\}\frac{\pi^{2}}{6}-
\frac{1}{2}\left\{\frac{\frac{\frac{\pi^{2}}{3}}{\frac{\pi^{2}}{6}+
\frac{1}{2}\ln^{2}(a)}}{\frac{\frac{\pi^{2}}{3}}{\frac{\pi^{2}}{6}+\frac{1}{2}\ln^{2}(a)}+1}\right\}\ln^{2}(a). $$
 We thus find that $\text{Li}_{2}(-a)$ may be expressed as 
\begin{multline*}
-\biggl\{\frac{\pi^{2}-\ln^{2}(a)}{\pi^{2}+\ln^{2}(a)}\biggl\}\frac{\pi^{2}}{6}-\frac{1}{2}\biggl\{\frac{2\pi^{2}}{3\ln^{2} + 
 3\pi^{2}}\biggl\}\ln^{2}(a) = \\ 
\frac{\pi^{2}}{6}\biggl\{\frac{\ln^{2}(a)-\pi^{2}-2\ln^{2}(a)}{\ln^{2}(a)+\pi^{2}}\biggl\}
=-\frac{\pi^{2}}{6}\biggl\{\frac{\ln^{2}(a)+\pi^{2}}{\ln^{2}(a)+\pi^{2}}\biggl\}= 
-\frac{\pi^{2}}{6}. 
\end{multline*}
Let us still examine the $\operatorname{\gemini_{a}}(x)$-function a bit
more detailed, whose total area is $\frac{\pi^{2}}{3}$. The formula for the total area is given by 
 $$A_{tot}=\int_{0}^{\infty}\gemini_{a}(x) \, dx=\frac{\pi^{2}}{6}-\operatorname{Li_{2}}(-a)=\frac{\pi^{2}}{3} \Rightarrow \operatorname{Li_{2}}(-a)=-
\frac{\pi^{2}}{6}.$$ 
 Hence, we can also write
 $$\operatorname{Li_{2}}(-a)=-\frac{\pi^{2}}{6} \Rightarrow -\operatorname{Li_{2}}(-\frac{1}{a})-\frac{\pi^{2}}{6}-\frac{1}{2}\ln^{2}(a)=
-\frac{\pi^{2}}{6} \Rightarrow \operatorname{Li_{2}}(-\frac{1}{a})=-\frac{1}{2}\ln^{2}(a). $$
 This shape 
 factor $a\approx+2.393308$ or more generally the argument value -2.393308 for a dilogarithm is special. If the argument of the
dilogarithm is $-2.393308$ then its exact value consists only the $\pi^{2}$-term. On the other hand, if the argument of the dilogarihm is a
reciprocal of $-2.393308$, i.e., -0.417831 then the exact value consists only the logarithmic term. Is this value -2.393308 the
only one with this property, or are there other respective values with this same property, or can there be an infinite amount of this kind
of values? The property is simply as follows: Let $\operatorname{Li_{2}}(-a)=-n\pi^{2}$ and
 $\operatorname{Li_{2}}(-\frac{1}{a})=-m\ln^{2}(a)$
 in such a way that $m$ and $n$ are arbitrary rational numbers, i.e., $m, n \in \mathbb{Q}$. 

 By setting the right-hand side of the inversion identity equal to zero, then the outcome is given by $-\frac{\pi^{2}}{6} - 
 \frac{1}{2}\ln^{2}(x)=0 \Rightarrow x = e^{\frac{i\pi}{\sqrt{3}}}$.
This is naturally a multi-valued result, because $x$ is periodic. On the other hand, the LHS of the inversion formula must also be zero at
$x=e^{\frac{i\pi}{\sqrt{3}}}$. Hence, the identity is given by

\begin{eqnarray}
\operatorname{Li_{2}}\left(-e^{\frac{i\pi}{\sqrt{3}}}\right)+\operatorname{Li_{2}}\left(-e^{-\frac{i\pi}{\sqrt{3}}}\right)=0.
\end{eqnarray}

\section{Properties related to the scale factor}\label{sectionscale}
 In this section, we present a few special cases related to the operation with a scale factor.

\subsection{The shape factor versus the scale factor}
 We investigate how the scale factor affects the total area when both factors increase at the same rate. The valid domain for the 
 shape factor is such that $a \in[-1,\infty)$. Hence, we have to set the scale factor $\frac{1}{b}$ in such a way that it is $1$ at $a = 
 -1$, and it starts to scale down from now on such that $b=a+2$ for $b\ge1$. Hence, the formula for the total area is given by
\begin{multline*}
 A_{tot}=\int_{0}^{\infty}\gemini_{a}^{\frac{1}{b}}(x) \, dx=\int_{0}^{\infty}\gemini_{a}^{\frac{1}{a+2}}(x) \, dx = \\ 
\int_{0}^{\infty}\frac{1}{(a+2)}\ln\left(\frac{1+ae^{-x(a+2)}}{1-e^{-x(a+2)}}\right) \, 
 dx=\frac{\frac{\pi^{2}}{6}-\operatorname{Li_{2}}(-a)}{(a + 2)^2}. 
\end{multline*} 

\begin{example}
 We can calculate the critical point $a_{c}$, where the above formula reaches its maximum. At this point, the scale factor starts to 
 dominate the total area. The area increases monotonically up to this point and starts to decrease asymptotically from there. Let us 
 first evaluate a general solution. Instead of inserting the number 2 in the denominator, let this value be an arbitrary parameter $p$. Next, 
 we take the derivative of this function for determining the critical shape factor $a$ as a function of the parameter $p$. Hence, we can 
 write $$\frac{d}{da}A_{tot}(a,p)=\frac{d}{da}\frac{\frac{\pi^{2}}{6}-\operatorname{Li_{2}}(-a)}{(a+p)^2}=\frac{6 
 a\operatorname{Li_{2}}(-a) - a\pi^{2}+(3a+3p)\ln(a+1)}{3a(a+p)^{3}}=0, $$ which implies that 
 $$\operatorname{Li_{2}}(-a)=\frac{\pi^{2}}{6}-\frac{a+p}{2a}\ln(a+1).$$
 By inserting $p=2$, we get
 $$\operatorname{Li_{2}}(-a)=\frac{\pi^{2}}{6}-\frac{a+2}{2a}\ln(a+1),$$
 with $ b=a_{c}+2$. 
 This leads us toward the 
 following question. What will the parameter $p$ be if the critical shape factor $a_{c}=0$? In practice, this means that the
maximum of the $A_{tot}(a,p)$-function is also at $a_{c}=0$. Hence, we are dealing with the $\gemini_{0}(x)$-function, since $a=a_{c}=0$. Now, we can write
 $$\lim\limits_{a \to 0} \left[\operatorname{Li_{2}}(-a)-\frac{\pi^{2}}{6}+\frac{a+p}{2a}\ln(a+1)\right]=\frac{p}{2}-\frac{\pi^{2}}{6}=0 \Rightarrow
p=\frac{\pi^{2}}{3}, $$
 and this implies that 
$ b=a_{c}+p=0+\frac{\pi^{2}}{3}=\frac{\pi^{2}}{3}$.

\begin{figure}[htbp!]
\begin{center}
\includegraphics[width=8cm, height=6cm]{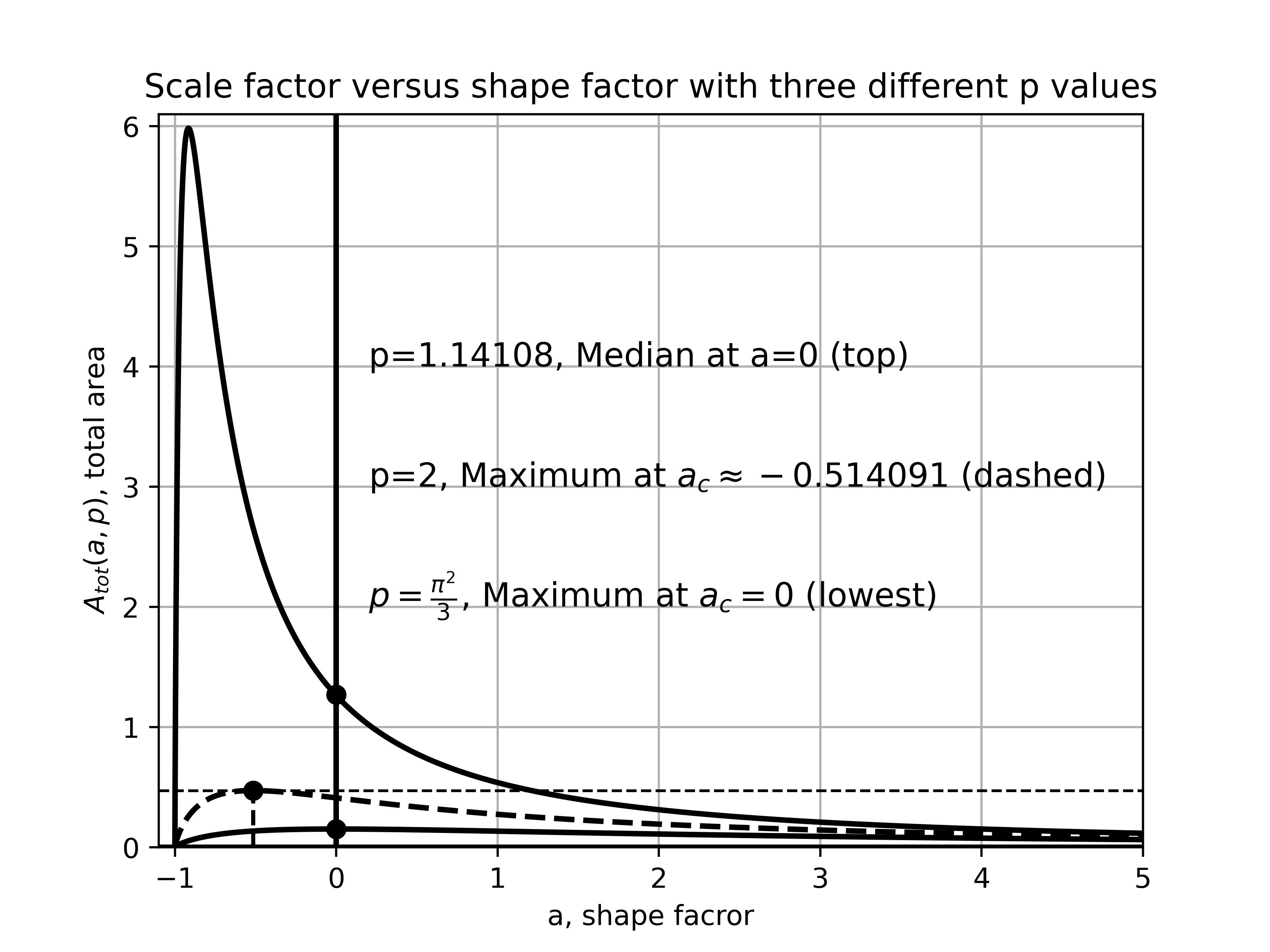}
\end{center}
\caption{\label{Figure13} In this figure, three total area curves are plotted with the different $p$-parameters. The dashed graph
corresponds to the case $p=2$. The scale factor starts to dominate the total area at the critical point $a_{c}\approx-0.514091$. The
top graph has its median at $a=0$ and the lowest graph has its maximum at $a=0$.}
\end{figure}

 When $b=\frac{\pi^{2}}{3}$, 
 the corresponding 
 $A_{tot}(a,p)$-function reaches its maximum at $a_{c}=0$, with a scale factor equal to $b$. The
 corresponding maximum total area of the $\gemini_{0}^{\frac{3}{\pi^{2}}}(x)$-function is given by
\begin{multline*}
 A_{max}(a_{c},\frac{1}{b})=A_{max}(0,\frac{3}{\pi^{2}})=\int_{0}^{\infty}\gemini_{0}^{\frac{1}{b}}(x) \, dx=
\frac{1}{b}\int_{0}^{\infty}\ln\left(\frac{1}{1-e^{-bx}}\right) \, dx = \\ 
 -\frac{1}{b^{2}}\bigg|_{0}^{\infty}\operatorname{Li_{2}}(e^{-bx})=\frac{9}{\pi^{4}}\cdot\frac{\pi^{2}}{6}=\frac{3}{2\pi^2}. 
\end{multline*} 
This obtained $A_{tot}(a,p)$-function is itself interesting, because its total area is also finite. In the above, we evaluated the case, where
the parameter $p$ was fixed equal to 2. By replacing this value 2 to an arbitrary parameter $p$ and integrating the $A_{tot}(a,p)$-function
from $-1$ to infinity, we obtain a new function $A(p)$. The result of this improper integral leads us to 
 define 
 $$ A(p) := \frac{\pi^{2}}{2p}+\frac{\ln^{2}(p-1)}{2p}, $$
 for $p>1$. 
\end{example}

\begin{example}
We can derive a single-term dilogarithm representation using $A(p)$-function. Let us define the parameter $p$ in such a way
that the median of the respective $A_{tot}(a,p)$-function is at $a=0$. Thus, the half of the total area can be formulated in two different ways.
Hence, we can write 
\begin{multline*}
 \frac{1}{2}A_{tot}(a,p)=\int_{-1}^{0}\frac{\frac{\pi^2}{6}-\operatorname{Li_{2}}(-a)}{(a+p)^2}da = \\ 
 -\frac{1}{p}\operatorname{Li_{2}}\left(\frac{1}{1-p}\right) = 
 \frac{1}{p}\left[\operatorname{Li_{2}}(1-p)+\frac{\pi^2}{6}+\frac{1}{2}\ln^{2}(p-1)\right],
\end{multline*} 
so that $ \frac{1}{2}A_{tot}(a,p)$ may be expressed as 
\begin{multline*}
 \frac{1}{p}\left[\operatorname{Li_{2}}\left(\frac{1}{p}\right)+\frac{1}{2}\ln^{2}(p-1)+\frac{1}{2}\ln(p)\ln\left(\frac{p}{(p-1)^{2}}\right)\right] = \\ 
 \frac{1}{2}A(p)=
\frac{1}{2}\left[\frac{\pi^{2}}{2p}+\frac{\ln^{2}(p-1)}{2p}\right], 
\end{multline*}
 so that 
 $$\operatorname{Li_{2}}\left( 
 \frac{1}{p} \right) + 
 \frac{1}{2}\ln(p)\ln\left[\frac{p}{(p-1)^{2}}\right]=\frac{\pi^{2}}{4}-\frac{\ln^{2}(p-1)}{4}, $$
 and this implies that 
 $$\operatorname{Li_{2}}\left(\frac{1}{p}\right)=\frac{\pi^{2}}{4}-\ln^{2}\sqrt{p-1}-\ln(p)\ln\left(\frac{\sqrt{p}}{p-1}\right). $$
\end{example}

\subsection{Fitting total areas with scale factors}

 In this Section, we represent some simple operations with the scale factor $b$. If one wants to fit an arbitrary gemini function so that 
 its total area becomes the same as another gemini function, we have to either scale down or scale up the function with the aid of the 
 scale factor. In this case, the scale factor can be written by $$b=\sqrt{\frac{A_{tot_1}}{A_{tot_2}}}=\sqrt{\frac{\frac{\pi^{2}}{6} - 
 \operatorname{Li_{2}}(-a_{1})}{\frac{\pi^{2}}{6}-\operatorname{Li_{2}}(-a_{2})}}.$$

Let us next compare the total areas of the degenerate $\gemini_{0}(x)$ and fundamental $\gemini_{1}(x)$ forms of a gemini function. The respective
areas are such that $A_{tot_0}=\frac{\pi^2}{6}$ and $A_{tot_1}=\frac{\pi^2}{4}$. To make the areas equal, one can magnify the area of a
$\gemini_{0}(x)$-function, with the scale factor $b$ greater than one. Respectively, one can scale down the $\gemini_{1}(x)$-function 
 with the
scale factor $b$ less then one. It is worth to note once again that the area increases and decreases proportionally to $b^{2}$. Now, the
required scale factor $b$ is given by $b = \frac{3}{2}$, 
 and the corresponding gemini function is such that 
 $$\gemini_{0}^{\sqrt{\frac{3}{2}}}(x)=\sqrt{\frac{3}{2}}\ln\left(\frac{1}{1-e^{-x\sqrt{\frac{2}{3}}}}\right).$$ 
 Consequently, we find that 
 $$\int_{0}^{\infty}\gemini_{0}^{b}(x) \, dx=b\int_{0}^{\infty}\ln\left(\frac{1}{1-e^{-\frac{x}{b}}}\right) \, dx=
-b^{2}\bigg|_{0}^{\infty}\operatorname{Li_{2}}(e^{-\frac{x}{b}})=\frac{3}{2} \cdot \frac{\pi^{2}}{6}=\frac{\pi^{2}}{4}.$$

\begin{figure}[htbp!]
\begin{center}
\includegraphics[width=10cm, height=7cm]{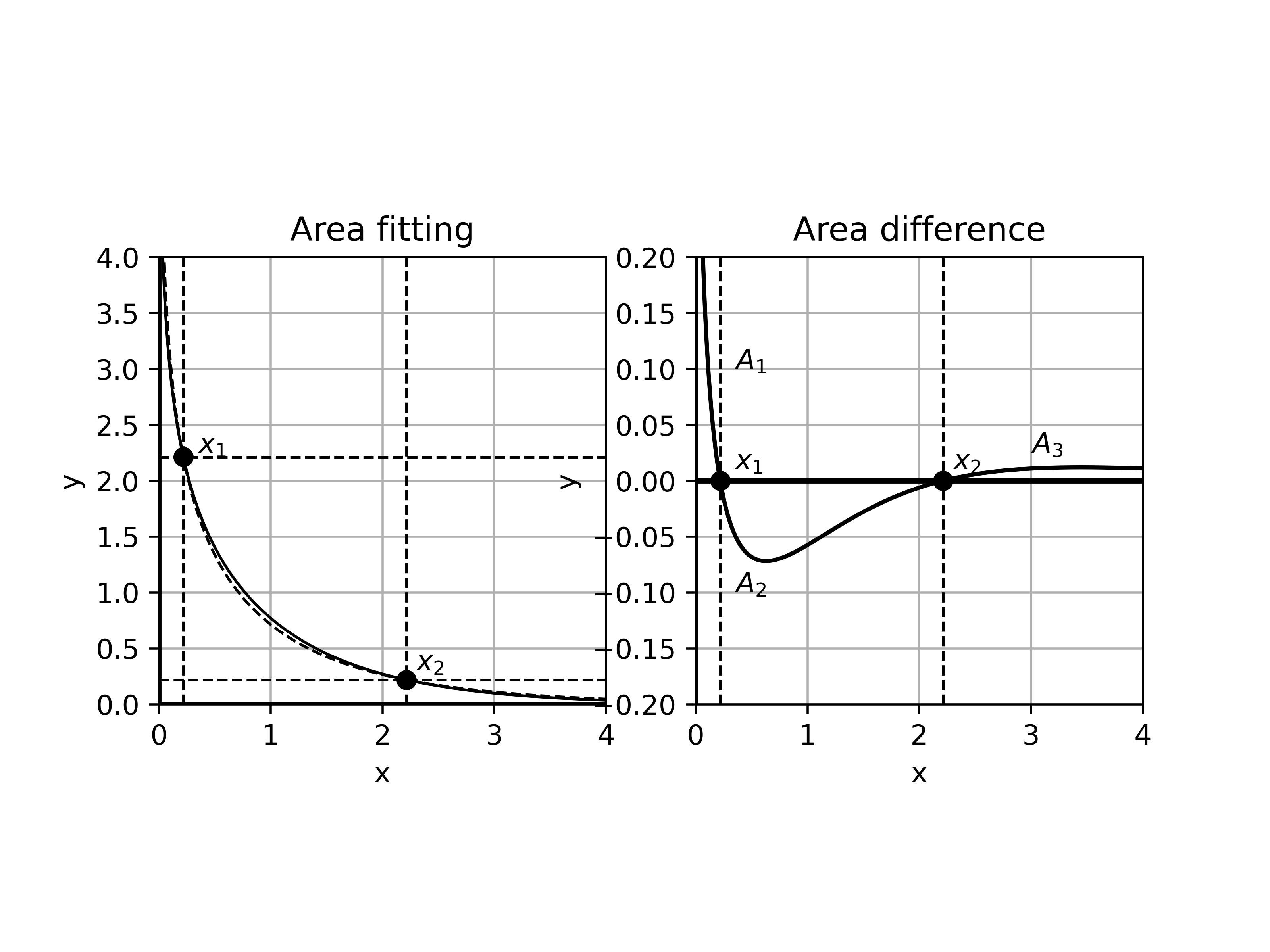}
\end{center}
\caption{\label{Figure14} The left hand side plot depicts the fitted dashed graph of the $\gemini_{0}^{\sqrt{\frac{3}{2}}}(x)$-function and
the solid graph of the $\gemini_{1}(x)$-function. The total areas are equal. The right hand side plot depicts the area difference between
$\gemini_{0}^{\sqrt{\frac{3}{2}}}(x)$- and $\gemini_{1}(x)$-functions. The total area above the $x$-axis is equal to the area
below the $x$-axis, i.e., $A_{1}+A_{3}={\lvert}A_{2}{\rvert}$.}
\end{figure}

Observe that the dashed graph of the fitted function in Figure \ref{Figure14} is always lower, relative to the domain of the intersection
points $x_{1}$ and $x_{2}$, when $b>1$. In a closely related way, the fitted dashed function attains strictly lesser values
elsewhere, compared to the solid target function. The area difference function has an interesting
property. The sum of the two areas above the $x$-axis is equal to the area below the $x$-axis, i.e., $A_{1}+A_{3}=|A_{2}|$ and $A_{1}=A_{3}$.

\section{Further applications}
 Campbell recently introduced a numerically discovered identity \cite{Campbell2025} that is shown below, 
 and the problem of proving this 
result was left as an open problem in the work of Campbell. We succeed in obtaining a complete proof of this experimentally discovered 
and conjectured result. An alternative approach toward proving this result was recently given by Hakimoglu-Brown \cite{Hak_25}. 

\begin{theorem}
The closed form evaluation 
\begin{multline*}
\operatorname{Li_{2}}\left(\frac{1}{2\phi^{2}}-\frac{1}{2}\sqrt{-1-\frac{1}{\phi^{2}}}\right)-
\operatorname{Li_{2}}\left(\frac{1-\sqrt{(1-2\phi)(1+2\phi)}}{2}\right) = \\ 
\frac{\ln^{2}(\phi)}{2}+\frac{3\pi\ln(\phi)i}{5}+\frac{\pi^{2}}{150} 
\end{multline*}
holds true (cf.\ \cite{Campbell2025}). 
\end{theorem}

\textit{Proof.} We begin by manipulating the first dilogarithm term with the aid of a reflection identity, so that 
$$\operatorname{Li_{2}}\left(\frac{1}{2\phi^{2}}-\frac{1}{2}\sqrt{-1-\frac{1}{\phi^{2}}}\right) = 
\operatorname{Li_{2}}\left(\frac{1}{2\phi^{2}}-\frac{i}{2}\sqrt{1+\frac{1}{\phi^{2}}}\right), $$ yielding the equality 
\begin{multline*}
\operatorname{Li_{2}}\left(\frac{1}{2\phi^{2}}-\frac{1}{2}\sqrt{-1-\frac{1}{\phi^{2}}}\right)=-\operatorname{Li_{2}}\left(\frac{\phi}{2}+\frac{i}{2}\sqrt{1+\frac{1}{\phi^{2}}}\right)+\frac{\pi^{2}}{6} - \\ 
\ln\left(\frac{\phi}{2}+\frac{i}{2}\sqrt{1+\frac{1}{\phi^{2}}}\right)\ln\left(\frac{1}{2\phi^{2}}-\frac{i}{2}\sqrt{1+\frac{1}{\phi^{2}}}\right), 
\end{multline*}
and the right-hand side then reduces to $-\operatorname{Li_{2}}\left(e^{\frac{i\pi}{5}}\right)+\frac{13\pi^{2}}{150}+\frac{i\pi\ln(\phi)}{5}$. 
Next, we apply Theorem \ref{mainthreeterm} by setting the argument of the second term to be equal to $\phi$, writing $\frac{1}{a^{2} 
- a + 1}=\phi \Rightarrow a=\frac{1 \pm i\sqrt{5-2\sqrt{5}}}{2}$. By then inserting the root with a positive imaginary part into 
Theorem \ref{mainthreeterm}, we obtain that 
\begin{multline*}
\operatorname{Li_{2}}\left(\frac{a-1}{a^2}\right)=\operatorname{Li_{2}}\left(\frac{1+i\sqrt{5+2\sqrt{5}}}{2}\right) = \\ 
\operatorname{Li_{2}}\left(\frac{1+i\sqrt{\sqrt{5}\phi^{3}}}{2}\right)=\operatorname{Li_{2}}\left(\frac{1+\sqrt{(1-2\phi)(1+2\phi)}}{2}\right)
\end{multline*}
 and that $$\operatorname{Li_{2}}\left(\frac{1}{a^{2}-a+1}\right)=\operatorname{Li_{2}}(\phi),$$ with 
$-\operatorname{Li_{2}}(\frac{a}{a^{2}-a+1})=-\operatorname{Li_{2}}(e^{\frac{i\pi}{5}})$ and $\ln\left(\frac{a}{a - 
1}\right)\ln\left(\frac{a}{a^{2}-a+1}\right)=\frac{6\pi^{2}}{25}+\frac{3i\pi\ln(\phi)}{5}$. Next, we compose the identity Theorem 
\ref{mainthreeterm} with the above evaluated terms. This yields $$\operatorname{Li_{2}}\left(\frac{1+i\sqrt{\sqrt{5}\phi^{3}}}{2}\right) + 
\operatorname{Li_{2}}(\phi)-\operatorname{Li_{2}}(e^{\frac{i\pi}{5}})+\frac{6\pi^{2}}{25}+\frac{3i\pi\ln(\phi)}{5} = \frac{2\pi^{2}}{5}. $$ 
If the argument of a dilogarithm is of the form
$\frac{1}{2} \pm iu$ for $u \in \mathbb{R}$ then its real part can always be determined by the formula 
$$ \mathfrak{Re}\biggl\{\operatorname{Li_{2}}\left(\frac{1}{2}+iu\right)\biggl\}=\frac{\pi^{2}}{12} - 
 \frac{1}{8}\ln^{2}\left(\frac{1+4u^{2}}{4}\right)-
\frac{\operatorname{arctan^{2}}(2u)}{2}, $$
 which implies that 
 $$ \mathfrak{Re} \left\{\operatorname{Li_{2}}\left(\frac{1+i\sqrt{\sqrt{5}\phi^{3}}}{2}\right) \right\} =
 \frac{\pi^{2}}{300}-\frac{\ln^{2}(\phi)}{2}. $$ 
 The exact value for $\operatorname{Li_{2}}(\phi)$ is $\frac{7\pi^{2}}{30}+\frac{\ln^{2}(\phi)}{2}-i\pi\ln(\phi)$. Hence, the real part is 
 given by
 $$\mathfrak{Re}\biggl\{\operatorname{Li_{2}}\left(\phi\right)\biggl\}=\frac{7\pi^{2}}{30}+\frac{\ln^{2}(\phi)}{2}.$$
 The real part for the 
 third term may be evaluated via
 Kummer's rule, which is given by
 $$\mathfrak{Re}\biggl\{\operatorname{Li_{2}}\left(e^{i\theta}\right)\biggl\}=\frac{\pi^{2}}{6}-\frac{2\pi\theta-\theta^{2}}{4} \Rightarrow
-\mathfrak{Re}\biggl\{\operatorname{Li_{2}}\left(e^{\frac{i\pi}{5}}\right)\biggl\}=-\frac{23\pi^{2}}{300}.$$ 
 By adding all the evaluated real parts together, we get
 $$\frac{\pi^{2}}{300}-\frac{\ln^{2}(\phi)}{2}+\frac{7\pi^{2}}{30} + 
 \frac{\ln^{2}(\phi)}{2}-\frac{23\pi^{2}}{300}=+\frac{4\pi^{2}}{25}, $$ 
 as desired. 
 Hence, the left-hand side 
 real part constant is necessarily $\frac{6\pi^{2}}{25}-\frac{2\pi^{2}}{5}=-\frac{4\pi^{2}}{25}$, which is canceled and
makes the identity zero, i.e., $+\frac{4\pi^{2}}{25}-\frac{4\pi^{2}}{25}=0$.
 Next, we proceed by inserting the analytically evaluated real part terms into the obtained identity and we can write
$$\operatorname{Li_{2}}\left(\frac{1+i\sqrt{\sqrt{5}\phi^{3}}}{2}\right)+
\operatorname{Li_{2}}(\phi)-\operatorname{Li_{2}}(e^{\frac{i\pi}{5}})-\frac{4\pi^{2}}{25}+\frac{3i\pi\ln(\phi)}{5}=0, $$
 and this implies that 
$$\operatorname{Li_{2}}\left(\frac{1+i\sqrt{\sqrt{5}\phi^{3}}}{2}\right)-\operatorname{Li_{2}}(e^{\frac{i\pi}{5}})+\frac{11\pi^{2}}{150}+
\frac{\ln^{2}(\phi)}{2}-\frac{2i\pi\ln(\phi)}{5}=0.$$
 Next, the first term of this obtained identity must be transformed by a reflection identity, making it the same as the first term of
Campbell's identity. Consequently, we have that 
\begin{multline*}
 \operatorname{Li_{2}}\left(\frac{1-\sqrt{(1-2\phi)(1+2\phi)}}{2}\right) = \\ 
\operatorname{Li_{2}}\left(1-\left(\frac{1+i\sqrt{\sqrt{5}\phi^{3}}}{2}\right)\right) = 
\operatorname{Li_{2}}\left(\frac{1-i\sqrt{\sqrt{5}\phi^{3}}}{2}\right), 
\end{multline*}
 so that 
\begin{multline*}
 \operatorname{Li_{2}}\left(\frac{1-i\sqrt{\sqrt{5}\phi^{3}}}{2}\right) = 
 -\operatorname{Li_{2}}\left(\frac{1+i\sqrt{\sqrt{5}\phi^{3}}}{2}\right) + \\ 
\frac{\pi^2}{6}-\ln\left(\frac{1-i\sqrt{\sqrt{5}\phi^{3}}}{2}\right)\ln\left(\frac{1+i\sqrt{\sqrt{5}\phi^{3}}}{2}\right), 
\end{multline*}
 so that 
 $$\operatorname{Li_{2}}\left(\frac{1-i\sqrt{\sqrt{5}\phi^{3}}}{2}\right)=-\operatorname{Li_{2}}\left(\frac{1+i\sqrt{\sqrt{5}\phi^{3}}}{2}\right)+
\frac{\pi^2}{6}-\frac{4\pi^2}{25}-\ln^{2}(\phi), $$
 and hence 
\begin{multline*}
 \operatorname{Li_{2}}\left(\frac{1-i\sqrt{\sqrt{5}\phi^{3}}}{2}\right)=-\operatorname{Li_{2}}\left(e^{\frac{i\pi}{5}}\right)+\frac{7\pi^2}{30}
+\frac{\ln^{2}(\phi)}{2}-i\pi\ln(\phi) - \\ \frac{4\pi^2}{25} + \frac{3i\pi\ln(\phi)}{5}+\frac{\pi^2}{6} - 
 \frac{4\pi^2}{25}-\ln^{2}(\phi). 
\end{multline*}
As a consequence, we have that 
\begin{multline*}
-\operatorname{Li_{2}}\left(\frac{1-i\sqrt{\sqrt{5}\phi^{3}}}{2}\right)=-\operatorname{Li_{2}}\left(\frac{1-\sqrt{(1-2\phi)(1+2\phi)}}{2}\right) = \\ 
\operatorname{Li_{2}}\left(e^{\frac{i\pi}{5}}\right)-\frac{2\pi^2}{25}+\frac{\ln^{2}(\phi)}{2}+\frac{2i\pi\ln(\phi)}{5}.
\end{multline*}
 By inserting our equivalent expressions for the dilogarithmic
 terms appearing in Campbell's numerically discovered, 
 this gives us an equivalent version of the desired result. 
 
 In addition, we can determine analytical values for both terms of Campbell's identity, since the exact value of 
 $\operatorname{Li_{2}}\left(e^{\frac{i\pi}{5}}\right)$ can be evaluated in terms of the trigamma function. 

\subsection{Rederiving two known dilogarithmic ladders in the base $\frac{1}{2}$}
 We can apply our addinacci identity to derive two identities introduced in the work of Bailey et al.\ \cite{Bai_97}. Let us set $x = 
 4$ then the respective exponent constant $N$ is given by
 $$x-2=+\frac{1}{x^{N}} \Rightarrow 4-2=\frac{1}{4^{N}} \Rightarrow N=-\frac{\ln(2)}{\ln(4)}=-\frac{1}{2}.$$ 
 We proceed to construct 
 the corresponding five-term addinacci-identity in such a way that $x=4$ and $N=-\frac{1}{2}$, with 
\begin{multline*}
4\operatorname{Li_{2}}\left(\frac{1}{4}\right)-2\operatorname{Li_{2}}\left(\frac{1}{4^{-\frac{3}{2}}}\right) 
 +4\operatorname{Li_{2}}\left(\frac{1}{4^{-\frac{1}{2}}}\right)+\operatorname{Li_{2}}\left(\frac{1}{4^{-3}}\right) - \\ 
 2\operatorname{Li_{2}}\left(\frac{1}{4^{-1}}\right)-\frac{\pi^{2}}{3}+2\ln^{2}(4)= 
 0, 
\end{multline*}
 so that 
 $$4\operatorname{Li_{2}}\left(\frac{1}{4}\right)-2\operatorname{Li_{2}}\left(8\right)+
4\operatorname{Li_{2}}\left(2\right)+\operatorname{Li_{2}}\left(64\right)-
2\operatorname{Li_{2}}\left(4\right)-\frac{\pi^{2}}{3}+8\ln^{2}(2)=0. $$
 Now, we apply \eqref{complex}, which 
 enables us to convert the dilogarithm terms with arguments greater than one into their reciprocals, so that 
$$ 4\operatorname{Li_{2}}\left(\frac{1}{2}\right)-6\operatorname{Li_{2}}\left(\frac{1}{4}\right)-2\operatorname{Li_{2}}\left(\frac{1}{8}\right)
+\operatorname{Li_{2}}\left(\frac{1}{64}\right)=\ln^{2}(2). 
 $$ This identity can be can be converted so that 
$$ 36\operatorname{Li_{2}}\left(\frac{1}{2}\right)-36\operatorname{Li_{2}}\left(\frac{1}{4}\right) - 
 12\operatorname{Li_{2}}\left(\frac{1}{8}\right)
+6\operatorname{Li_{2}}\left(\frac{1}{64}\right)=\pi^{2}. $$

\subsection{A further derivation related to the addinacci-identity}
 Next, we apply Theorem \ref{theoremaddin} to derive a five-term ladder, by setting $N=-\frac{3}{4}$, giving us that $$ x = 1+\sqrt{1 + 
 \frac{1}{x^{-\frac{3}{4}-1}}}.$$ The corresponding ladder is given by
\begin{multline*}
 4\operatorname{Li_{2}}\left(\frac{1}{x}\right)-2\operatorname{Li_{2}}\left(\frac{1}{x^{-\frac{3}{4}-1}}\right)+
4\operatorname{Li_{2}}\left(\frac{1}{x^{-\frac{3}{4}}}\right) + \\ 
 \operatorname{Li_{2}}\left(\frac{1}{x^{-\frac{3}{4}\cdot2-2}}\right)
-2\operatorname{Li_{2}}\left(\frac{1}{x^{-\frac{3}{4}\cdot2}}\right)- 
 \frac{\pi^{2}}{3}-2\ln^{2}(x) = 0, 
\end{multline*}
 and this implies that 
\begin{multline*}
 4\operatorname{Li_{2}}\left(\frac{1}{x}\right)-2\operatorname{Li_{2}}\left(x^{\frac{7}{4}}\right) + 
 4\operatorname{Li_{2}}\left(x^{\frac{3}{4}}\right)+\operatorname{Li_{2}}\left(x^{\frac{7}{2}}\right) - \\ 
 2\operatorname{Li_{2}}\left(x^\frac{3}{2}\right)-\frac{\pi^{2}}{3}-2\ln^{2}(x) = 0.
\end{multline*}
 We then set $y=\sqrt[4]{x}=\sqrt[4]{2\mathcal{T}_{tri}+2}$, and hence 
 $$4\operatorname{Li_{2}}\left(\frac{1}{y^{4}}\right)-2\operatorname{Li_{2}}\left(y^{7}\right)
+4\operatorname{Li_{2}}\left(y^{3}\right)+\operatorname{Li_{2}}\left(y^{14}\right)
-2\operatorname{Li_{2}}\left(y^{6}\right)-\frac{\pi^{2}}{3}-2\ln^{2}(y^{4}) = 0.$$
 We proveed to again apply \eqref{complex}, 
 to convert the dilogarithm terms with arguments greater than one into their reciprocals.
Hence, we get
\begin{eqnarray}\label{triaddi1}
4\operatorname{Li_{2}}\left(\frac{1}{y^{3}}\right)-4\operatorname{Li_{2}}\left(\frac{1}{y^{4}}\right)
-2\operatorname{Li_{2}}\left(\frac{1}{y^{6}}\right)
-2\operatorname{Li_{2}}\left(\frac{1}{y^{7}}\right)+\operatorname{Li_{2}}\left(\frac{1}{y^{14}}\right)-\ln^{2}(y)=0. \quad \quad 
\end{eqnarray}

By substituting $y=\frac{\mathcal{T}_{tri}+1}{\mathcal{T}_{tri}}$ into \eqref{triaddi1}, the same identity becomes, as shown
in \eqref{triaddi2}.

\begin{eqnarray}
4\operatorname{Li_{2}}\left(\frac{1}{2\mathcal{T}_{tri}}\right)-4\operatorname{Li_{2}}\left(\frac{1}{2\mathcal{T}_{tri}+2}\right)
-2\operatorname{Li_{2}}\left(\frac{1}{4\mathcal{T}_{tri}^{2}}\right)
-2\operatorname{Li_{2}}\left(\frac{2-\mathcal{T}_{tri}}{4\mathcal{T}_{tri}-4}\right)+
\nonumber
\end{eqnarray}
\begin{eqnarray}\label{triaddi2}
 \operatorname{Li_{2}}\left(\frac{5\mathcal{T}_{tri}-9}{64\mathcal{T}_{tri}-32}\right)-
 \ln^{2}\left(\frac{\mathcal{T}_{tri}+1}{\mathcal{T}_{tri}}\right)=0.
\end{eqnarray}

\section{Conclusion}
 The methods we have applied raise questions as to how such methods or similar methods could be applied to obtain new proofs of 
 Watson's dilogarithmic identities \cite{Wat_29}. We are convinced that they can be derived using gemini-identities, as long as 
 we find the right initial values. We leave it to a future project to explore this.

The primary result in this publication is the five-term gemini identity, which can be used to derive many new dilogarithm evaluations.
Ideally, our results will lead to new areas of research concerning higher logarithm functions as well as number-theoretic applications
of such functions. 
 
\subsection*{Acknowledgments}
 The experimental discovery of our results relied heavily on the use of the Wolfram language. It has been an invaluable help, and without it, 
 the author would not have been able to perform the calculations related to this publication. I would also like to express special thanks to 
 Dr.\ Lauri Jetsu for his editorial help and encouragement. Otherwise, writing this publication has been a lonely journey through the forest 
 of dilogarithm functions. The author thanks John Maxwell Campbell for help with the preparation of this manuscript and for
 useful feedback.

\end{document}